\input ./style/arxiv-general.cfg
\documentclass[aap,MSNbibl,dvips]{arximspdf}
\makeatletter
   \@ifpackageloaded{graphicx}{}{\usepackage{graphicx}}
\makeatother

\doi{10.1214/14-AAP1052} 
\volume{25}
\issue{5}
\pubyear{2015}
\firstpage{2416}
\lastpage{2461}
\docsubty{FLA}

\makeatletter
\renewcommand{\mid}{|}
\newcommand{\rrVert}{\Vert}
\newcommand{\llVert}{\Vert}

\newtheorem{theorem}{Theorem}[section]
\newtheorem{lemma}[theorem]{Lemma}
\newtheorem{proposition}[theorem]{Proposition}
\newtheorem{conjecture}{Conjecture}
\newproclaim{remark}{Remark}[section]
\newcommand{\critical}{\mathrm{c}}
\newcommand{\detpath}{\mathrm{t}}
\newcommand{\compl}{\mathrm{c}}
\newcommand{\veci}{i}
\newcommand{\origin}{0}
\makeatother

\begin{document}
\begin{frontmatter}

\title{Space--time percolation and detection by mobile~nodes}
\runtitle{Space--time percolation and detection by mobile nodes}

\begin{aug}
\author{\fnms{Alexandre}~\snm{Stauffer}\corref{}\ead[label=e1]{a.stauffer@bath.ac.uk}\thanksref{T1}}
\runauthor{A. Stauffer}
\affiliation{University of Bath}
\address{Department of Mathematical Sciences\\
University of Bath\\
Claverton Down\\
Bath, BA1 2BG\\
United Kingdom\\
\printead{e1}} 
\end{aug}
\thankstext{T1}{This work was done while the author was at Microsoft
Research, Redmond WA, USA.}

\received{\smonth{3} \syear{2012}}
\revised{\smonth{7} \syear{2014}}

%
\begin{abstract}
Consider the model where nodes are initially distributed as a Poisson
point process with intensity $\lambda$ over $\mathbb{R}^d$
and are moving in continuous time according to independent Brownian motions.
We assume that nodes are capable of detecting all points within
distance $r$ of their location and
study
the problem of determining the first time at which a target particle,
which is initially placed at the origin of~$\mathbb{R}^d$, is detected
by at least one node.
We consider the case where the target particle can move according to
any continuous function and can adapt its motion based on the location
of the nodes.
We show that there exists a sufficiently large value of~$\lambda$ so that
the target will eventually be detected almost surely.
This means that the target cannot evade detection even if it has full
information about the past, present and future locations of the nodes.
Also, this establishes a phase transition
for $\lambda$ since, for small enough~$\lambda$, with positive
probability the target can avoid detection forever.
A~key ingredient of our proof is to use fractal percolation and
multi-scale analysis to show that cells with a small density of nodes
do not percolate in space and time.
\end{abstract}

%
\end{frontmatter}

\section{Introduction}\label{secintro}
Let $\Pi_0$ be a Poisson point process over $\mathbb{R}^d$ of
intensity \mbox{$\lambda>0$}.
We refer to the points of~$\Pi_0$ as \emph{nodes}, and let each node
of~$\Pi_0$ move as an independent Brownian motion.
Define $\Pi_s$ to be the point process obtained after the nodes
of~$\Pi
_0$ have moved for time $s$. More formally,
for each $x\in\Pi_0$, let $(\zeta_x(s))_{s\geq0}$ be a standard
Brownian motion and define $\Pi_s = \{x+\zeta_x(s) \dvtx x \in\Pi
_0 \}$.
It is well known that
Brownian motion is a measure-preserving transformation of Poisson point
processes~\cite{vandenberg}, Proposition~1.3, which gives that, for any
fixed $s$,
$\Pi_s$ is also distributed as a Poisson point process. However, $\Pi
_s$ and $\Pi_0$ are \emph{not}
independent, and it is this feature that makes this model challenging
to analyze.

We consider a fixed constant $r$ so that,
at any time, a node is able to detect all points inside the ball of
radius~$r$ centered at its location. Then, letting $B(x,r)$ stand for
the ball of
radius $r$ centered at $x$, we have that
%
%
\begin{equation}
\mbox{at time }s, \mbox{ the nodes of }\Pi_s \mbox{ detect the
region }\bigcup_{x\in\Pi_s} B(x,r). \label{eqdetectregion}
\end{equation}
The region in~(\ref{eqdetectregion}) is related to the random
geometric graph model~\cite{MR,Penrose}, where (nonmobile) nodes are
given by a Poisson point process and edges are added between
pairs of nodes whose distance is at most~$r$.

This model is closely related to a model of mobile graphs introduced by
van den Berg, Meester and White~\cite{vandenberg}. This and other
similar models of mobile graphs
have been considered as natural models for mobile wireless
networks~\cite{elgamal,DiazFinal,clementi,PPPU11,LLMSW12}.

\subsection*{Detection time}
Consider an additional target particle $u$ that is located at the
origin of~$\mathbb{R}^d$ at time $0$. Let $u$ move according to a
continuous function $g(s)$, and define the \emph{detection time}
$T_\mathrm{det}
$ as the first time at which a node is within distance
$r$ from $u$. More formally, we have
\[
T_\mathrm{det}= \inf \biggl\{t\geq0 \dvtx g(t) \in\bigcup
_{x\in\Pi
_t} B(x,r) \biggr\}. %
\]
Detection is a fundamental problem in wireless networks and it appears,
for example, in the contexts of area surveillance and disaster recovery,
where mobile sensors are randomly deployed to explore a region which,
due to some natural disaster, is unsafe for humans~\cite{YG10}.

We consider the case of a target $u$ that wants to evade detection and
can adapt its motion according to the position of the nodes of~$\Pi_0$.
A fundamental question is whether there exists a phase transition on
$\lambda$ so that, for sufficiently large $\lambda$,
the target has no way to avoid detection almost surely as $t\to\infty$.
More formally, define a \emph{trajectory} $h\dvtx\mathbb{R}_+\to
\mathbb
{R}^d$ as a continuous function such that $h(0)$ is the origin
of~$\mathbb{R}^d$.
Then we say that $h$ is \emph{not detected} from time $0$ to $t$ if,
for each $s\in[0,t]$,
all nodes of~$\Pi_s$ are at distance larger than $r$ from $h(s)$.
The existence of such $h$ implies that, if the target chooses its
location at time $s$ to be $g(s)=h(s)$, then it avoids detection up to
time $t$.
We then define
\[
\rho_t(\lambda) = \mathbf{P} \bigl(\exists\mbox{ a trajectory that
is not detected by }(\Pi_s)_{s}\mbox{ from time }0\mbox{
to }t \bigr) %
\]
and, since $\rho_t(\lambda)$ is nonincreasing with $t$ and
nonnegative, the limit exists and we let
\[
\rho= \rho(\lambda) = \lim_{t\to\infty} \rho_t.
\]

Let $V_t=\mathbb{R}^d\setminus\bigcup_{x\in\Pi_t} B(x,r)$ be the
subset of~$\mathbb{R}^d$ that is not detected by the nodes of~$\Pi
_t$, which is usually referred to as the \emph{vacant} region.
Well-known results on random geometric graphs and mobile graphs~\cite
{MR,vandenberg} give that there exists a critical value $\lambda
_\critical$ so that, if $\lambda< \lambda_\critical$, then $V_t$
contains an infinite
connected component (also called the infinite vacant cluster) at all times.
Therefore,
since time is continuous and the target is allowed to move with
arbitrary speed, if $\lambda<\lambda_\critical$, the target can avoid
detection if the origin belongs to the infinite component of~$V_0$,
an event that occurs with positive probability when $\lambda>\lambda
_\critical$.
This gives that $\rho(\lambda)>0$ for all $\lambda<\lambda
_\critical$.
Our Theorem~\ref{thmdetection} below establishes that, if $\lambda$ is
larger than some value $\lambda_\detpath$, where $\detpath$ stands for
trajectory, then $\rho=0$, which means that the target cannot
avoid detection almost surely even if it is able to foresee the
locations of all nodes at \emph{all} times (including future times).
Since $\rho$ is monotone in $\lambda$,
this gives a phase transition on the value of~$\lambda$ for the
existence of a trajectory that is
not detected by the nodes.
%

\begin{theorem}\label{thmdetection}
In dimensions $d\geq2$, there exists a value $\lambda_\detpath
=\lambda
_\detpath(d)\in[\lambda_\critical,\infty)$ such that $\rho=0$ for all
$\lambda>\lambda_\detpath$, and $\rho>0$ for all $\lambda<\lambda
_\detpath$.
Furthermore, there exist $\lambda_\detpath'\geq\lambda_\detpath$, an
explicit positive constant $c=c(d)$ and a positive $C$ independent
of~$t$ such that,
for all large enough $t$ and all $\lambda>\lambda_\detpath'$,
%
%
\begin{equation}
\rho_t(\lambda) \leq\cases{ \displaystyle\exp \biggl(-C
\frac{t}{(\log t)^{c}} \biggr), &\quad for $d=2$,
\cr
\noalign{\vspace*{3pt}}
\displaystyle\exp(-Ct ), &\quad for $d\geq3$.} \label{eqtaildetection}
\end{equation}
\end{theorem}

%
\begin{remark}\label{remdetection}
(i)~For $d=1$, we have that $\rho=0$ for all $\lambda>0$. This holds
since, for any two nodes $v_1,v_2\in\Pi_0$, after a finite time, $v_1$
and $v_2$ will meet
almost surely. Therefore, by considering $v_1$ and $v_2$ such that $u$
is between them at time $0$, we have that $u$ will be detected in
finite time almost surely.

(ii) Note that $\rho<1$ for all $\lambda>0$ since, with constant
probability, the origin of~$\mathbb{R}^d$ is detected at
time~$0$.

(iii)~Based on results by Kesidis, Konstantopoulos and Phoha~\cite
{Kesidis,Konst} (see also the discussion in~\cite{PSSS11}), we have the
following
lower bound for $\rho_t(\lambda)$, which is obtained by considering a
nonmobile target:
\[
\rho_t(\lambda) \geq\cases{ \displaystyle\exp \bigl(-C'
\sqrt{t} \bigr), &\quad for $d=1$,
\cr
\noalign{\vspace*{3pt}} \displaystyle\exp
\biggl(-C' \frac{t}{\log t} \biggr), &\quad for $d=2$,
\cr
\noalign{ \vspace*{3pt}} \displaystyle\exp \bigl(-C' t \bigr), &
\quad for $d \geq3$,} %
\]
for all large enough $t$ and where $C'>0$ does not depend on $t$.
Therefore, comparing this lower bound with the tail bound in~(\ref
{eqtaildetection}) for $\lambda>\lambda_\detpath'$ reveals that,
disregarding logarithmic factors for $d=2$ and
constant factors for $d\geq3$, a target that is able to choose its
motion strategically in response to the
past, present and even future positions of the nodes and is also
capable of moving with arbitrary speed cannot
do much better in terms of avoiding detection than a target that does
not move at all.
\end{remark}

We believe that there
exists a regime for $\lambda$ so that $V_t$ contains no infinite
component at every $t$, but the target is still able to avoid
detection; that is, $\lambda_\detpath>\lambda_\critical$.
If this is true, when $\lambda\in(\lambda_\critical,\lambda
_\detpath
)$, we have that the target is completely surrounded by nodes at all
times; that is, at any time,
the positions at which the target can be at that time form a bounded set.
However, this set may not collapse to the empty set as time proceeds,
and thus the target may be able to avoid detection with positive probability.
Also, we believe that $\lambda_\detpath=\lambda_\detpath'$; that
is, as
soon as $\lambda$ is large enough so that $\rho=0$, then the
exponential tail bounds in~(\ref{eqtaildetection}) should hold.
We formalize these ideas in the following conjecture.
%

\begin{conjecture}
For all dimensions $d\geq2$, we have $\lambda_\critical< \lambda
_\detpath$ and $\lambda_\detpath=\lambda_\detpath'$.
\end{conjecture}

\subsection*{Space--time percolation}
In order to prove Theorem~\ref{thmdetection}, which is the main
motivation of this paper, we develop a framework
that gives a more general result regarding space--time percolation of
increasing events.
We consider a tessellation of~$\mathbb{R}^d$ into cubes of side length
$\ell$,
and a tessellation of the time interval $[0,t]$ into subintervals of
length $\beta$.
Let $\veci\in\mathbb{Z}^d$ be
an index for the cubes of the tessellation
and $\tau\in\mathbb{Z}_+$ be an index for the time intervals.
Let $\varepsilon\in(0,1)$ be fixed and let $E(\veci,\tau)$ be the
indicator random variable
for the event that $\Pi_{\tau\beta}$ contains
at least $(1-\varepsilon)\lambda\ell^d$ nodes in the cube indexed by
$\veci$.

The pairs $(\veci,\tau)$ can be seen as a tessellation of the \emph
{space--time} region $\mathbb{R}^{d+1}$ and
$E(\veci,\tau)$ is a process over this region. We refer to the
space--time subregion indexed by $(\veci,\tau)$ as a \emph{cell}.
For this process, we say that a cell $(\veci,\tau)$ is \emph{adjacent}
to a cell $(\veci',\tau')$ if
$\llVert \veci-\veci'\rrVert _\infty\leq1$ and $\vert \tau-\tau
'\vert \leq1$.
We say that a cell $(i,\tau)$ is
\emph{bad} if $E(\veci,\tau)=0$ and, in this case, define $K(\veci
,\tau
)$ as the set of bad cells from which there exists a path of adjacent
bad cells to
$(\veci,\tau)$; if $E(\veci,\tau)=1$ we let $K(\veci,\tau
)=\varnothing$.
$K(\veci,\tau)$ is usually referred to as the \emph{bad cluster}
of~$(\veci,\tau)$.
Our Theorem~\ref{thmmainresult} below establishes an upper bound for
the bad cluster of the origin $(\origin,0)$, which implies that the bad
cells do not percolate in space and time.
We note that the theorem below is carried out for all dimensions $d\geq
1$, unlike Theorem~\ref{thmdetection} for which the case $d=1$ is
trivial [see Remark~\ref{remdetection}(i)].
%

\begin{theorem}\label{thmmainresult}
Let $\varepsilon>0$, $\ell>0$ and $\beta>0$ be fixed.
Let $Q_t$ be the space--time region $(-t,t)^d\times[0,t/2)$.
Then there exist a positive constant $c=c(d)$, a positive value $C$,
and values $\lambda_0,t_0>0$ such that,
for all $\lambda>\lambda_0$ and $t>t_0$, we have
\begin{eqnarray*}
&&\mathbf{P} \bigl(\exists\mbox{ a cell of }K(\origin,0)\mbox{ not contained in
}Q_t \bigr)\\
&&\qquad \leq\cases{ \displaystyle\exp \biggl(-C \frac{\sqrt{t}}{(\log t)^{c}}
\biggr), &\quad for $d=1$,
\cr
\noalign{\vspace*{3pt}} \displaystyle\exp
\biggl(-C \frac{t}{(\log t)^{c}} \biggr), &\quad for $d=2$,
\cr
\noalign{\vspace*{3pt}}
\displaystyle\exp(-Ct ), &\quad for $d\geq 3$.} %
\end{eqnarray*}
\end{theorem}

In Section~\ref{secfp}, we prove a more general form of Theorem~\ref
{thmmainresult} (which we state in Theorem~\ref{thmfp}).
In this general form,
$E$ is not restricted to be the event defined above, but
can be taken to be the indicator random variable of
\emph{any increasing} event that depends only on a \emph{bounded}
neighborhood of cells and whose marginal probability is large enough.
Under these assumptions, we prove that the cells $(\veci,\tau)$ for
which $E(\veci,\tau)=1$ percolate in space and time, which implies
that the
cells for which $E(\veci,\tau)=0$ do not percolate. We prove
Theorem~\ref{thmdetection} from this result by defining $E$ with
respect to a specific event that
implies that if the target enters a cell $(\veci,\tau)$ for which
$E(\veci,\tau)=1$ then the target is detected.
Since the statement of Theorem~\ref{thmfp} requires some extra
notation, we defer it to Section~\ref{secfp}.

\subsection*{Proof overview}
The proof of Theorem~\ref{thmmainresult} proceeds via a multi-scale
argument inspired by fractal percolation.
We will consider tessellations of the space--time region $\mathbb
{R}^{d+1}$ into cells of various scale. We start with
very large cells so that, by standard arguments, we can show that, with
sufficiently large probability, all such cells contain
a sufficiently high density of nodes at the beginning of their time interval.
Then we consider only the cells that were seen to be sufficiently dense and
partition them into smaller cells; the cells that were observed not to
be dense are simply disregarded, similarly to a fractal percolation process.
We use that the large cells were dense to infer that, with sufficiently
large probability,
the smaller cells also contain a high density of nodes.
We then repeat this procedure until we obtain cells of side length
$\ell
$, for which the density requirement translates
to the event $E(\veci,\tau)=1$.
We show that, despite the procedure of removing nondense cells
described above, the cells of side length $\ell$ that are observed to
be dense
percolate in space and time. The proof of this result requires a
delicate construction that allows us to control dependencies among
cells of various scales.
The details are given in Section~\ref{secfp}.

With this framework of space--time percolation, the proof of
Theorem~\ref{thmdetection} follows rather easily from Theorem~\ref
{thmmainresult}.
The intuition is that, whenever
$E(\veci,\tau)=1$, the cube $\veci$ contains sufficiently many nodes
that can prevent
the target to cut through the cube $\veci$ during the interval $[\tau
\beta,(\tau+1)\beta]$.
Indeed, by setting the parameters $\ell$ and $\beta$ small enough, we
can ensure that the target can only avoid detection if it never enters
a cell for which $E(\veci,\tau)=1$. However, since
the cells for which $E(\veci,\tau)=0$ do not percolate in space and
time by Theorem~\ref{thmmainresult}, we are assured that the target
must be detected. The details are given in Section~\ref{secdetection}.

\subsection*{Related work}
The detection of a target that moves independently of the nodes
of~$(\Pi
_s)_s$ is by now well understood.
Kesidis, Konstantopoulous and Phoha~\cite{Kesidis,Konst} observed that,
for the case $g\equiv0$ (i.e., $u$ does not move),
a very precise asymptotic expression for $\mathbf{P} (T_\mathrm
{det}\geq t )$ can be
obtained using ideas from stochastic geometry~\cite{SKM95}.
This was later extended by Peres et~al.~\cite{PSSS11} and Peres and
Sousi~\cite{PeresSousi11} for
the case when $g$ is any function independent of the nodes of~$(\Pi_s)_s$.
They establish~the interesting fact that
the best strategy for a target that moves independently of
the nodes and wants to avoid detection is to stay put and not to move.
A result of similar flavor was proved for
continuous-time random walks on the square lattice by Drewitz
et~al.~\cite{DGRS}
and Moreau et~al. \cite{MOBC}. The detection problem has also being
analyzed in different models of static and mobile networks~\cite
{Liu,Dousse,Bal}.

Multi-scale arguments as developed in our proof of Theorem~\ref
{thmmainresult} are not uncommon.
The most related references are the papers by Kesten and
Sidoravicius \mbox{\cite{KS04,KS05}} for the study of the spread of infection
in a moving population,
the work of Peres et~al.~\cite{PSSS11} for the so-called percolation
time of mobile geometric graphs, and
a paper by Chatterjee et~al.~\cite{CPPR10} for the gravitational
allocation of Poisson point processes.
However, the techniques in these papers are tailored to specific
problems.
Our Theorem~\ref{thmmainresult} (or, more precisely, the detailed
version of Theorem~\ref{thmfp})
provides a flexible multi-scale analysis that establish space--time
percolation of increasing events. We believe this technique can be
extended to handle other problems as well.

\subsection*{Organization of the paper}
The remainder of the paper is organized as follows.
In Section~\ref{secmixing}, we generalize a coupling argument
developed in~\cite{SinclairStauffer10,PSSS11} for the mixing of moving nodes,
which we will need in our main proofs. Then, in Section~\ref{secfp},
we give the precise statement and
proof of Theorem~\ref{thmmainresult}.
Finally, we apply this result in Section~\ref{secdetection}, where we
prove Theorem~\ref{thmdetection}.

\section{Mixing of mobile graphs}\label{secmixing}
In this section, we extend a coupling argument developed by Sinclair
and Stauffer~\cite{SinclairStauffer10}
(see also~\cite{PSSS11}).
In a high-level description, the result in~\cite{SinclairStauffer10}
establishes that, if a point process contains sufficiently many nodes
in each cell of a suitable tessellation, then
after the nodes have moved as \emph{independent Brownian motions} for
some time interval $\varrho$ that depends on the size of the cells of
the tessellation, the nodes will contain an
independent Poisson point process with
high probability. This argument gives a way to handle dependencies on
mobile geometric graphs, but is not enough in our proof of Theorem~\ref
{thmmainresult}. The reason is that,
when nodes are moving as independent Brownian motions, it may be the
case that a node moves atypically far away during the time
interval~$\varrho$. Therefore, a node
may affect a large region in space during $\varrho$, causing large
dependencies among the cells of the tessellation.
In order to better control dependencies, we only consider nodes that do
not move very far away during
the interval $\varrho$. These nodes are then moving as independent
Brownian motions \emph{conditioned} on not moving very far away during
$\varrho$.
We carry out a more careful analysis of the coupling argument in~\cite
{SinclairStauffer10} in order to derive a corresponding result to this
more general setting for the motion of the nodes.

We note that, very recently, Benjamini and Stauffer~\cite
{BenjaminiStauffer11} employed similar techniques to show
that, for a particular point process,
after the nodes have moved for some time, they will
be \emph{contained} in an independent Poisson point process.

Now we describe the setting.
Fix $\ell>0$ and consider the cube $Q_\ell=[-\ell/2,\ell/2]^d$.
Consider a node that, at time $0$, is located at some arbitrary
position \mbox{$x\in Q_\ell$}. We assume that the position of this
node at
time~$\Delta$ is
distributed according to a translation-invariant function $f_\Delta$,
which means that the probability density function for this node to move from
$x$ to some arbitrary position $y\in\mathbb{R}^d$ after time~$\Delta$
is $f_\Delta(y-x)$.
Let $M>\ell$.
We say that a subdensity function $g\dvtx\mathbb{R}^d \to\mathbb
{R}_+$ is $(\varepsilon,\ell)$-\emph{indistinguishable from
$f_\Delta$
over the cube
$Q_M$} if
%
%
\begin{eqnarray}
\label{eqindistinguishable} g(y) \leq f_\Delta(y-x)
\nonumber
\\[-8pt]
\\[-8pt]
\eqntext{\displaystyle\mbox{for all }x \in Q_\ell\mbox{ and } y\in
Q_M,\mbox{ and also } \int_{Q_{M-\ell}}g(y) \,dy \geq1-
\varepsilon.}
\end{eqnarray}
[Recall that a function $h$ is called a \emph{subdensity} function if
$h(x)\geq0$ for all $x\in\mathbb{R}^d$ and
$\int_{\mathbb{R}^d}h(x) \,dx\leq1$.]
We will apply~(\ref{eqindistinguishable}) in the following context.
Consider two nodes $u$ and $v$ that are initially in
arbitrary positions inside $Q_{\ell/2}$; so their distance vector is in
$Q_\ell$ and is represented by $x$ in~(\ref{eqindistinguishable}).
Assume that, at time $0$, $u$ is at position $a_0$ and that, at time
$\Delta$, the position of~$u$ is $a_\Delta$, where $a_\Delta-a_0$ is
distributed
according to the density function $f_\Delta$. Similarly, assume that
$v$ is at position $b_0$ at time $0$ and that $b_\Delta$ is $v$'s
position at time $\Delta$, where
$b_\Delta-b_0$ is distributed according to the density function given
by a renormalization of~$g$.
Suppose we are given $a_0$ and $b_0$, and we want to sample $b_\Delta$
and $a_\Delta$ in a coupled way.
Then, for any $b_\Delta$ such that
$b_\Delta-b_0\in Q_M$, we have from~(\ref{eqindistinguishable}) that
$g(b_\Delta-b_0) \leq f_\Delta(b_\Delta-a_0)$. Thus, we can obtain a
coupling so that
$a_\Delta=b_\Delta$ whenever $b_\Delta-b_0\in Q_M$. Also, since the
integral of~$g$ over $Q_{M-\ell}$ is at least $1-\varepsilon$, we obtain
that this coupling
gives $\mathbf{P} (a_\Delta=b_\Delta
)\geq1-\varepsilon$. In words, an
indistinguishable function $g$ is such that
a motion according to $g$ inside $Q_M$ is very similar to a motion
according to $f_\Delta$ under any perturbation of the starting point
that is within $Q_\ell$.

We now state Proposition~\ref{procoupling}, which is the most general
version of our extension to~\cite{SinclairStauffer10}, Proposition~4.1.
Later, in Proposition~\ref{procouplingapplied},
we give a special case of Proposition~\ref{procoupling}, which we will
use in our proofs for space--time percolation.
In~\cite{SinclairStauffer10}, a proof is carried out by constructing a
$(\varepsilon,\ell)$-indistinguishable function for the special case when
$f_\Delta$
is a (standard) Brownian motion run for time $\Delta$.
Below, for any two sets $A,A'\subset\mathbb{R}^d$, we define $A+A'$ to
be the set $\{x+x' \dvtx\mbox{for all } x\in A \mbox{ and } x'\in
A'\}$.
%

\begin{proposition}\label{procoupling}
Let $S \supset S'$ be two bounded regions of~$\mathbb{R}^d$ and define
$R = \sup\{k>0 \dvtx S'+Q_k \subseteq S\}$.
Consider any partition of~$S$ into sets $S_1,S_2,\ldots,S_m$ that we
call cells such that the diameter of each cell is at most $\gamma$ for
some fixed $\gamma>0$.
Let $\Phi_0$ be an arbitrary point process at time $0$
such that, for each $i=1,2,\ldots,m$, $\Phi_0$ contains at least
$\beta\operatorname{vol} (S_i )$ nodes in $S_i$ for some
$\beta>0$.
Let $\Phi_\Delta$ be the point process obtained at time $\Delta$ from
$\Phi_0$ after the nodes
have moved independently according to a translation-invariant density
function $f_\Delta$.
Fix $\varepsilon\in(0,1)$.
If there exists a translation-invariant subdensity function $g\dvtx
\mathbb{R}^d\to\mathbb{R}_+$ that is
$(\varepsilon/2,2\gamma)$-indistinguishable from $f_\Delta$ over the cube
$Q_{R+2\gamma}$, then
we can couple the nodes of~$\Phi_{\Delta}$ with those of a Poisson
point process $\Xi$ that is
independent of~$\Phi_0$ and has intensity $(1-\varepsilon)\beta$
so that the nodes of~$\Xi$ are a subset of the nodes of~$\Phi_\Delta$
inside $S'$ with probability at least
\[
1-\sum_{i=1}^m\exp \bigl(-c
\varepsilon^2\beta\operatorname{vol} (S_i ) \bigr)\qquad
\mbox{for some positive constant }c=c(d). %
\]
\end{proposition}

\begin{pf}
We will construct $\Xi$ via three Poisson point processes: $\Xi_0$,
$\Xi_0'$ and~$\Xi_\Delta'$.
We start by defining $\Xi_0$ as a Poisson point process over $S$
with intensity $(1-\varepsilon/2)\beta$.
Then, for each $i=1,2,\ldots,m$,
$\Xi_0$ has fewer nodes than $\Phi_0$ in $S_i$ if
$\Xi_0$ has less than $\beta\operatorname{vol} (S_i )$
nodes in that cell, which by a
standard Chernoff bound (cf. Lemma~\ref{lemcbpoisson}) occurs
with probability larger than
$1-\exp(-\frac{{\delta}^2(1-\varepsilon/2)\beta\operatorname
{vol} (S_i )}{2}(1-{\delta}/3) )$
for $\delta$ such that $(1+\delta)(1-\varepsilon/2)=1$. Note that
$\delta
\in(\varepsilon/2,1)$, so
the probability above can be bounded below by
$1-\exp(-c_1\varepsilon^2\beta\operatorname{vol} (S_i
) )$ for some universal
constant $c_1$.
Let $\{\Xi_0 \preceq\Phi_0\}$ be the event that, for every
$i=1,2,\ldots,m$, $\Xi_0$ has fewer nodes
than $\Phi_0$ in every cell.
Using the union bound we obtain
%
%
\begin{equation}
\mathbf{P} (\Xi_0 \preceq\Phi_0 ) \geq1- \sum
_{i=1}^m\exp \bigl(-c_1
\varepsilon^2\beta\operatorname{vol} (S_i ) \bigr).
\label{eqcbcoupling}
\end{equation}

If $\{\Xi_0 \preceq\Phi_0\}$ holds, then we can pair each node $u
\in
\Xi_0$ to a unique node of~$v\in\Phi_0$ in the same cell. We call $u$
the corresponding pair of~$v$ (and vice versa).
We will now show that we can couple
the motion of the nodes in $\Xi_0$ with the motion of their
corresponding pairs in
$\Phi_0$ so that the probability that an arbitrary pair is at the same
location at time $\Delta$
is sufficiently large.

To describe the coupling, let
$v'$ be a node of~$\Xi_0$ located at $y' \in S$, and let $v$ be the
corresponding pair of~$v'$ in $\Phi_0$. Let $y$ be the location of~$v$
in $S$, and note that
since $v$ and
$v'$ belong to the same cell we have $\llVert y-y'\rrVert _2 \leq
\gamma$.
Since $g$ is $(\varepsilon/2,2\gamma)$-indistinguishable from
$f_\Delta$,
we assume that the motion of~$v'$ from time $0$ to $\Delta$ is such
that its density function is a
renormalization of~$g$. Then, for any point~$z$ such that
$z-y' \in Q_{R}$, we have
\[
g \bigl(z-y' \bigr) = g \bigl(z-y- \bigl(y'-y \bigr)
\bigr) \leq f_\Delta(z-y), %
\]
since $y'-y\in Q_{2\gamma}$ and $z-y = (z-y')-(y-y')\in Q_{R+2\gamma}$.
Then we have that the function $g$ is smaller than the densities
for the motions of~$v$ and $v'$ to the location~$z$ for all $z$ such
that $z-y'\in Q_{R}$ and $y-y'\in Q_{2\gamma}$.

Define $\tilde g(x) = g(x)\mathbf{1} (x \in Q_{R} )$ and
%
%
\begin{equation}
\psi= \int_{Q_{R}} \tilde g(x) \,dx =\int_{Q_{R}}g(x)
\,dx \geq1-\frac{\varepsilon}{2}. \label{eqboundg}
\end{equation}
Hence, with probability $\psi$ we
can use the density function $\frac{\tilde g(x)}{\psi}$ to sample a
single location $y'+x$
for the position of both $v$ and $v'$ at time $\Delta$.
Then the second Poisson point process in the construction of~$\Xi$,
which we denote by $\Xi'_0$, is defined as
a Poisson point process of
intensity $\psi(1-\varepsilon/2)\beta$
obtained by \textit{thinning} $\Xi_0$ (i.e., deleting each node
of~$\Xi_0$ with probability $1-\psi$).
At this step, we have used the fact that the function
$\tilde g(x)$ is oblivious of the location of~$v$ and, consequently,
is independent of the
point process $\Phi_0$.

Then the third Poisson point process $\Xi'_\Delta$ is obtained from
$\Xi'_0$
after the nodes have moved
according to the density function $\frac{\tilde g(z)}{\psi}$.
Thus, since $g$ is translation invariant, we have that $\Xi'_\Delta$
is a Poisson point process and, due to the coupling described above,
the nodes of~$\Xi'_\Delta$
are a subset of the nodes of~$\Phi_\Delta$ and are independent of the
nodes of~$\Phi_0$, where $\Phi_\Delta$ is
obtained by letting the nodes of~$\Phi_0$ move from time $0$ to time
$\Delta$ according to the density function $f_\Delta$.

Note that $\Xi'_\Delta$ is a \emph{nonhomogeneous} Poisson point process.
It remains to show that the intensity of~$\Xi'_\Delta$ is strictly
larger than $(1-\varepsilon)\beta$ in $S'$ so
that $\Xi$ can be obtained from $\Xi'_\Delta$ via thinning; since
$\Xi
'_\Delta$ is independent of~$\Phi_0$, so is $\Xi$.

For $z \in\mathbb{R}^d$, let $\mu(z)$ be the intensity of~$\Xi
'_{\Delta}$.
Since $\Xi'_0$ has no node outside $S$, we obtain, for any $z\in S'$,
\[
\mu(z) \geq\psi(1-\varepsilon/2)\beta\int_{z+Q_{R}}
\frac{\tilde
g(z-x)}{\psi} \,dx = (1-\varepsilon/2)\beta\int_{Q_{R}}
\tilde g(x) \,dx, %
\]
where the inequality follows since $z+Q_{R} \subset S$ for all $z \in S'$.
Using~(\ref{eqboundg}), we have
$\mu(z) \geq(1-\varepsilon/2)(1-\varepsilon/2)\beta\geq
(1-\varepsilon
)\beta
$, which is the intensity of~$\Xi$.
Therefore, when $\{\Xi_0 \preceq\Phi_0\}$ holds, which occurs with
probability given by (\ref{eqcbcoupling}),
the nodes of~$\Xi$ are a subset of the nodes of~$\Phi_\Delta$, which
completes the proof of Proposition~\ref{procoupling}.
\end{pf}

Now we illustrate an application of Proposition~\ref{procoupling} by
giving an example that we will use later in our proofs.
Consider a node that, at time $0$, is located at an arbitrary position
$x\in Q_\ell$.
Let this
node move for time $\Delta$ according to a Brownian motion
and, for each $z>0$, define $F_\Delta(z)$ to be the event that this node
never leaves the cube $x+Q_z$ during\vspace*{1pt} the interval
$[0,\Delta]$.
Let $M\geq\ell$ and,
for $y=(y_1,y_2,\ldots,y_d)\in Q_{2M}$, define $\tilde f_\Delta(y)$ to
be the
probability density function for the location of this node at position
$x+y$ at time $\Delta$ conditioned on $F_\Delta(3M)$.
It follows from the reflection principle of Brownian motion that
%
%
\begin{eqnarray}
\label{eqdeff} && \tilde f_\Delta(y) \mathbf{P} \bigl(F_\Delta(3M)
\bigr)
\nonumber
\\
&&\qquad\geq\prod_{i=1}^d \biggl(
\frac{1}{\sqrt{2\pi\Delta}}\exp \biggl(-\frac
{y_i^2}{2\Delta} \biggr) -\frac{1}{\sqrt{2\pi\Delta}}\exp
\biggl(-\frac{(3M-y_i)^2}{2\Delta
} \biggr)
\nonumber
\\
&&\hspace*{143pt} {}-\frac{1}{\sqrt{2\pi\Delta}}\exp \biggl(- \frac{(3M+y_i)^2}{2\Delta
} \biggr)
\biggr)
\nonumber
\\[-8pt]
\\[-8pt]
\nonumber
&&\qquad= \frac{1}{(2\pi\Delta)^{d/2}}\exp \biggl(-\frac{\llVert
y\rrVert
_2^2}{2\Delta
} \biggr)
\\
&&\quad\qquad{}\times\prod_{i=1}^d
\biggl[1-\exp \biggl(-\frac{9M^2}{2\Delta
} \biggr) \biggl( \exp \biggl(
\frac{6My_i}{2\Delta} \biggr) +\exp \biggl(- \frac{6My_i}{2\Delta
} \biggr) \biggr)
\biggr]
\nonumber
\\
&&\qquad\geq\frac{1}{(2\pi\Delta)^{d/2}}\exp \biggl(-\frac
{\llVert y\rrVert
_2^2}{2\Delta
} \biggr)
\biggl(1-2d\exp \biggl(-\frac{3M^2}{2\Delta} \biggr) \biggr),
\nonumber
\end{eqnarray}
where the last step follows since the function $e^a+e^{-a}$ is
increasing in $\vert a\vert $ and $y_i\in[-M,M]$ for all $y\in Q_{2M}$.

Our goal is to apply Proposition~\ref{procoupling} with $f_\Delta
=\tilde f_\Delta$. In order to do this, we will use the technical lemma
below, which
constructs a subdensity function that is indistinguishable from
$\tilde f_\Delta$.
%

\begin{lemma}\label{lemindistinguishable}
Let $\xi\in(0,1)$ and $m>0$. Let
\[
\Delta\geq\frac{d^3m^2}{\xi^2}\quad\mbox{and}\quad M\geq\sqrt {8\Delta\log
\bigl(8d \xi^{-1} \bigr)}.
\]
For $z\in Q_M$, define
\[
g(z) = \frac{1}{(2\pi\Delta)^{d/2}}\exp \biggl(-\frac{(\llVert
z\rrVert
_2+m\sqrt
{d}/2)^2}{2\Delta} \biggr) \biggl(1-2d\exp
\biggl(-\frac{3M^2}{2\Delta
} \biggr) \biggr), %
\]
and,\vspace*{2pt}
for $z\notin Q_M$, set $g(z)=0$.
Then, $g$ is $(\xi,m)$-indistinguishable from $\tilde f_\Delta$ over
$Q_M$, where $\tilde f_\Delta$ is the probability density function for
the location of a Brownian
motion at time $\Delta$ given that it never leaves the cube $Q_{3M}$
during the interval $[0,\Delta]$.
\end{lemma}

\begin{pf}
Note that, for any $x\in Q_m$, the triangle inequality gives that
$\llVert
z-x\rrVert _2 \leq\llVert x\rrVert _2+\llVert z\rrVert _2 \leq
\frac{m\sqrt{d}}{2}+\llVert z\rrVert _2$.
Thus, for all $x\in Q_m$ and $z\in Q_M$, we have that $z-x\in Q_{M+m}
\subseteq Q_{2M}$ and, from~(\ref{eqdeff}), we obtain that
$g(z) \leq\tilde f_\Delta(z-x)$ as required by the first condition
in~(\ref{eqindistinguishable}).
Now, let $\rho=m\sqrt{d}/2$ and
\[
\upsilon(x) =\frac{1}{\sqrt{2 \pi\Delta}} \exp \biggl(-\frac
{(\vert x\vert +\rho
)^2}{2 \Delta} \biggr),
\]
for $x \in{{\mathbb{R}}}$.
Note that $\sum_{i=1}^{d} (\vert z_i\vert +\rho)^2 = \llVert
z\rrVert _2^2 + 2 \rho\llVert z \rrVert _1
+ d\rho^2\geq(\llVert z \rrVert _2 + \rho)^2$,
so
%
%
\begin{equation}
\label{eqprod} \quad \biggl(1-2d\exp \biggl(-\frac{3M^2}{2\Delta} \biggr) \biggr)
\prod_{i=1}^{d} \upsilon(z_i)
\leq g(z)\qquad\mbox{for } z= (z_1,\ldots,z_d) \in
Q_M.
\end{equation}
Next, observe that
\begin{eqnarray*}
\int_{-\infty}^{\infty} \upsilon(x) \,dx &=& 1 - \int
_{-\rho}^{\rho} \frac{1}{\sqrt{2\pi\Delta}} \exp \biggl(-
\frac
{y^2}{2\Delta} \biggr) \,dy
\\
&\geq& 1- \frac{2\rho}{\sqrt{2\pi\Delta}} \geq1 - \frac{\rho
}{\sqrt{\Delta}} \geq1 -
\frac{\xi}{2d},
\end{eqnarray*}
where the last step follows from $\Delta\geq\frac{d^3m^2}{\xi
^2}=\frac
{4d^2\rho^2}{\xi^2}$.
Now,
since $\frac{M-m}{2}+\rho\geq\frac{M}{2}\geq\sqrt{\Delta}$, we
apply the Gaussian tail bound (cf.~Lemma~\ref{lemgaussiantail}) to obtain
\begin{eqnarray*}
\int_{(M-m)/2}^{\infty} \upsilon(x) \,dx &\leq&
\frac{1}{\sqrt{2\pi}}\frac{\sqrt{\Delta}}{(M-m)/2+\rho
}\exp \biggl(-\frac{((M-m)/2+\rho)^2}{2\Delta} \biggr)
\\
&\leq&\frac{1}{\sqrt{2\pi}}\frac{\sqrt{\Delta}}{M/2}\exp \biggl(-\frac
{M^2}{8\Delta}
\biggr) \leq\frac{\xi}{8d},
\end{eqnarray*}
for any $\xi\in(0,1)$ since $M\geq\sqrt{8\Delta\log(8d\xi^{-1})}$.
Thus, $\int_{-(M-m)/2}^{(M-m)/2} \upsilon(x) \,dx \geq1 - \frac{\xi
}{2d} - 2\frac{\xi}{8d} = 1-\frac{3\xi}{4d}$.
We can now deduce from (\ref{eqprod}) that
\begin{eqnarray*}
\int_{Q_{M-m}} g(z) \,dz &\geq& \biggl(1-2d\exp \biggl(-
\frac{3M^2}{2\Delta} \biggr) \biggr) \int_{Q_{M-m}} \prod
_{i=1}^{d} \upsilon(z_i) \,dz
\\
&\geq& \biggl(1-2d\exp \biggl(-\frac{3M^2}{2\Delta} \biggr) \biggr) \biggl(1-
\frac{3\xi}{4d} \biggr)^d
\\
&\geq& \biggl(1-2d \biggl(\frac{\xi}{8d} \biggr)^{12} \biggr)
\biggl(1-\frac
{3\xi}{4} \biggr)
\\
&\geq& \biggl(1-\frac{\xi}{4} \biggr) \biggl(1-\frac{3\xi}{4} \biggr)
\geq1 -\frac{\xi}{4}-\frac{3\xi}{4}= 1 - \xi.
\end{eqnarray*}
\upqed
\end{pf}

The next proposition is a special case of Proposition~\ref
{procoupling} that we will use in our proofs.
%

\begin{proposition}\label{procouplingapplied}
Fix $K > \ell>0$ and consider the cube $Q_K$ tessellated into
subcubes of side length $\ell$.
Let $\Phi_0$ be an arbitrary point process at time $0$
that contains at least $\beta\ell^d$ nodes in each subcube for some
$\beta>0$.
For any $z> 0$, let $\Phi_\Delta(z)$ be the point process obtained by
letting the nodes of~$\Phi_0$ move for time $\Delta$
according to independent Brownian motions that are conditioned on
being inside $Q_z$ throughout the interval $[0,\Delta]$.
Fix $\varepsilon\in(0,1)$.
There are constants $c_1,c_2,c_3$ depending only on $d$ such that, if
$\Delta\geq\frac{c_1 \ell^2}{\varepsilon^2}$ and $K' \leq K - c_2
\sqrt
{\Delta\log(16d\varepsilon^{-1})}>0$,
we can couple the nodes of~$\Phi_{\Delta} (3 (K-K'+2\sqrt
{d}\ell
) )$
with those of a Poisson point process $\Xi$ that is independent
of~$\Phi_0$ and has intensity $(1-\varepsilon)\beta$
so that\vspace*{1pt}
the nodes of~$\Xi$ are a subset of the nodes of~$\Phi_\Delta
(3
(K-K'+2\sqrt{d}\ell) )$ inside the cube $Q_{K'}$ with
probability at least
\[
1-\frac{K^d}{\ell^d}\exp \bigl(-c_3\varepsilon^2\beta
\ell^d \bigr). %
\]
\end{proposition}

\begin{pf}
Denote by $\tilde f_\Delta(y)$ the probability density function for
the location of a Brownian motion at time $\Delta$ conditioned on
the motion never leaving $Q_{3(K-K'+2\sqrt{d}\ell)}$ during the whole
of~$[0,\Delta]$.
%
%
Now, if $c_1$ and $c_2$ are large enough with respect to $d$,
we can apply Lemma~\ref{lemindistinguishable} with $\xi=\varepsilon/2$,
$M=K-K'+2\sqrt{d}\ell$ and $m=2\sqrt{d}\ell$, which gives that
$g$ is $(\varepsilon/2,2\sqrt{d}\ell)$-indistinguishable from
$\tilde
f_\Delta$ over $Q_{K-K'+2\sqrt{d}\ell}$.
Thus, we apply Proposition~\ref{procoupling} with $S=Q_K$ and
$S'=Q_{K'}$. This gives that $R=K-K'$ and $\gamma=\ell\sqrt{d}$.
Thus, we have that a Poisson point process $\Xi$ with intensity
$(1-\varepsilon)\beta$ over $Q_{K'}$
can be coupled with $\Phi_\Delta(3(K-K'+2\sqrt{d}\ell)
)$ so
that $\Xi$ is a subset of~$\Phi_\Delta(3(K-K'+2\sqrt{d}\ell) )$
with probability
at least
\[
1-\frac{K^d}{\ell^d}\exp \bigl(-c_3\varepsilon^2\beta
\ell^d \bigr)\qquad\mbox{for some positive constant
}c_3=c_3(d). %
\]
\upqed
\end{pf}

\section{Space--time percolation}\label{secfp}

In this section, we develop the main technical result of this paper,
which we stated in a simplified form in Theorem~\ref{thmmainresult}.
In order to state this theorem in full generality, we need to
introduce some notation, which will extend the setting introduced in
the part \emph{space--time percolation} in Section~\ref{secintro}.

We tessellate $\mathbb{R}^d$ into cubes of side length $\ell$.
We index the cubes by integer vectors $\veci\in\mathbb{Z}^d$ such that
\[
\mbox{cube }\veci=(i_1,i_2,\ldots,i_d)
\mbox{ corresponds to the region }\prod_{j=1}^d
\bigl[i_j\ell,(i_j+1)\ell \bigr]\subset
\mathbb{R}^d. %
\]
Let $t>0$, and tessellate the time interval $[0,t]$ into subintervals of
length $\beta$. We index the subintervals by $\tau\in\mathbb
{Z}_+$, where
\[
\mbox{subinterval }\tau\mbox{ represents the time interval } \bigl[\tau\beta,(
\tau+1)\beta \bigr]\subset\mathbb{R}. %
\]
We call each pair $(\veci,\tau)$ of the space--time tessellation
a \emph{cell}, and define the region of a cell as
\[
R_1(\veci,\tau)=\prod_{j=1}^d
\bigl[i_j\ell,(i_j+1)\ell \bigr] \times \bigl[\tau
\beta,(\tau+1)\beta \bigr]. %
\]

The space--time tessellation defined above will be sufficient when we
apply in Section~\ref{secdetection}
the technique developed in this section to the detection problem.
However, in previous works, being able to handle overlapping cells
turned out to be very useful; for example, in the study
of the spread of infection~\cite{KS05} and in the estimation of the
length of the shortest path between nodes~\cite{FSS13}.
For this reason, we generalize our framework slightly by introducing
bigger, overlapping cells.
Consider an integer parameter
$\eta\geq1$ and, for each cube $\veci=(i_1,i_2,\ldots,i_d)$ and time
interval $\tau$,
\begin{eqnarray*}
&\displaystyle\mbox{define the \emph{super cube} }\veci\mbox { as }\prod
_{j=1}^d \bigl[i_j
\ell,(i_j+\eta)\ell \bigr]\quad\mbox{and}&
\\
&\displaystyle\mbox{the \emph{super interval} } \tau\mbox{ as } \bigl[\tau
\beta,(\tau+\eta)\beta \bigr]. &
\end{eqnarray*}
Then define the \emph{super cell} $(\veci,\tau)$ as the Cartesian
product of the super cube $\veci$ and the super interval $\tau$.

For a sequence of point processes $(\Pi_s)_{s\geq0}$,
we say that
\[
\begin{tabular} {p{310pt}} an\vspace*{2pt} event $E$ is \emph
{increasing} for $( \Pi_s)_{s\geq0}$ if the fact that $E$
holds for $( \Pi_s)_{s\geq0}$ implies that it holds for all $
( \Pi'_s )_{s\geq0}$ for which $
\Pi'_s \supseteq\Pi_s$ for all $s\geq0$.
\end{tabular} %
\]
We also say that an event $E$ is \emph{restricted} to a region
$X\subset\mathbb{R}^d$ and a time interval $[t_0,t_1]$
if it is measurable with respect to the $\sigma$-field generated
by the nodes that are inside $X$ at time $t_0$ and their positions from
time $t_0$ to $t_1$.
For an increasing event $E$ that is restricted to a region $X$ and a
time interval $[t_0,t_1]$, we have the following definition:
%
%
\begin{equation}
\label{eqdefnu} %
\begin{tabular} {p{280pt}} $\nu_E$ is
called the \emph{probability associated to $E$} if, for any $\mu\geq0$ and region
$X'\subset\mathbb{R}^d$, $\nu_E (
\mu,X' )$ is the probability that $E$ happens given that, at
time $t_0$, the nodes in $X$ are given by a Poisson point process of
intensity $\mu$ and their motions from $t_0$ to $t_1$
are independent Brownian motions conditioned to have displacement inside
$X'$ from $t_0$ to $t_1$. \end{tabular}
\end{equation}
In other words, if $(x_t)_{t\geq0}$ are the locations of one such
node, then $x_{t_0}\in X$ and, for each $s\in[t_0,t_1]$, we have
$x_s-x_{t_0}\in X'$.

We say that:
\[
\begin{tabular} {p{280pt}} two distinct cells $(\veci,\tau)$ and $
( \veci', \tau' )$ are \emph{adjacent} if $
\llVert\veci- \veci' \rrVert_\infty\leq1$ and $
\vert\tau-\tau' \vert\leq1$. \end{tabular} %
\]
For each $(\veci,\tau)\in\mathbb{Z}^{d+1}$,
let $E(\veci,\tau)$ be the indicator random variable for an increasing
event restricted to the super cube $\veci$ and the super interval
$\tau$.
We say that a cell $(\veci,\tau)$ is \emph{bad} if $E(\veci,\tau)=0$,
and in this case,
we define the \emph{bad cluster} $K(\veci,\tau)$ of~$(\veci,\tau)$
as the following set of cells:
%
%
\begin{eqnarray}
\label{eqdefk} K(\veci,\tau) &=& \bigl\{ \bigl(\veci',
\tau' \bigr)\in\mathbb{Z}^{d+1} \dvtx E \bigl(\veci
',\tau' \bigr)=0\mbox{ and } \exists\mbox{ a path}
\nonumber
\\[-8pt]
\\[-8pt]
\nonumber
&&\hspace*{5pt}\mbox{of adjacent bad cells from }( \veci,\tau)\mbox{
to } \bigl(\veci',\tau' \bigr) \bigr\}.
\end{eqnarray}
If $E(\veci,\tau)=1$, we let $K(\veci,\tau)=\varnothing$.
Finally, define
%
%
\begin{equation}
\mathcal{R}^t_1= \bigl\{(\veci,\tau)\in
\mathbb{Z}^{d+1}\dvtx R_1(\veci,\tau)\subset(-t,t)
\times[0,t) \bigr\}. \label{eqrt}
\end{equation}
Our main technical result, Theorem~\ref{thmfp} below, establishes a
bound on the bad cluster of the origin.
%

\begin{theorem}\label{thmfp}
For each $(\veci,\tau)\in\mathbb{Z}^{d+1}$,
let $E(\veci,\tau)$ be the indicator random variable for an increasing
event that is restricted to the super cube $\veci$ and the super interval
$\tau$, and let $\nu_E$ be the probability associated to $E$ as
defined in~(\ref{eqdefnu}).
Fix a constant $\varepsilon\in(0,1)$, an integer $\eta\geq1$ and the
ratio $\beta/\ell^2>0$.
Fix also $w$ such that
\[
w\geq\sqrt{18\eta\frac{\beta}{\ell^2}\log \biggl(\frac
{8d}{\varepsilon
} \biggr)}.
\]
Then there exist positive numbers $\alpha_0$ and $t_0$ so that if
\[
\alpha=\min \biggl\{\varepsilon^2\lambda\ell^d, \log
\biggl(\frac{1}{1-\nu_E((1-\varepsilon)\lambda,Q_{w\ell})} \biggr) \biggr\} \geq\alpha_0, %
\]
we have a positive constant $c=c(d)$ and a positive value $C$
independent of~$t$ such that, for all $t\geq t_0$,
\[
\mathbf{P} \bigl(K(\origin,0) \nsubseteq\mathcal{R}^t_1
\bigr) \leq\cases{ \displaystyle\exp \biggl(-C \frac{\sqrt
{t}}{(\log t)^{c}} \biggr), &\quad
for $d=1$,
\cr
\noalign{\vspace*{3pt}} \displaystyle\exp \biggl(-C
\frac{t}{(\log t)^{c}} \biggr), &\quad for $d=2$,
\cr
\noalign{\vspace*{3pt}}
\displaystyle\exp(-C t ), &\quad for $d\geq3$.} %
\]
\end{theorem}

%
%
\begin{remark}
(i)~Theorem~\ref{thmmainresult} is a special case of Theorem~\ref
{thmfp} with $\eta=1$ and $E(\veci,\tau)$ corresponding to
the event that the cell $\veci$ has at least $(1-\varepsilon)\lambda
\ell
^d$ nodes at time $\tau\beta$. Note that, in this case, $w$ can be
taken to
be any positive constant since $\nu_E$ does not depend on $w$.

(ii)~In Theorem~\ref{thmfp}, $C$ may depend on all the parameters
(including $\lambda$), but does not depend on $t$.

(iii)~Note that, in Theorem~\ref{thmfp}, we first fix the values
of~$\varepsilon$, $\eta$ and $\beta/\ell^2$. Given these values, we
fix $w$,
and only then we set $\lambda$ and $\ell$ so that the condition on
$\alpha$ is satisfied. Note that setting $\ell$ in this last step is
equivalent to
setting $\beta$, since the ratio $\beta/\ell^2$ is fixed. It is
important that $w$ does not depend on $\lambda$ or $\ell$. In typical
applications
of Theorem~\ref{thmfp}, $w$ is set to be sufficiently large so that
$\nu_E((1-\varepsilon)\lambda,Q_{w\ell})$ can be made arbitrarily
close to 1
by setting $\lambda$ or $\ell$ large enough after having fixed~$w$.
This is usually the case since, if $w$ is large enough with respect to
$\eta$ and $\beta/\ell^2$, we have that the probability that a
standard Brownian motion stays inside $Q_{w\ell}$ during $[0,\eta
\beta]$
is at least $1-\exp(-\frac{cw^2\ell^2}{\eta\beta} )$ for
some positive constant $c$. Then we have that most of the
nodes of the Poisson point process of intensity $(1-\varepsilon
)\lambda$
in the definition of~$\nu_E$
stay inside the cube of side length $w\ell$ that is centered at their
initial position.

(iv)~As observed in~\cite{Kesidis,Konst}, if we define $E(\veci,\tau
)=0$ when there is no node in the cube $\veci$ at time $\tau\beta$, then
$\mathbf{P} (K(\origin,0) \nsubseteq\mathcal
{R}^t_1 )$ is at least $\exp
(-c \sqrt{t} )$ for $d=1$,
$\exp(-c \frac{t}{\log t} )$ for $d=2$ and
$\exp(-c t )$ for $d\geq3$, where $c$ does not depend on
$t$. This lower bound is achieved by the event that
$E(\origin,\tau)=0$ for all $\tau=0,1,\ldots,\lceil t/\beta\rceil$.
Hence, with respect to $t$,
the exponents in the bound of Theorem~\ref{thmfp} are tight up to
logarithmic factors for
$d=1,2$ and up to constants for $d\geq3$.
\end{remark}

In order to prove Theorem~\ref{thmfp}, we will do a multi-scale analysis,
where, at each scale, we tessellate $\mathbb{R}^{d+1}$ into cells of
a given size. We then analyze the density of nodes in these cells using
the framework of fractal percolation.
The proof of Theorem~\ref{thmfp} is rather lengthy and, for this
reason, we split the proof in many parts.
We start by introducing the multi-scale tessellation in Section~\ref
{sectess}. Then, in
Section~\ref{secfpprocess}, we define the fractal percolation process
we will analyze.
In Section~\ref{sechighlevel}, we present a high-level sketch of the
proof using the notation established in the previous two sections, and
also discuss the relation with standard fractal percolation.
In Section~\ref{secsupport}, we define the support of a cell and prove
some geometric properties of the support.
In Section~\ref{secscpath}, we use the definition of the support of a
cell to introduce a new notion of path of cells, and also
prove some lemmas about the paths.
Finally, in Section~\ref{secproof}, we put these results together to
give the proof of Theorem~\ref{thmfp}.

\subsection{Multi-scale tessellation of space and time}\label{sectess}

We start with the tessellation of space.
Let $m$ be a sufficiently large integer.
For each scale $k\geq1$, tessellate $\mathbb{R}^d$
into cubes of length $\ell_k$ such that
%
%
\begin{equation}
\ell_1=\ell\quad\mbox{and}\quad\ell_k=m k^3
\ell_{k-1}=m^{k-1}(k!)^3 \ell. \label{eqell}
\end{equation}
We also set $\ell_0=\ell/m$.

%
%
\begin{figure}[b]

\includegraphics{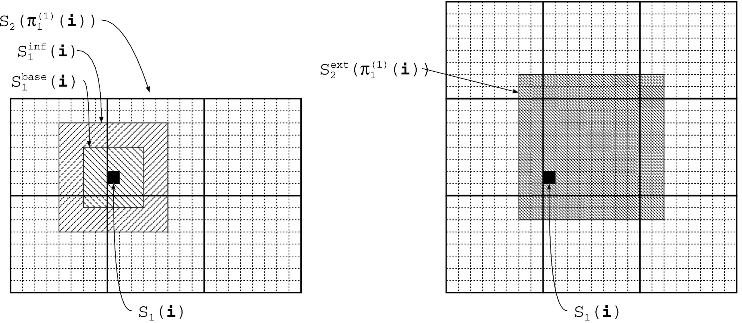}

\caption{Illustration of the tessellation of~$\mathbb{R}^d$. Different
scales are represented by the thickness of the lines; for example,\vspace*{1pt}
$S_2(\pi_1^{(1)}(\veci))$
is the square with thick borders that contains $S_1(\veci)$, which is
the black square.
Note that $S_1(\veci)$ is at the same position in both the left and
the right pictures above,
illustrating that $S^\mathrm{base}_1(\veci) \subset S^\mathrm
{ext}_{2}(\pi_1^{(1)}(\veci
))$ as given in~(\protect\ref{eqrelcubes}).}\label{figspacescale}
\end{figure}

%
%
\begin{remark}\label{remdepm}
We let $m$ be sufficiently large with respect to $\varepsilon$, $\eta$,
$w$ and the ratio $\beta/\ell^2$; however,
$m$ does not depend on $\lambda$ or $\ell$. Hence, after having $m$
fixed, we can make $\alpha$ arbitrarily large
by setting either $\lambda$ or $\ell$ large enough.
\end{remark}

It will be useful to refer to Figure~\ref{figspacescale} in the
discussion below.
We index the cubes by integer vectors $\veci\in\mathbb{Z}^d$ and denote
them by
$S_k(\veci)$. Therefore, for $\veci=(i_1,i_2,\ldots,i_d)$, we have
\[
S_k(\veci) = \prod_{j=1}^d
\bigl[i_j\ell_k,(i_j+1)\ell_k
\bigr]. %
\]
So, $S_1(\veci)$ is exactly the cube $\veci$ introduced in the
beginning of Section~\ref{secfp}.
Note that $S_k(\veci)$ is the union of~$(mk^3)^d$ cubes of scale $k-1$;
in particular, a cube of scale smaller than $k$ is contained inside a
unique cube of scale $k$. With this,
for each $k$, $j\geq0$ and $\veci\in\mathbb{Z}^d$, we can define
\[
\pi_{k}^{(j)}(\veci)=\veci'\quad\mbox{iff}
\quad S_k(\veci) \subseteq S_{k+j} \bigl(
\veci' \bigr). %
\]
Thus, $\pi_k^{(j)}$ defines a hierarchy over cubes since, for all
$j'\leq j$,
$\pi_k^{(j)}(\veci)=\pi_{k+j'}^{(j-j')}(\pi_k^{(j')}(\veci))$.
Using this notion, for $\veci,\veci'\in\mathbb{Z}^d$ and $k\geq0$, we
say that $(k+1,\veci')$ is the \emph{parent} of~$(k,\veci)$ if $\pi
_k^{(1)}(\veci)=\veci'$; in this case, we also say
that $(k,\veci)$ is a \emph{child} of~$(k+1,\veci')$.
We define the set of descendants of~$(k,\veci)$ as $(k,\veci)$ and
the union of the descendants of
the children of~$(k,\veci)$, or only $(k,\veci)$ in case $(k,\veci)$
has no child.

We now introduce a new variable $n$ that satisfies
%
%
\begin{equation}
n^d=\frac{m}{7\eta}. \label{eqdefn}
\end{equation}
The variable $m$ must be large enough so that $n> 1$.
We also assume that $m$ is specified in a way that makes $n$ an integer.

We define some larger cubes based on $S_k(\veci)$.
For $k\geq0$, define the \emph{base} $S^\mathrm{base}_k(\veci)$
of~$S_k(\veci)$
and the \emph{area of influence} $S^\mathrm{inf}_k(\veci)$
of~$S_k(\veci)$ by
\[
S^\mathrm{base}_k(\veci) = \bigcup
_{\veci'\dvtx\llVert \veci-\veci'\rrVert
_\infty\leq
\eta mn(k+1)^3} S_k \bigl(\veci' \bigr)
\]
and
\[
S^\mathrm{inf}_k(\veci)=\bigcup_{\veci'\dvtx\llVert \veci-\veci
'\rrVert
_\infty\leq
2\eta mn(k+1)^3}
S_k \bigl(\veci' \bigr). %
\]
We also define the \emph{extended} cube $S^\mathrm{ext}_k(\veci)$ by
\[
S^\mathrm{ext}_k(\veci)=\bigcup_{\veci'\dvtx\pi
_{k-1}^{(1)}(\veci')=\veci
}S^\mathrm{base}_{k-1}
\bigl(\veci' \bigr). %
\]
Note that $S^\mathrm{ext}_k(\veci)$ is the union of the bases of the
children of~$(k,\veci)$, which are the $(k-1)$-cubes contained in
$S_k(\veci)$.
It is easy to see that $S_k(\veci)\subset S^\mathrm{base}_k(\veci)
\subset S^\mathrm{inf}
_k(\veci)$ and
%
%
\begin{equation}
S^\mathrm{ext}_{k+1} \bigl(\pi_{k}^{(1)}(
\veci) \bigr) = \bigcup_{\veci'\dvtx\pi_{k}^{(1)}(\veci')=\pi
_k^{(1)}(\veci
)}S^\mathrm{base}
_{k} \bigl(\veci' \bigr) \supset S^\mathrm{base}_k(
\veci). \label{eqrelcubes}
\end{equation}

%
\begin{remark}\label{remextcontainssuper}
One important property obtained from these definitions is that an
extended cube of scale 1 has side length
$\ell+ 2\eta m n\ell_0 = (1+2\eta n)\ell$. Therefore, for any $\veci
\in\mathbb{Z}^d$, the extended cube
$S^\mathrm{ext}_1(\veci)$ contains the super cell $\veci$ defined in the
beginning of Section~\ref{secfp}.
\end{remark}

\begin{figure}

\includegraphics{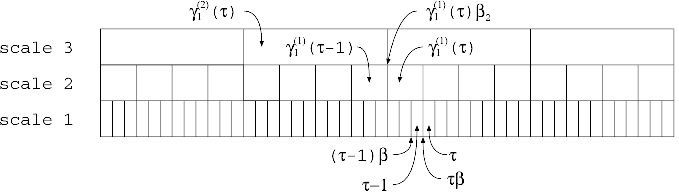}

\caption{Time scale. The horizontal axis represents time and the
vertical axis represents the scale.
Note that $\gamma_1^{(1)}(\tau)=\gamma_1^{(1)}(\tau+1)=\gamma
_1^{(1)}(\tau+2)$.}
\label{figtimescale}
\end{figure}

Now, we define a multi-scale tessellation of time.
In this discussion, it will be useful to refer to Figure~\ref{figtimescale}.
We define
\[
\varepsilon_1=\varepsilon\quad\mbox{and}\quad\varepsilon_k=
\varepsilon_{k-1} - \frac{\varepsilon}{k^2}\qquad\mbox{for all }k\geq1.
\]
There will be no scale $0$ for time, but we have defined $\varepsilon
_0=2\varepsilon$ for consistency.
Now we define
%
%
\begin{equation}
\beta_k=C_\mathrm{mix}\frac{\ell_{k-1}^2}{(\varepsilon
_{k-1}-\varepsilon_k)^2} =
C_\mathrm{mix}\frac{\ell_{k-1}^2k^4}{\varepsilon^2}\qquad\mbox {for all }k\geq1,
\label{eqdefbeta}
\end{equation}
where $C_\mathrm{mix}\geq4c_1$ and $c_1$ is the constant in
Proposition~\ref{procouplingapplied} that depends on $d$ only.
Note that, for $k=1$, we have
%
%
\begin{equation}
\beta=\beta_1 = C_\mathrm{mix}\frac{(\ell/m)^2}{\varepsilon^2}.
\label{eqdefcmix}
\end{equation}
Hence, given $\beta/\ell^2$ and $\varepsilon$, we can set $m$
sufficiently large
so that $C_\mathrm{mix}\geq4c_1$.
Also, note that
%
%
\begin{equation}
\frac{\beta_{k+1}}{\beta_k} =\frac{\ell_k^2(k+1)^4}{\ell
_{k-1}^2k^4} =m^2k^2(k+1)^4
\qquad\mbox{for all }k\geq1. \label{eqbetaratio}
\end{equation}

Now, for scale $k\geq1$, we tessellate time into intervals of length
$\beta_k$.
We index the time intervals by $\tau\in\mathbb{Z}$ and denote them by
$T_k(\tau)$, where
\[
T_k(\tau)= \bigl[\tau\beta_k,(\tau+1)\beta_k
\bigr]. %
\]
%
We allow time to be negative, and note that $\beta_{k+1}/\beta_{k}$ is
always an integer by~(\ref{eqbetaratio}), which gives that
%
%
\begin{equation}
\label{eqpropertimetess} %
\begin{tabular} {p{285pt}} a time interval of scale
$k$ is contained inside a unique time interval of scale $k+1$. \end{tabular}
\end{equation}
At scale $1$, we will consider the time intervals that intersect $[0,t]$.

%
%

We also introduce a hierarchy over time, but it is important to
emphasize that this time hierarchy is
conceptually different than the spatial hierarchy induced by~$\pi$.
For all $k$ and $\tau$, let $\gamma_{k}^{(0)}(\tau)=\tau$, and, for
$j\geq1$, define
\[
\gamma_{k}^{(j)}(\tau)=\tau'\qquad\mbox{if }
\gamma_k^{(j-1)}(\tau)\beta_{k+j-1}\in
T_{k+j} \bigl(\tau'+1 \bigr). %
\]
In the spatial tessellation, we constructed the hierarchy $\pi$ based
on whether a large cell contains a small cell. For the time tessellation,
if $\tau'=\gamma_k^{(1)}(\tau)$, then the interval at scale $k+1$ that
contains $T_k(\tau)$ [which is uniquely defined as noted in~(\ref
{eqpropertimetess})] is
$T_{k+1}(\tau'+1)$.
In other words, as we move from scale $k+1$ to $k$, this hierarchy
shifts time forward. The reason for this
definition is that
this shift in time will be used to allow nodes to
mix inside their cells using the techniques of Section~\ref{secmixing}.
This will become clearer in Section~\ref{secscpath}.

Now note that $\gamma$ indeed establishes a hierarchy over time since,
for any $j'\leq j$, we have $\gamma_{k}^{(j)}=\gamma
_{k+j'}^{(j-j')}(\gamma_{k}^{(j')})$.
Thus, for $\tau,\tau'\in\mathbb{Z}$ and $k\geq1$, we say that
$(k+1,\tau')$ is the \emph{parent} of~$(k,\tau)$, if $\gamma
_k^{(1)}(\tau)=\tau'$; in this case we also say
that $(k,\tau)$ is a \emph{child} of~$(k+1,\tau')$.
We also define the set of descendants of~$(k,\tau)$ as $(k,\tau)$ and
the union of the descendants of the
children of~$(k,\tau)$, or only $(k,\tau)$ in case $(k,\tau)$ has no child.

Now, for any $\veci\in\mathbb{Z}^d,k\geq1,\tau\in\mathbb{Z}$, we define
the space--time parallelogram
\[
R_k(\veci,\tau)=S_k(\veci)\times T_k(\tau),
\]
and note that these parallelograms induce a tessellation over space and
time. From now on, we reserve the word \emph{cube} to refer to the
spatial region
defined by $S_k(\veci)$, \emph{interval} to refer to the time region
defined by $T_k(\tau)$ and \emph{cell} to refer to the space--time region
$R_k(\veci,\tau)$.

For $k\geq1$, define $\mathcal{S}_k$ to be the set of indices $\veci
\in
\mathbb{Z}^d$ given by
\[
\mathcal{S}_k = \bigl\{\veci\in\mathbb{Z}^d \dvtx
S_k(\veci)\mbox{ intersects } [-t,t]^d \bigr\}.
\]
Equivalently, we can see $\mathcal{S}_k$ as the
set of cubes of scale $k$ that have a descendant of scale $1$
intersecting $[-t,t]^d$. Similarly, we define
$\mathcal{T}_k$ as the set of indices $\tau\in\mathbb{Z}$ for time
intervals of scale $k$ that have a descendant at scale $1$ intersecting
$[0,t]$. More formally, we
have
\[
\mathcal{T}_k = \bigl\{\tau\in\mathbb{Z} \dvtx\exists
\tau' \mbox{ s.t. } \gamma_1^{(k-1)} \bigl(
\tau' \bigr)=\tau\mbox{ and } T_1 \bigl( \tau
' \bigr)\mbox{ intersects } [0,t] \bigr\}. %
\]
Note that an interval in $\mathcal{T}_k$ with $k\geq2$ may not
intersect $[0,t]$.
Using these definitions, we set
\[
\mathcal{R}_k = \mathcal{S}_k \times
\mathcal{T}_k. %
\]
Note that $\mathcal{R}_1 \supset\mathcal{R}_1^t$, where $\mathcal
{R}^t_1$ is defined in~(\ref{eqrt}). Also, for any cell $(\veci,\tau
)\in\mathcal{R}_1^t$, the cells
that are adjacent to $(\veci,\tau)$ belong to $\mathcal{R}_1$.
Now define
\[
\mathcal{R} = \bigl\{(k,\veci,\tau) \dvtx1\leq k \leq\kappa\mbox { and } (\veci,
\tau)\in\mathcal{R}_k \bigr\}; %
\]
that is, $\mathcal{R}$ is set of tuples of the form $(k,\veci,\tau)$
giving the index set of all cells of all scales we consider.
Later in the proof of Theorem~\ref{thmfp}
we will set $\kappa=c\frac{\log t}{\log\log t}$ for
some sufficiently large constant $c$, but we define the tessellation
now for all scales.
We extend the hierarchy of space and time to the cells of~$\mathcal
{R}$. Then, for $(k,\veci,\tau)\in\mathcal{R}$, we define
the \emph{descendants} of~$(k,\veci,\tau)$ as the cells $(k',\veci
',\tau
')$ so that $(k',\veci')$ is a descendant of~$(k,\veci)$ and
$(k',\tau')$ is a descendant of~$(k,\tau)$.
As usual, we say that $(k,\veci,\tau)$ is an \emph{ancestor}
of~$(k',\veci',\tau')$ if $(k',\veci',\tau')$ is a descendant
of~$(k,\veci
,\tau)$.

\subsection{A fractal percolation process}\label{secfpprocess}
Here, we define the percolation process we will analyze.
For $k\geq1$, define $S_k(\veci)$ to be $k$-\emph{dense} if each
cube of scale $k-1$ that is contained in $S_k(\veci)$ has
at least $(1-\varepsilon_{k})\lambda\ell_{k-1}^d$ nodes.
For $(k,\veci,\tau)\in\mathcal{R}$, let $D_k(\veci,\tau)$ be the
indicator random variable such that
\[
D_k(\veci,\tau)=1\quad\mbox{iff}\quad S_k(\veci)\mbox{
is $k$-dense at time $\tau\beta_k$.} %
\]
Now, take a node $v\in\Pi_0$ and let $(x_t)_t$ be the locations of~$v$
in $\mathbb{R}^d$.
For a time interval $[t_0,t_1]\subset\mathbb{R}$ and a region
$X\subset
\mathbb{R}^d$, we say that:
\[
\mbox{the \emph{displacement} of~$v$ throughout }[t_0,t_1]
\mbox{ is in $X$ if } \bigcup_{s=t_0}^{t_1}(x_s-x_{t_0})
\subset X. %
\]
In other words, $v$ never leaves the region $x_{t_0}+X$ during the
interval $[t_0,t_1]$.

We need to introduce a more restrictive notion of density. For this, we
define the
indicator random variable $D^\mathrm{ext}_k(\veci,\tau)$
for each $(k,\veci,\tau)\in\mathcal{R}$ so that
\[
\begin{tabular} {p{300pt}} $D^\mathrm{ext}_k(
\veci,\tau)=1$ iff, at time $\tau\beta_k$, all cubes of scale $k-1$
contained in $S^\mathrm{ext}_k( \veci)$ have at least $(1-
\varepsilon_{k})\lambda\ell_{k-1}^d$ nodes whose
displacement throughout $ [\tau\beta_k,(\tau+2)
\beta_k ]$ is in $Q_{\eta m n k^3\ell_{k-1}}$, \end{tabular} %
\]
where $Q_z$ denotes the cube of side length $z$ given by $[-z/2,z/2]^d$.
Clearly,
$D^\mathrm{ext}_k(\veci,\tau) \leq D_k(\veci,\tau)$ for all
$(k,\veci,\tau)\in
\mathcal{R}$.
%

\begin{remark}\label{remdextfillbase}
An important property of this definition is that, when\break $D^\mathrm
{ext}_{k}(\veci,\tau)=1$, if $(k-1,\veci',\tau')$ is a child
of~$(k,\veci,\tau)$, then
we know that there are enough nodes in $S^\mathrm{base}_{k-1}(\veci
')$ at time
$\tau\beta_k$ and these nodes never leave the cube $S^\mathrm
{inf}_{k-1}(\veci
')$ during
the interval $[\tau\beta_k,\tau'\beta_{k-1}]$. This will allow us to
apply the mixing technique of Section~\ref{secmixing} to show that if
$D^\mathrm{ext}_{k}(\veci,\tau)=1$ then $D^\mathrm{ext}_{k-1}(\veci
',\tau')$ is likely to
be $1$.
\end{remark}

Now, for $k\leq\kappa-1$, let
\[
\begin{tabular} {p{300pt}} $D^\mathrm{base}_k(
\veci,\tau)=1$ iff, at time $\gamma_k^{(1)}(\tau)
\beta_{k+1}$, all cubes of scale $k$ inside $S^\mathrm{base}_k(
\veci)$ contain at least $(1-\varepsilon_{k+1})\lambda
\ell_{k}^d$ nodes whose displacement throughout $ [
\gamma_k^{(1)}( \tau)\beta_{k+1},\tau
\beta_k ]$ is in $Q_{\eta m n (k+1)^3\ell_{k}}$. \end{tabular} %
\]
Note that if $D^\mathrm{ext}_{k+1}(\pi_k^{(1)}(\veci),\gamma
_k^{(1)}(\tau))=1$
then $D^\mathrm{base}_k(\veci,\tau)=1$. This gives that
%
%
\begin{equation}
D^\mathrm{base}_k(\veci,\tau)\geq D^\mathrm{ext}_{k+1}
\bigl(\pi_k^{(1)}(\veci),\gamma_k^{(1)}(
\tau) \bigr)\qquad\mbox{for all }(k,\veci,\tau)\in\mathcal{R}. \label{eqdbasedext}
\end{equation}

For scale $\kappa$, we define
\[
A_\kappa(\veci,\tau) = D^\mathrm{ext}_\kappa(\veci,\tau),
\]
and, for $k$ satisfying $2\leq k \leq\kappa-1$, we set
\[
A_k(\veci,\tau) = \max \bigl\{D^\mathrm{ext}_k(
\veci,\tau),1- D^\mathrm{base}_k(\veci,\tau) \bigr\}.
\]
Finally, for scale $1$ we set
\[
A_1(\veci,\tau) = \max \bigl\{E(\veci,\tau),1-D^\mathrm
{base}_1(\veci,\tau) \bigr\} %
\]
and define
%
%
\begin{equation}
A(\veci,\tau)=\prod_{k=1}^{\kappa}
A_k \bigl(\pi_1^{(k-1)}(\veci),\gamma
_1^{(k-1)}(\tau) \bigr). \label{eqa}
\end{equation}
We will analyze the random variables $A(\veci,\tau)$ instead
of~$E(\veci,\tau)$.
We say that a cell $(\veci,\tau)\in\mathcal{R}_1$ has a \emph{bad
ancestry} if $A(\veci,\tau)=0$, and in this case we define
\begin{eqnarray*}
K'(\veci,\tau)&=& \bigl\{ \bigl(\veci',
\tau' \bigr)\in\mathbb{Z}^{d+1} \dvtx A \bigl(\veci
',\tau' \bigr)=0\mbox{ and $\exists$ a path of
adjacent cells}
\\
&&\hspace*{5pt}\mbox{from }(\veci,\tau)\mbox{ to } \bigl(\veci',
\tau' \bigr)\mbox{ where each cell of the path has}\\
&&\hspace*{3pt}\mbox{ a bad ancestry}
\bigr\}.
\end{eqnarray*}
If $A(\veci,\tau)=1$ then $K'(\veci,\tau)=\varnothing$.
In words, $K'(\veci,\tau)$ is the set of
cells of scale $1$ with bad ancestry that can be reached from $(\veci
,\tau)$ via a path of cells of scale $1$ with bad ancestry.
The lemma below shows that we can bound $K(\veci,\tau)$ by $K'(\veci
,\tau)$.
%

\begin{lemma}\label{lemkkp}
For each cell $(\veci,\tau)\in\mathbb{Z}^{d+1}$ of scale $1$, we have
that $E(\veci,\tau)\geq A(\veci,\tau)$.
This implies that $K(\veci,\tau) \subseteq K'(\veci,\tau)$.
\end{lemma}

\begin{pf}
We first fix $(\veci,\tau)\in\mathbb{Z}^{d+1}$. Then, for $k=1$,
define $X_1=E(\veci,\tau)$ and, for $k\geq2$, define
$X_k=D^\mathrm{ext}_k(\pi_1^{(k-1)}(\veci),\gamma_1^{(k-1)}(\tau
))$. Also, let
$Y_k=\break D^\mathrm{base}_k(\pi_1^{(k-1)}(\veci),\gamma
_1^{(k-1)}(\tau))$.
Therefore, by the definition of~$A$ in~(\ref{eqa}), we have
%
%
\begin{equation}
A(\veci,\tau)= \Biggl(\prod_{k=1}^{\kappa-1}\max
\{X_k,1-Y_k\} \Biggr)X_\kappa. \label{eqsimpla}
\end{equation}
By~(\ref{eqdbasedext}), we have that $Y_k\geq X_{k+1}$ for all $k$. Thus,
for any $k\leq\kappa-1$, we have
\[
\max\{X_{k},1-Y_k\}X_{k+1} \leq\max
\{X_{k},1-X_{k+1}\}X_{k+1} =X_{k}X_{k+1}.
\]
Therefore, applying the inequality above repetitively in~(\ref
{eqsimpla}) we have
\[
A(\veci,\tau) \leq \Biggl(\prod_{k=1}^{\kappa-2}
\max\{X_k,1-Y_k\} \Biggr)X_{\kappa
-1}X_\kappa
\leq\prod_{k=1}^\kappa X_k \leq
X_1 = E(\veci,\tau). %
\]
\upqed
\end{pf}

\subsection{High-level overview}\label{sechighlevel}

In this section, we give some intuition for the definitions
of Sections~\ref{sectess} and~\ref{secfpprocess}, and give a
high-level overview of the
proof of Theorem~\ref{thmfp}.

\begin{figure}

\includegraphics{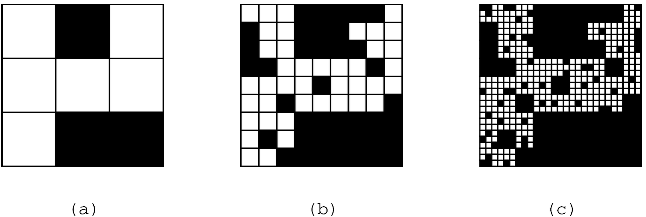}

\caption{Illustration of a fractal percolation process with $\mu=3$
and its $1$-cubes~\textup{(a)}, 2-cubes~\textup{(b)} and
3-cubes~\textup{(c)}. Black squares represent
closed cubes and white squares represent open cubes.}
\label{figfractalpercolation}
\end{figure}

First, we give the definition of standard fractal percolation. Take the
cube $[0,1]^d$ and partition it into $\mu^d$ identical subcubes of
side length
$1/\mu$, where $\mu>1$ is an integer.
We refer to such subcubes as $1$-cubes, and
declare each $1$-cube to be \emph{open} independently with probability
$p\in(0,1)$; otherwise, declare it to be \emph{closed}.
Now repeat this process independently for
each open $1$-cube, splitting it into $\mu^d$ identical subcubes of
size $1/\mu^2$ that we call $2$-cubes and declaring each $2$-cube to
be open independently with probability $p$.
The $1$-cubes that are closed are not partitioned again, and all the
region spanned by these cubes is considered to be closed (see
Figure~\ref{figfractalpercolation}).
We continue this procedure until we obtain $z$-cubes of side length
$1/\mu^{z}$.
Many results in this area~\cite{CCD88,CC89,CCGS91,FG92,O96,MPV01}
study properties of the set of open $z$-cubes
as $z\to\infty$.
%
%

Now, we discuss the intuition behind our definitions and the connection
with fractal percolation.
We start at scale $\kappa$.
We tessellate space and time into very large cells.
These are the cells indexed by the tuples in $\mathcal{R}_\kappa$,
and each cell corresponds to a cube in space and a time interval.
Then, for each cell $(\veci,\tau)\in\mathcal{R}_\kappa$,
we check whether the cell contains sufficiently many nodes at the
beginning of its time interval; that is, we check whether $A_\kappa
(\veci,\tau)=D^\mathrm{ext}_\kappa(\veci,\tau)=1$.
If $A_\kappa(\veci,\tau)=1$, we do a finer tessellation of the cell in
both space and time. In terms of fractal percolation,
this corresponds to the event that a large cube is open and then is
subdivided into smaller cubes. On the other hand,
if $A_\kappa(\veci,\tau)=0$, we skip that cell and tessellate it no
further, similarly to what happens to cubes that are closed in a
fractal percolation process.
We iterate this procedure until we obtain cells of volume $\ell^d\beta$
(i.e., cells of scale $1$).
The main reason for employing this idea
instead of analyzing the events $D_k(\veci,\tau)$ directly is that the
$D_k(\veci,\tau)$ are highly dependent.

In the analysis, we start with the random variables $A_k(\veci,\tau)$
of scale $k=\kappa$,
where cells are so large that we can easily obtain that $A_\kappa
(\veci,\tau)=1$ for
all $(\veci,\tau)\in R_\kappa$.
Then we move from scale $k+1$ to scale $k$. Let $(\veci,\tau)\in
\mathcal
{R}_k$. In order to analyze $A_k(\veci,\tau)$, we need to observe
$A_{k+1}(\veci',\tau')$ such that $\pi_k^{(1)}(\veci)=\veci'$ and
$\gamma_k^{(1)}(\tau)=\tau'$; that is,
$(k+1,\veci',\tau')$ is the parent of~$(k,\veci,\tau)$
with respect to the hierarchies $\pi$~and~$\gamma$.
If $A_{k+1}(\veci',\tau')=0$, then we do not need to observe
$A_k(\veci,\tau)$ since we will not do the finer tessellation
of~$R_{k+1}(\veci
',\tau')$ that
produces the cell $(k,\veci,\tau)$. In this case, we will consider all
the descendants at scale $1$ of the cell $(k+1,\veci',\tau')$ as
``bad,'' and hence we will not need to observe any other
descendant of~$(k+1,\veci',\tau')$ such as $(k,\veci,\tau)$.
On the other hand, if $A_{k+1}(\veci',\tau')=1$, we know that there is
a sufficiently large density of nodes in the region
$S^\mathrm{base}_k(\veci)\subset S^\mathrm{ext}_{k+1}(\veci')$ that
surrounds
$S_k(\veci)$ at time $\tau'\beta_{k+1}$. Then, by allowing these nodes
to move from $\tau'\beta_{k+1}$ to $\tau\beta_k$,
we obtain that many of these nodes move inside $S_k(\veci)$, giving
that the probability that $A_k(\veci,\tau)=0$,
which corresponds to the event $D^\mathrm{base}_{k}(\veci,\tau)=1$
and \mbox{$D^\mathrm{ext}
_{k}(\veci,\tau)=0$},
is small. We then apply this reasoning
for all $(k,\veci,\tau)\in\mathcal{R}$. The key fact is that a dense
cell at scale $k$ makes the children of this cell likely to be dense.

We now give the intuition behind the different types of cubes.
Let $(k,\veci,\tau)\in\mathcal{R}$ and assume that $(k+1,\veci
',\tau
')$ is the parent of~$(k,\veci,\tau)$.
We consider the extended cube $S^\mathrm{ext}_{k+1}(\veci')$ instead
of just
$S_{k+1}(\veci')$ to assure that, when\break 
$D^\mathrm{ext}_{k+1}(\veci',\tau')=1$, then there is a large
density of nodes
around $S_k(\veci)$ at time $\tau'\beta_{k+1}$ even if
$S_k(\veci)$ lies near the boundary of~$S_{k+1}(\veci')$; this
happens since
$\{D^\mathrm{ext}_{k+1}(\veci',\tau')=1\}$ guarantees that there are
sufficiently many nodes in $S^\mathrm{base}_k(\veci)\subset S^\mathrm
{ext}_{k+1}(\veci')$.
We then let the nodes move for time $\tau\beta_k-\tau'\beta
_{k+1}\geq
\beta_{k+1}$,
thereby allowing them to mix in $S^\mathrm{base}_k(\veci)$ and move inside
$S_k(\veci)$.
While these nodes move in the interval $[\tau'\beta_{k+1},\tau\beta
_k]$, they never leave the \emph{area of influence} $S^\mathrm
{inf}_k(\veci)$.
This allows us to argue that cells that are sufficiently far apart in
space are ``roughly independent'' since we only observe nodes
that stay inside the area of influence of their cells.

Now, we give a brief sketch of the proof. We want to give an upper
bound for the probability that $K(0,0)$ is not contained in $\mathcal{R}_1^t$.
When $K(0,0) \nsubseteq\mathcal{R}_1^t$, then there exists a
\emph{very large} path of adjacent bad cells of scale $1$.
A natural strategy is to consider any fixed path of adjacent cells from
the cell $(0,0)$ to
a cell outside $\mathcal{R}_1^t$ and show that the probability that all
cells in this path are bad is exponentially small, and then take the
union bound over all such paths.
However, this strategy seems
challenging
due to the dependencies among the events that the cells of a given path
are bad.
We use two ideas to solve this problem: path of cells of varying scales
and well separated cells.

We start with cells of scale $\kappa$, which are so large that we can
show that, with very large probability, $A_\kappa(\veci,\tau)=1$ for
all $(\veci,\tau)\in\mathcal{R}_\kappa$.
Therefore, if a cell $(\veci,\tau)$ of scale 1 has $A(\veci,\tau)=0$,
we know that there exists an ancestor $(k',\veci',\tau')$
of~$(1,\veci,\tau)$ such that $(k',\veci',\tau')$ is bad but its
parent is good [i.e.,
$A_{k'}(\veci',\tau')=0$].
With this,
we have that if a path of bad cells of scale $1$ exists, then there is
a path of bad cells of varying scales. This path must contain
sufficiently many cells
because it must connect the cell $(0,0)$ to a cell outside $\mathcal{R}_1^t$.
We take any fixed path of cells of varying scale and show that, given
that this path contains sufficiently many cells, we can obtain a
subset of the cells of the path so that
these cells are ``well separated'' in space and time. This implies that
the $A_k(\veci,\tau)$ are ``roughly'' independent for the well
separated cells.
Using this, we can show that
the probability that all cells in this subset are bad is very small.
Then, by applying the union bound over all sets of well separated cells
that can be obtained from a path of
cells of varying scales, we establish
Theorem~\ref{thmfp}.

\subsection{The support of a cell}\label{secsupport}
We define the \emph{time of influence} $T^\mathrm{inf}_k(\tau)$
of~$(k,\tau)$ as
\[
T^\mathrm{inf}_1(\tau) = \bigl[\gamma_1^{(1)}(
\tau)\beta_{2}, \bigl(\tau+\max\{\eta,2\} \bigr)\beta_1
\bigr]
\]
and
\[
T^\mathrm{inf}_k(\tau) = \bigl[ \gamma_k^{(1)}(
\tau)\beta_{k+1},(\tau+2) \beta_k \bigr]\qquad\mbox{for }k
\geq2 %
\]
and set the \emph{region of influence} as
\[
R^\mathrm{inf}_k(\veci,\tau) = S^\mathrm{inf}_k(
\veci)\times T^\mathrm{inf}_k(\tau). %
\]
Intuitively, the event $\{A_k(\veci,\tau)=1\}$ depends on the motion of
nodes during a subinterval of~$T^\mathrm{inf}_k(\tau)$ and these
nodes never
leave the region
$S^\mathrm{inf}_k(\veci)$. Thus, we will be able to argue later that
two cells
with disjoint regions of influence are ``roughly independent.''

We assume that $m$ is sufficiently large with respect to $\eta$ so that\break 
$\max\{\eta,2\} \beta\leq\beta_2=16 m^2\beta$, which gives that
%
%
\begin{equation}
T^\mathrm{inf}_k(\tau)\subseteq T_{k+1} \bigl(
\gamma_{k}^{(1)}(\tau) \bigr) \cup T_{k+1} \bigl(
\gamma_{k}^{(1)}(\tau)+1 \bigr) \cup T_{k+1} \bigl(
\gamma_{k}^{(1)}(\tau)+2 \bigr). \label{eqboundtinf}
\end{equation}

We define the \emph{time support} $T^\mathrm{sup}_k(\tau)$
of~$(k,\tau)$ as
%
%
\begin{eqnarray}
\label{eqtsup} T^\mathrm{sup}_k(\tau) &=& \bigcup
_{i=0}^8T_{k+1} \bigl(\gamma
_{k}^{(1)}(\tau)-3+i \bigr)
\nonumber
\\[-8pt]
\\[-8pt]
\nonumber
&=& \bigl[ \bigl( \gamma_{k}^{(1)}(\tau)-3 \bigr)
\beta_{k+1}, \bigl( \gamma_{k}^{(1)}(\tau)+6
\bigr)\beta_{k+1} \bigr],
\end{eqnarray}
and note that, by~(\ref{eqboundtinf}),
\[
T^\mathrm{inf}_k(\tau) 
\subset T^\mathrm{sup}_k(\tau). %
\]
We also define the \emph{spatial support} $S^\mathrm{sup}_k(\veci)$
of~$(k,\veci)$ as
%
%
\begin{equation}
S^\mathrm{sup}_k(\veci) = \bigcup
_{\veci'\dvtx\llVert \veci'-\pi
_k^{(1)}(\veci
)\rrVert _\infty\leq m} S_{k+1} \bigl(\veci' \bigr),
\label{eqssup}
\end{equation}
and, for any $(k,\veci,\tau)\in\mathcal{R}$, we define
%
%
\begin{equation}
R^{\mathrm{sup}}_k(\veci,\tau) = S^\mathrm{sup}_k(
\veci)\times T^\mathrm{sup}_k(\tau). \label{eqrsup}
\end{equation}
The main idea behind the definition of the support is that, in order to
control dependencies, we will not only need to consider path of cells of
varying scales, but we will also need to restrict our attention to
paths of cells that are sufficiently far apart in both space and time;
we will define
later two cells to be far apart if the support of one cell does not
contain the support of the other.
We now prove some useful geometric properties of the support.

The lemma below gives that if the regions of influence of two cells intersect,
then the region of influence of the cell of smaller scale is contained
in the support of the other cell.
In other words, the support of a given cell contains all the cells of
smaller scale whose region of influence intersects the region of
influence of the given cell.
%

\begin{lemma}\label{lemsupport}
For any sufficiently large $m$ the following is true. 
For any $(k,\veci,\tau),(k',\veci',\tau')\in\mathcal{R}$ with
$k\geq k'$,
if $R^\mathrm{inf}_{k'}(\veci',\tau')\nsubseteq R^{\mathrm
{sup}}_k(\veci,\tau)$ then
$R^\mathrm{inf}_{k'}(\veci',\tau')\cap R^\mathrm{inf}_k(\veci,\tau
)=\varnothing$.
\end{lemma}

\begin{pf}
Note that, if $R^\mathrm{inf}_{k'}(\veci',\tau')\nsubseteq
R^{\mathrm{sup}}_k(\veci,\tau)$, then either
$T^\mathrm{inf}_{k'}(\tau')\nsubseteq\break  T^\mathrm{sup}_k(\tau)$ or
$S^\mathrm{inf}_{k'}(\veci
')\nsubseteq S^\mathrm{sup}_k(\veci)$.
We first assume that $T^\mathrm{inf}_{k'}(\tau')\nsubseteq T^\mathrm
{sup}_k(\tau)$
and show that this implies
\[
T^\mathrm{inf}_{k'} \bigl(\tau' \bigr)\cap
T^\mathrm{inf}_k(\tau)=\varnothing, %
\]
which yields
$R^\mathrm{inf}_{k'}(\veci',\tau')\cap R^\mathrm{inf}_k(\veci,\tau
)=\varnothing$.

Note that the interval $T^\mathrm{inf}_{k'}(\tau')$ has length at
most $3\beta
_{k'+1}$ by~(\ref{eqboundtinf}).
Then, since $T^\mathrm{inf}_{k'}(\tau')\nsubseteq T^\mathrm
{sup}_k(\tau)$,
%
%
\begin{equation}
\label{eqtinfknotkp} T^\mathrm{inf}_{k'} \bigl(\tau'
\bigr) \cap \bigl[ \bigl(\gamma_k^{(1)}(\tau)-3 \bigr)\beta
_{k+1} + 3 \beta_{k'+1}, \bigl(\gamma_{k}^{(1)}(
\tau)+6 \bigr) \beta_{k+1}-3\beta_{k'+1} \bigr]=\varnothing.
\hspace*{-25pt}
\end{equation}
Using the fact that $\beta_{k'} \leq\beta_k$, we obtain
\begin{eqnarray*}
&& \bigl[ \bigl(\gamma_k^{(1)}(\tau)-3 \bigr)
\beta_{k+1} + 3\beta_{k'+1}, \bigl(\gamma_{k}^{(1)}(
\tau)+6 \bigr)\beta_{k+1}-3\beta_{k'+1} \bigr]
\\
&&\qquad\supseteq \bigl[\gamma_k^{(1)}(\tau)
\beta_{k+1}, \bigl(\gamma_{k}^{(1)}(\tau)+3 \bigr)
\beta_{k+1} \bigr]
\\
&&\qquad= T_{k+1} \bigl(\gamma_k^{(1)}(\tau)
\bigr) \cup T_{k+1} \bigl(\gamma_k^{(1)}(\tau)+1
\bigr)\cup T_{k+1} \bigl(\gamma_k^{(1)}(\tau)+2
\bigr)
\\
&&\qquad\supseteq T^\mathrm{inf}_k(\tau),
\end{eqnarray*}
where the last step follows from~(\ref{eqboundtinf}). This, together
with~(\ref{eqtinfknotkp}), implies that $T^\mathrm{inf}_{k'}(\tau
')$ does not
intersect $T^\mathrm{inf}_k(\tau)$.

Now, we turn to the case $S^\mathrm{inf}_{k'}(\veci')\nsubseteq
S^\mathrm{sup}
_k(\veci)$, for which we want to show
\[
S^\mathrm{inf}_{k'} \bigl(\veci' \bigr)\cap
S^\mathrm{inf}_k(\veci)=\varnothing. %
\]
Let $x_1,x_2,\ldots,x_d$ be defined so that $S_k(\veci)=\prod_{j=1}^d
[x_j,x_j+\ell_k]$. Then
we can write
%
%
\begin{equation}
S^\mathrm{inf}_k(\veci) = \prod_{j=1}^d
\bigl[x_j-2\eta mn(k+1)^3\ell_k,x_j+
\ell_k+2\eta mn(k+1)^3\ell_k \bigr].
\label{eqsinfset}
\end{equation}
Now, let $y_1,y_2,\ldots,y_d$ be defined so that $S^\mathrm
{sup}_k(\veci
)=\prod_{j=1}^d [y_j,y_j+(2m+1)\ell_{k+1}]$.
Since $S^\mathrm{inf}_{k'}(\veci')$ is a cube of side length
$(1+4\eta
mn(k'+1)^3)\ell_{k'} \leq(1+4\eta mn(k+1)^3)\ell_{k}$ and
$S^\mathrm{inf}_{k'}(\veci')$ is not contained in $S^\mathrm
{sup}_k(\veci)$, we have that
%
%
\begin{eqnarray}
\label{eqsinf1} && S^\mathrm{inf}_{k'} \bigl(
\veci' \bigr) \cap\prod_{j=1}^d
\bigl[y_j+ \bigl(1+4 \eta mn(k+1)^3 \bigr)\ell
_k,
\nonumber
\\[-8pt]
\\[-8pt]
\nonumber
&&\hspace*{62pt} y_j+(2m+1) \ell_{k+1}-
\bigl(1+4\eta mn(k+1)^3 \bigr)\ell_k \bigr]=\varnothing.
\end{eqnarray}
Now, we use the fact that $m\ell_{k+1} \leq x_j-y_j\leq(m+1)\ell
_{k+1}-\ell_k$ for all $j=1,2,\ldots,d$. This and~(\ref{eqsinfset})
give
%
%
\begin{eqnarray}
\label{eqsinf2} && S^\mathrm{inf}_k(\veci) \subseteq\prod
_{j=1}^d \bigl[y_j+m\ell
_{k+1}-2\eta mn(k+1)^3\ell_k,
\nonumber
\\[-8pt]
\\[-8pt]
\nonumber
&&\hspace*{61pt} y_j+(m+1) \ell_{k+1}+2\eta
mn(k+1)^3\ell_k \bigr].
\end{eqnarray}
Now, using the relation between $m$ and $n$ in~(\ref{eqdefn}), we
have that
%
%
\begin{equation}
\qquad m\ell_{k+1} = m^2(k+1)^3
\ell_k = 7 \eta m n^d (k+1)^3
\ell_k \geq \bigl(1+6\eta mn(k+1)^3 \bigr)
\ell_k. \label{eqspsuprel}
\end{equation}
Using this result in~(\ref{eqsinf1}) we get that $S^\mathrm
{inf}_{k'}(\veci')$
does not intersect
%
%
\begin{equation}
\label{eqreg} \prod_{j=1}^d
\bigl[y_j + \bigl(1+4\eta mn(k+1)^3 \bigr)
\ell_k, y_j+(m+1) \ell_{k+1}+2\eta
mn(k+1)^3\ell_k \bigr].\hspace*{-25pt}
\end{equation}
Similarly, plugging~(\ref{eqspsuprel}) into~(\ref{eqsinf2})
we see that $S^\mathrm{inf}_k(\veci)$ is contained in the space--time region
given by~(\ref{eqreg}).
These two facts
establish that $S^\mathrm{inf}_{k'}(\veci')$ does not intersect
$S^\mathrm{inf}_k(\veci)$.
\end{pf}

The second important property we will use is given in the next lemma,
which establishes that
the support of a cell contains all its descendants. The main use of
this lemma is that, once we encounter a
cell $(k,\veci,\tau)$ for which $A_k(\veci,\tau)=0$, then we want to
regard all its descendants as bad. However, the set of descendants
of~$(k,\veci,\tau)$ may be a complicated set. With the lemma below, we
then just consider all cells that are contained inside the support
of~$(k,\veci,\tau)$ as bad [which includes all descendants
of~$(k,\veci
,\tau)$].
%

\begin{lemma}\label{lemsupport2}
Assume $m\geq3$.
For any $(k,\veci,\tau)\in\mathcal{R}$, if $(k',\veci',\tau')$ is a
descendant of~$(k,\veci,\tau)$ then
\[
R_{k'} \bigl(\veci',\tau' \bigr)
\subseteq R^{\mathrm{sup}}_k(\veci,\tau). %
\]
Moreover, $R^{\mathrm{sup}}_k(\veci,\tau)$ contains all the
neighbors of~$(k',\veci',\tau')$; that is,
\[
\bigcup_{(\veci'',\tau'')\dvtx\llVert (\veci'',\tau'')-(\veci
',\tau
')\rrVert
_\infty\leq1 }R_{k'} \bigl(
\veci'',\tau'' \bigr)
\subseteq R^{\mathrm{sup}}_k(\veci,\tau). %
\]
\end{lemma}

\begin{pf}
Fix $(\veci'',\tau'')$ such that $(k',\veci'',\tau'')$ is adjacent to
$(k',\veci',\tau')$ and assume that the ancestor of~$(k',\veci
'',\tau
'')$ of scale $k$ is
not $(k,\veci,\tau)$, otherwise the second part of the lemma follows
from the first part. We prove this lemma first for space and then for time.
For\vspace*{1pt} space, since $(k',\veci',\tau')$ is a descendant
of~$(k,\veci
,\tau
)$ we have that
$S_{k'}(\veci')\subseteq S_k(\veci)\subseteq S^\mathrm{sup}_k(\veci)$.
Also, $(k',\veci'')$ is adjacent to $(k',\veci')$ which implies that
the ancestor of~$(k',\veci'')$ of
scale $k$ is adjacent to $(k,\veci)$. Since $S^\mathrm{sup}_k(\veci)$
contains all cells of scale $k$ that are adjacent to $(k,\veci)$, it
also contains $S_{k'}(\veci'')$.

It remains to establish the lemma for the time dimension. For the
first part of the lemma, this corresponds to showing that $T_{k'}(\tau
') \subseteq T^\mathrm{sup}_k(\tau)$.
Recall that $T_{k'}(\tau')=[\tau'\beta_{k'},(\tau'+1)\beta_{k'}]$,
which is contained in $[\tau\beta_{k},(\tau'+1)\beta_{k'}]$
since $(k',\veci',\tau')$ is a descendant of~$(k,\veci,\tau)$. Now,
note that
\[
\tau\beta_k=\gamma_{k'}^{(k-k')} \bigl(
\tau' \bigr)\beta_{k} \geq\gamma_{k'}^{(k-k'-1)}
\bigl(\tau' \bigr)\beta_{k-1}-2\beta_k \geq
\tau'\beta_{k'}-2\sum_{j=k'+1}^k
\beta_j. %
\]
Then, since $k'\geq1$, we can use the bound
\[
\sum_{j=2}^k \beta_j =
C_\mathrm{mix}\sum_{j=2}^k
\frac{\ell_{j-1}^2j^4}{\varepsilon^2} \leq C_\mathrm{mix}\frac
{2\ell_{k-1}^2k^4}{\varepsilon^2} = 2
\beta_k, %
\]
where the last inequality can be proved by induction on $k$.
Therefore, we conclude that
%
%
\begin{equation}
\tau\beta_{k}\geq\tau'\beta_{k'}-4
\beta_k. \label{eqbetaxscale}
\end{equation}
Since $k>k'\geq0$,
we have $k\geq1$ and
\[
4\beta_k+\beta_{k'} \leq5\beta_k = 5
\frac{\beta_{k+1}}{m^2k^2(k+1)^4} \leq\beta_{k+1}. %
\]
This and the inequality in~(\ref{eqbetaxscale}) yield
\[
T_{k'} \bigl(\tau' \bigr) \subseteq[\tau
\beta_k, \tau\beta_k+4\beta_k+
\beta_{k'}] \subseteq[\tau\beta_k, \tau\beta_k+
\beta_{k+1}] \subseteq T^\mathrm{sup} _k(\tau).
\]
This establishes the first part of the lemma. For the second part,
using the fact that $(k',\tau'')$ is adjacent to $(k',\tau')$ together
with the result above, we have
\[
T_{k'} \bigl(\tau'' \bigr) \subseteq[\tau
\beta_k-\beta_{k'}, \tau\beta_k+
\beta_{k+1}+\beta_{k'}] \subseteq T^\mathrm{sup}_k(
\tau). %
\]
\upqed
\end{pf}

\subsection{Support connected paths}\label{secscpath}
We start defining the \emph{extended support} of a cell.
Given a cell $(k,\veci,\tau)\in\mathcal{R}$, define
\[
T^\mathrm{2sup}_k(\tau) = \bigcup
_{i=0}^{26}T_{k+1} \bigl(\gamma
_{k}^{(1)}(\tau)-12+i \bigr)
\]
and
\[
S^\mathrm{2sup}_k(\veci) = \bigcup
_{\veci'\dvtx\llVert \veci'-\pi
_k^{(1)}(\veci
)\rrVert _\infty\leq3m+1} S_{k+1} \bigl(\veci' \bigr).
\]
Also, we let
\[
R^{\mathrm{2sup}}_k(\veci,\tau) = S^\mathrm{2sup}_k(
\veci)\times T^\mathrm{2sup}_k(\tau). %
\]
The extended support is defined so that the following property is satisfied.
Let $(k_1,\veci_1,\tau_1),(k_2,\veci_2,\tau_2)\in\mathcal{R}$ with
$k_1\geq k_2$. Then
%
%
\begin{equation}
\label{eqpropsupport} \mbox{if } R^{\mathrm{sup}}_{k_1}(\veci_1,
\tau_1)\mbox{ intersects } R^{\mathrm{sup}} _{k_2}(
\veci_2,\tau_2),\mbox{ we have } R^{\mathrm{sup}}_{k_2}(
\veci_2,\tau_2) \subseteq R^{\mathrm
{2sup}}_{k_1}(
\veci_1,\tau_1).\hspace*{-25pt}
\end{equation}
Note that the extended support of~$(k,\veci,\tau)$ is three times
larger than the support of~$(k,\veci,\tau)$ [defined in~(\ref
{eqrsup})] since
$S^\mathrm{sup}_k(\veci)$ is a cube of
side length $(2m+1)\ell_{k+1}$ while
$S^\mathrm{2sup}_k(\veci)$ is a cube of
side length $3(2m+1)\ell_{k+1}$, and
$T^\mathrm{sup}_k(\tau)$ is an interval of length
$9\beta_{k+1}$ while
$T^\mathrm{2sup}_k(\tau)$ is an interval of length
$27\beta_{k+1}$.

We define the extended support because on the one hand we want to look
at cells that are well separated (in the sense that the support of one cell
does not contain the support of the other), but on the other hand we
want to have a notion of a path of well separated cells. In such a
path, cells must be
well separated but should not be excessively far from one another. We
will use the extended support to say that two cells are adjacent (in
this new
notion of path) if their extended supports intersect; we will call this
a \emph{support connected path}. We make this notion rigorous in the following.

Recall that a cell $(\veci,\tau)\in\mathcal{R}_1$ is said to have a bad
ancestry if $A(\veci,\tau)=0$.
Also, in the beginning of Section~\ref{secfp}, we defined a cell
$(\veci,\tau)$ of scale 1 to be bad if $E(\veci,\tau)=0$.
Here, we change this definition slightly and extend it to cells of
arbitrary scales:
we say that a cell $(k,\veci,\tau)\in\mathcal{R}$ is bad if
$A_k(\veci,\tau)=0$.
Our goal is to show that, if $t$ is sufficiently large, then the
probability that $K(\origin,0)\nsubseteq\mathcal{R}^t_1$ is small.
First, recall that two cells $(k,\veci_1,\tau_1)$ and $(k,\veci
_2,\tau
_2)$ at the same scale are adjacent if
$\llVert \veci-\veci'\rrVert _\infty\leq1$ and $\vert \tau-\tau
'\vert \leq1$.
For arbitrary scales $k_1>k_2$, we define that
\[
\begin{tabular} {p{310pt}} $(k_1,\veci_1,
\tau_1)$ and\vspace*{2pt} $(k_2,\veci_2,
\tau_2)$ are adjacent if $(k_1,\veci_1,
\tau_1)$ is adjacent to $ (k_1,
\pi_{k_2}^{(k_1-k_2)}( \veci_2), \gamma_{k_2}^{(k_1-k_2)}(
\tau_2) )$. \end{tabular} %
\]
In other words, $(k_1,\veci_1,\tau_1)$ and $(k_2,\veci_2,\tau_2)$ are
adjacent if the cell at scale
$k_1$ that is the ancestor of~$(k_2,\veci_2,\tau_2)$ is adjacent to
$(k_1,\veci_1,\tau_1)$.
Note that
two adjacent cells of different scales may be disjoint, whereas two
adjacent cells of the same scale must at
least intersect at a point.

We refer to a \emph{path} as a sequence of distinct cells for which any
two consecutive cells in the sequence are adjacent,
and we say that a sequence of cells is
a \emph{cluster} if each cell of the sequence is adjacent
to \emph{some} other cell in the sequence. Note that, unlike in a path,
the order of the cells of
a cluster is not important, so
we regard a cluster as a \emph{set} of cells.

For any two cells $(k_1,\veci_1,\tau_1)$ and $(k_2,\veci_2,\tau
_2)$, we
say that
%
%
\begin{equation}
\label{eqwellseparated} %
\begin{tabular} {p{245pt}} $(k_1,
\veci_1,\tau_1)$\vspace*{2pt} and $(k_2,
\veci_2,\tau_2)$ are \emph{well separated}~if\\
$R_{k_1}(\veci_1, \tau_1) \nsubseteq
R^{\mathrm{sup}}_{k_2}(\veci_2,\tau_2)$ and
$R_{k_2}(\veci_2,\tau_2) \nsubseteq
R^{\mathrm{sup}}_{k_1}(\veci_1, \tau_1)$.
\end{tabular} %
\end{equation}
We will look at path of cells that are mutually well separated. (Recall
that from Lemma~\ref{lemsupport}, well separated cells have disjoint
regions of influence and,
therefore, we will be able to argue that they behave roughly
independently from one another.)
We define that
\[
\begin{tabular} {p{230pt}} $(k_1,\veci_1,
\tau_1)$\vspace*{2pt} and $(k_2,\veci_2,\tau
_2)$ are \emph{support adjacent} if $R^{\mathrm{2sup}}_{k_1}(
\veci_1,\tau_1)\cap R^{\mathrm
{2sup}}_{k_2}(
\veci_2,\tau_2) \neq\varnothing$. \end{tabular}
\]
Finally, we say that a sequence of cells $P=((k_1,\veci_1,\tau
_1),(k_2,\veci_2,\tau_2),\ldots,\break (k_z,\veci_z,\tau_z))$ is
a \emph{support connected path} if the cells in $P$ are mutually well
separated and, for each $j=1,2,\ldots, z-1$,
$(k_j,\veci_j,\tau_j)$ is support adjacent to $(k_{j+1},\veci
_{j+1},\tau
_{j+1})$.
We also define a sequence of cells $P=((k_1,\veci_1,\tau
_1),(k_2,\veci
_2,\tau_2),\break \ldots,(k_z,\veci_z,\tau_z))$ to be
a \emph{support connected cluster} if the cells in $P$ are mutually
well separated and, for each $j=1,2,\ldots, z$,
there exists a $j'\in\{1,2,\ldots,z\}\setminus\{j\}$ such that
$(k_j,\veci_j,\tau_j)$ is support adjacent to $(k_{j'},\veci
_{j'},\tau_{j'})$.

Now, define
$\Omega$ as the set of all paths of cells of scale $1$ (i.e., cells
of~$\mathcal{R}_1$) so that the first cell of the path is $(\origin
,0)$ and
the last cell of the path is the only cell not contained in $\mathcal{R}^t_1$.
Also, define $\Omega^\mathrm{sup}_\kappa$ as the set of all support connected
paths of cells of scale at most
$\kappa$ (i.e., cells in $\bigcup_{k=1}^\kappa\mathcal{R}_k$) so
that the
extended support of the
first cell of the path contains $R_1(\origin,0)$ and the last cell of
the path is the only cell whose
extended support is not contained in $\bigcup_{(\veci,\tau)\in
\mathcal
{R}^t_1}R_1(\veci,\tau)$. The lemma below will allow us to turn our attention
to support connected paths of bad cells instead of
paths of cells with bad ancestry, whose
dependencies seem challenging to control.

%

\begin{lemma}\label{lemsupportpath}
We have that
\begin{eqnarray*}
&& \mathbf{P} (\exists P\in\Omega\mbox{ s.t. all cells of~$P$ have a bad
ancestry} )
\\
&&\qquad\leq\mathbf{P} \bigl(\exists P\in\Omega^\mathrm
{sup}_\kappa \mbox{ s.t. all cells of~$P$ are bad} \bigr). %
\end{eqnarray*}
\end{lemma}

\begin{pf}
The proof is split into two stages. In the first stage,
we show that, if there exists a path $P\in\Omega$ such that each cell
of~$P$ has a bad ancestry,
then there exists a path of bad cells of
arbitrary scales. In the second stage, we show that, given the
existence of such a path of bad cells of arbitrary scales,
then there exists a path of~$\Omega^\mathrm{sup}_\kappa$ such that
all cells
of the path are bad.

We now prove the first stage.
Let $\Omega_\kappa$ be the set of all paths of cells of arbitrary
scale (i.e., cells in $\mathcal{R}$) such that the first cell of the
path is an ancestor
of~$(\origin,0)\in\mathcal{R}_1$
and the last cell of the path is the only cell
whose support is not contained in $\bigcup_{(\veci,\tau)\in\mathcal
{R}^t_1}R_1(\veci,\tau)$.
In this stage, we establish that
\begin{eqnarray}
\label{eqstage1} && \mathbf{P} (\exists P\in\Omega\mbox{ s.t. all cells of~$P$
have a bad ancestry} )
\nonumber
\\[-8pt]
\\[-8pt]
\nonumber
&&\qquad\leq\mathbf{P} (\exists P\in\Omega_\kappa\mbox{ s.t.
all cells of~$P$ are bad} ).
\nonumber
\end{eqnarray}
Let $P=((1,\veci_1,\tau_1),(1,\veci_2,\tau_2),\ldots,(1,\veci
_z,\tau
_z))\in\Omega$ be a path of cells with bad ancestries; hence $(\veci
_1,\tau_1)=(\origin,0)$
and $(\veci_z,\tau_z)\notin\mathcal{R}^t_1$.
For each $j$, since $A(\veci_j,\tau_j)=0$, we
know by the definition of~$A$ in~(\ref{eqa})
that there\vspace*{-2pt} exists a $k_j'$ so that, if we set $\veci_j'
= \pi
_{1}^{(k_j'-1)}(\veci_j)$ and
$\tau_j' = \gamma_{1}^{(k_j'-1)}(\tau_j)$, we obtain $A_{k_j'}(\veci
_j',\tau_j')=0$.
Define $J\subseteq\{1,2,\ldots,z\}$ such that
$j\in J$ if and only if there exists no $j'<j$ with $(k_j',\veci
_j',\tau_j')=(k_{j'}',\veci_{j'}',\tau_{j'}')$ and there exists no
$j' \in\{1,2,\ldots,z\}\setminus\{j\}$ for which $(k_j',\veci
_j',\tau
_j')$ is a descendant of~$(k_{j'}',\veci_{j'}',\tau_{j'}')$.
In other words, $J$ contains only distinct elements of the set $\{
(k'_j,\veci_j',\tau_j') \dvtx j=1,2,\ldots,z\}$ which have
no ancestor in the set.
With this, we define
\[
\tilde P = \bigl\{ \bigl(k_j',\veci_j',
\tau_j' \bigr) \dvtx j\in J \bigr\}, %
\]
and show that $\tilde P$ is a cluster.
Before, note that,
since each cluster contains a path, this establishes the existence of
a path of bad cells of arbitrary scales.
In particular, we obtain a path starting from an ancestor of~$(1,\veci
_1,\tau_1)$ and such that there exists a cell $(k',\veci',\tau')\in
\tilde P$
that is an ancestor of a cell of~$P$ that is not contained in~$\mathcal
{R}^t_1$.
Then, by Lemma~\ref{lemsupport2}, we know that this cell of~$P$ is contained
in $R^{\mathrm{sup}}_{k'}(\veci',\tau')$, which gives that the union
of the
support of the
cells in $\tilde P$ is not contained in $\bigcup_{(\veci,\tau)\in
\mathcal
{R}^t_1}R_1(\veci,\tau)$.
Such a path belongs to $\Omega_\kappa$,
so it only remains to show that $\tilde P$ is a cluster.
Note that, by construction, each cell of~$P$ has exactly one ancestor
in $\tilde P$.
Now, for any two adjacent cells $(1,\veci_j,\tau_j),(1,\veci
_{j+1},\tau
_{j+1})$ of~$P$, either they have the
same ancestor in $\tilde P$ or their ancestors are adjacent since two
nonadjacent cells cannot
have descendants at scale $1$ that are adjacent.
This shows that each cell in $\tilde P$ is adjacent to at least one
other cell in $\tilde P$ and, consequently,
$\tilde P$ is a cluster. Therefore, we obtain~(\ref{eqstage1}).

Now we turn to the second stage of the proof, where we establish that
%
%
\begin{eqnarray}
\label{eqstage2} && \mathbf{P} (\exists P\in\Omega_\kappa\mbox{ s.t. all
cells of~$P$ are bad} )
\nonumber
\\[-8pt]
\\[-8pt]
\nonumber
&&\qquad\leq\mathbf{P} \bigl(\exists P\in\Omega^\mathrm
{sup}_\kappa\mbox{ s.t. all cells of~$P$ are bad} \bigr).
\end{eqnarray}
Let $P=((k_1,\veci_1,\tau_1),(k_2,\veci_2,\tau_2),\ldots
,(k_z,\veci
_z,\tau_z))\in\Omega_\kappa$ be a path of bad cells; thus,
the support of~$(k_z,\veci_z,\tau_z)$ is not contained in $\bigcup_{(\veci
,\tau)\in\mathcal{R}^t_1}R_1(\veci,\tau)$.
We show the existence of a
support connected cluster $P'$ of bad cells.
First, we order the cells of~$P$ in the following way. If two cells
have the same scale, we order them by taking an arbitrary order
of~$\mathbb{Z}^{d+1}$; for two cells of
different scales, we say that the cell with the larger scale precedes
the other cell in the order.
This clearly establishes a total order of
the cells of~$P$. Then let $L$ be the list of cells of~$P$ following
this order, where the first cell of~$L$ is the cell that precedes
all the other cells of~$P$ in the order.
We construct $P'$ in a step-by-step manner, where at each step we add
the first element of~$L$ to $P'$,
remove some elements from $L$ and repeat until $L$ has no element.
During this procedure, we associate each cell of~$P$ to a cell of~$P'$;
we use this association
later to show that $P'$ is a support connected cluster.
Below we give the formal description of each step in the construction
of~$P'$, where we assume that
$(k',\veci',\tau')$ is the first element of~$L$:
\begin{longlist}[(ii)]
\item Add $(k',\veci',\tau')$ to $P'$ and remove it from $L$. Since
$(k',\veci',\tau')$ is both in $P$ and $P'$, associate $(k',\veci
',\tau
')$ to itself.
\item Remove from $L$ all the cells $(k'',\veci'',\tau'')$ that are not
well separated from $(k',\veci',\tau')$
[see the definition of well separated cells in~(\ref
{eqwellseparated})], and associate each such
$(k'',\veci'',\tau'')$ to $(k',\veci',\tau')$.
\end{longlist}
We repeat these steps until $L$ is empty.
Note that the set $P'$ obtained at the end contains only cells that
are mutually well separated.
Also, note that there exists a cell of~$P'$ such that the extended
support of this cell contains $R_1(\origin,0)$.
This follows since, by definition, $R^{\mathrm{sup}}_{k_1}(\veci
_1,\tau_1)$
contains $R_1(\origin,0)$ and, by construction of~$P'$, there exists a
cell of~$P'$ whose\vspace*{1pt} support contains $R_{k_1}(\veci
_1,\tau_1)$. So, by~(\ref
{eqpropsupport}), the extended support of this cell contains
$R^{\mathrm{sup}}
_{k_1}(\veci_1,\tau_1)$ which
contains $R_1(\origin,0)$.

It remains to show that $P'$ is support connected
and that
\[
\bigcup_{(k',\veci',\tau')\in P'}R^{\mathrm{2sup}}_{k'}
\bigl( \veci',\tau' \bigr) \nsubseteq\bigcup
_{(\veci'',\tau'')\in\mathcal{R}^t_1}R_1 \bigl(\veci'',
\tau'' \bigr). %
\]
The second property is easy to check and follows from~(\ref
{eqpropsupport}) by noting that the cell $(k_z,\veci_z,\tau_z)\in P$ is
contained
in the support of the cell to which it has been associated in the
construction of~$P'$ and,
moreover, the support of~$(k_z,\veci_z,\tau_z)$ is not contained in
$\bigcup_{(\veci'',\tau'')\in\mathcal{R}^t_1}R_1(\veci'',\tau'')$
by the definition of~$P$.
Then~(\ref{eqpropsupport}) gives that the extended support of the
cell to which $(k_z,\veci_z,\tau_z)$ has been associated contains the support
of~$(k_z,\veci_z,\tau_z)$ and, hence, cannot be contained in $\bigcup_{(\veci'',\tau'')\in\mathcal{R}^t_1}R_1(\veci'',\tau'')$.
Here, we used the order of~$L$, which guarantees that $k_z$ is no
larger than the scale of the cell
to which $(k_z,\veci_z,\tau_z)$ has been associated, thereby allowing
us to apply~(\ref{eqpropsupport}).

%
%
\begin{figure}

\includegraphics{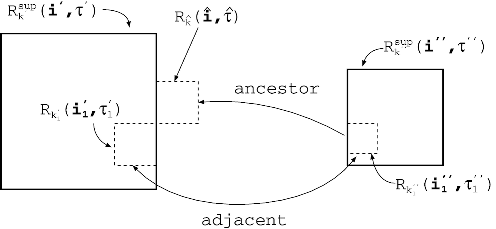}

\caption{Illustration for the proof of Lemma~\protect\ref{lemsupportpath}.
It shows that $R^{\mathrm{2sup}}_{k'}(\veci',\tau')$ intersects
$R^{\mathrm{2sup}}
_{k''}(\veci'',\tau'')$ since
$R^{\mathrm{sup}}_{\hat k}(\hat\veci,\hat\tau)$ intersects both
$R^{\mathrm{sup}}
_{k'}(\veci',\tau')$ and $R^{\mathrm{sup}}_{k''}(\veci'',\tau'')$.}\label{figsupportadjacent}
\end{figure}

Now, we assume, for the sake of establishing a contradiction, that
$P'$ is not support connected. Then
it must be the case that
$P'$ can be partitioned into two set of cells $Q$ and $Q'$ so that,
for any cell of~$Q$, the extended support of this cell
does not intersect the extended support of any cell of~$Q'$.
In this part, it will be useful to refer to Figure~\ref{figsupportadjacent}.
Let $P_0$ be the cells of~$P$ that are \emph{not} associated to any
cell of~$Q$ but are adjacent to at least one cell of~$P$ that is
associated to
a cell of~$Q$.
Since $P$ is a set of adjacent cells, $P_0$ has at least one cell.
Let $(k''_1,\veci''_1,\tau''_1)$ be a cell in $P_0$ and let
$(k'',\veci
'',\tau'')$ be the cell of~$P'$ that has been associated to
$(k''_1,\veci''_1,\tau''_1)$.
So $(k'',\veci'',\tau'')\in Q'$.
Let $(k'_1,\veci'_1,\tau'_1)$ be a cell of~$P$ that has been
associated to $(k',\veci',\tau')\in Q$ and is adjacent
to $(k''_1,\veci''_1,\tau''_1)$. Thus, $R_{k'_1}(\veci'_1,\tau'_1)
\subseteq R^{\mathrm{sup}}_{k'}(\veci',\tau')$ and
$R_{k''_1}(\veci''_1,\tau''_1) \subseteq R^{\mathrm
{sup}}_{k''}(\veci'',\tau'')$.
We now assume that $k'_1 \geq k''_1$; the other case follows by the
same argument.
Since $k'_1\geq k''_1$, we have a cell $(\hat k, \hat\veci,\hat\tau
)$ at scale
$\hat k=k'_1$ that is adjacent to $(k'_1,\veci'_1,\tau'_1)$ and is an
ancestor of~$(k''_1,\veci''_1,\tau''_1)$.
Then, by adjacency and Lemma~\ref{lemsupport2},
$R^{\mathrm{sup}}_{\hat k}(\hat\veci,\hat\tau)$ contains both
$R_{k'_1}(\veci'_1,\tau'_1)$ and $R_{k''_1}(\veci''_1,\tau''_1)$. Therefore,
$R^{\mathrm{sup}}_{\hat k}(\hat\veci,\hat\tau)$ intersects both
$R^{\mathrm{sup}}_{k'}(\veci',\tau')$ and $R^{\mathrm
{sup}}_{k''}(\veci'',\tau'')$.
Also, by the order of~$L$, we have that
$\hat k \leq k'$ since $\hat k=k'_1$ and $(k'_1,\veci'_1,\tau'_1)$ is
associated to $(k',\veci',\tau')$. Then,
by~(\ref{eqpropsupport}), we have that $R^{\mathrm{sup}}_{\hat
k}(\hat\veci
,\hat\tau)\subseteq R^{\mathrm{2sup}}_{k'}(\veci',\tau')$ which
gives that
$R^{\mathrm{2sup}}_{k'}(\veci',\tau')$ intersects $R^{\mathrm
{sup}}_{k''}(\veci
'',\tau'')$, thereby contradicting the fact that
$(k',\veci',\tau')$ and $(k'',\veci'',\tau'')$ are not support
adjacent. This establishes~(\ref{eqstage2}) and completes the proof of
the lemma.
\end{pf}

The next lemma is a technical result bounding the probability that a
Brownian motion stays inside a cube.
%

\begin{lemma}\label{lemeventF}
Let $\Delta>0$ and, for any $z>0$, define $F_\Delta(z)$ to be the
event that a Brownian motion starting from the origin stays inside
$Q_z$ throughout the time interval $[0,\Delta]$. Then, for any $z\geq
3\sqrt{\Delta}$, we have
\[
\mathbf{P} \bigl(F_{\Delta}(z) \bigr) \geq1-d\exp \biggl(-
\frac{z^2}{18\Delta} \biggr). %
\]
\end{lemma}

\begin{pf}
In order to bound $\mathbf{P} (F_{\Delta
}(z) )$, we use the bound for $f_\Delta
$ in~(\ref{eqdeff}) with $M=z/3$,
where $\frac{f_\Delta}{\mathbf{P} (F_{\Delta
}(z) )}$ is
the probability density function of the position of a Brownian motion
at time $\Delta$ given that the motion never leaves $Q_{z}$ during the
whole of~$[0,\Delta]$.
With this, we have
\begin{eqnarray*}
\mathbf{P} \bigl(F_{\Delta}(z) \bigr) &=& \int_{Q_{z}}
f_{\Delta}(y) \,dy \geq\int_{Q_{2z/3}} f_{\Delta}(y)
\,dy
\\
&\geq& \biggl(1-2d\exp \biggl(-\frac{z^2}{6\Delta} \biggr) \biggr) \biggl(1-2d
\int_{z/3}^\infty\frac{1}{\sqrt{2\pi\Delta}}\exp \biggl(-
\frac{y_1^2}{2\Delta} \biggr) \,dy_1 \biggr)
\\
&\geq& \biggl(1-2d\exp \biggl(-\frac{z^2}{6\Delta} \biggr) \biggr) \biggl(1-2d
\frac{3\sqrt{\Delta}}{\sqrt{2\pi}z}\exp \biggl(-\frac
{z^2}{18\Delta
} \biggr) \biggr)
\\
&\geq& 1-d\exp \biggl(-\frac{z^2}{18\Delta} \biggr),
\end{eqnarray*}
where we use the fact that $z\geq3\sqrt{\Delta}$ to apply the Gaussian
tail bound in the second inequality (cf. Lemma~\ref{lemgaussiantail})
and also
to obtain the simplifications in the last inequality.
\end{pf}

We now give a key lemma that we will use to argue that well separated
cells are roughly independent.
Let $\mathcal{F}_k(\veci,\tau)$ be the $\sigma$-field generated by all
$A_{k'}(\veci',\tau')$ for which
$T^\mathrm{inf}_{k'}(\tau')$ does not intersect $[\gamma
_{k}^{(1)}(\tau)\beta
_{k+1},\infty)$ or
both $\tau'\beta_{k'} \leq\tau\beta_k$ and $S^\mathrm
{inf}_k(\veci)\cap S^\mathrm{inf}
_{k'}(\veci')=\varnothing$.
An\vspace*{1pt} important property of this definition is that
an event in $\mathcal{F}_k(\veci,\tau)$ may reveal information about
nodes that affect $A_k(\veci,\tau)$ but only
regarding time steps that occur \emph{before} $\gamma_k^{(1)}(\tau
)\beta_{k+1}$.
Because of this, we are able to get an upper bound for the probability
that $A_k(\veci,\tau)=0$ given any event in $\mathcal{F}_k(\veci
,\tau)$.
The following quantity will be used in many of the subsequent lemmas to
simplify the equations
%
%
\begin{eqnarray} \label{eqpsi}
\psi_1 &=& \min \biggl\{\varepsilon^2\lambda
\ell^d,\log \biggl(\frac
{1}{1-\nu
_E((1-\varepsilon)\lambda,Q_{w\ell})} \biggr) \biggr\}\quad\mbox{and}
\nonumber\\[-8pt]\\[-8pt]\nonumber
\psi_k&=& \frac{\varepsilon^2\lambda\ell_{k-1}^d}{(k+1)^4} = \frac
{\varepsilon^2\lambda\ell^d
m^{d(k-2)}((k-1)!)^{3d}}{(k+1)^4}\qquad\mbox{for }k
\geq2.
\end{eqnarray}

%
%
\begin{lemma}\label{lempast}
Let $w\geq\sqrt{\frac{18\eta\beta}{\ell^2}\log(\frac
{8d}{\varepsilon
} )}$ and
\[
\alpha=\min\biggl\{\varepsilon^2\lambda\ell^d,\log
\biggl(\frac{1}{1-\nu_E((1-\varepsilon)\lambda,Q_{w\ell})} \biggr) \biggr\}
\]
as in
Theorem~\ref{thmfp}.
Fix any $(k,\veci, \tau)\in\mathcal{R}$ and any $F\in\mathcal
{F}_{k}(\veci,\tau)$.
If $m$ is sufficiently large with respect to $d$, $\beta/\ell^2$,
$\eta
$ and $\varepsilon$, then
there are positive constants \mbox{$c=c(d)\geq1$} and $\alpha_0$ so that,
for all $\alpha\geq\alpha_0$,
we have:
\begin{longlist}[(ii)]
\item $\mathbf{P} (A_k(
\veci,\tau)=0 ) \leq\exp(-c\psi_k )$ for all $k=1,2,
\ldots,\kappa$,

\item $\mathbf{P} (A_k(\veci,\tau)=0 \mid F ) \leq
\exp(-c\psi _k )$ for all $k=1,2,\ldots,\kappa-1$.
\end{longlist} %
\end{lemma}

\begin{pf}
Note that the $A_k$ are defined differently for $k=1$ and $2\leq k
\leq\kappa-1$.
In the sequel, we assume that $k\geq2$ and establish part~(ii) of the
lemma. At the end, we address both part~(i) and the case $k=1$ for part~(ii).
Since
\[
\mathbf{P} \bigl(A_k(\veci,\tau)=0 \mid F \bigr) =\mathbf{P} \bigl(
\bigl\{D^\mathrm{ext}_k( \veci,\tau)=0 \bigr\}\cap \bigl\{
D^\mathrm{base}_k(\veci, \tau)=1 \bigr\} \mid F \bigr),
\]
if $F \cap\{ D^\mathrm{base}_k(\veci,\tau)=1\}=\varnothing$, then
the lemma
clearly holds. So, we now assume that
$F \cap\{ D^\mathrm{base}_k(\veci,\tau)=1\}\neq\varnothing$ and write
\[
\mathbf{P} \bigl(A_k(\veci,\tau)=0 \mid F \bigr) \leq\mathbf {P}
\bigl(D^\mathrm{ext}_k( \veci,\tau)=0 \mid F\cap \bigl\{
D^\mathrm{base}_k(\veci,\tau)=1 \bigr\} \bigr). %
\]
Recall that $\{D^\mathrm{base}_k(\veci,\tau)=1\}$ gives that all
cubes of scale
$k$ contained
in~$S^\mathrm{base}_k(\veci)$ have at least $(1-\varepsilon
_{k+1})\lambda\ell_k^d$
nodes at\vspace*{1pt} time $\gamma_k^{(1)}(\tau)\beta_{k+1}$ and
the displacement of
these nodes throughout $[\gamma_k^{(1)}(\tau)\beta_{k+1},\tau\beta_k]$
is in $Q_{\eta m n(k+1)^3\ell_k}$.
Note that $F$
only reveals information about the location of these nodes \emph
{before} time $\gamma_{k}^{(1)}(\tau)\beta_{k+1}$ since these nodes
never leave
the cube $S^\mathrm{inf}_k(\veci)$ during the whole of~$[\gamma
_k^{(1)}(\tau
)\beta_{k+1},\tau\beta_k]$ (cf. Remark~\ref{remdextfillbase}).

Now, we apply Proposition~\ref{procouplingapplied}. To avoid
ambiguity, we add a bar to the variables appearing in the statement of
Proposition~\ref{procouplingapplied}. We apply this proposition with
$\bar K=(1+2\eta mn(k+1)^3)\ell_k$, $\bar\ell= \ell_k$, $\bar\beta
=(1-\varepsilon_{k+1})\lambda$,
$\bar\Delta=\tau\beta_k-\gamma_k^{(1)}(\tau)\beta_{k+1}\in
[\beta
_{k+1},2\beta_{k+1}]$,
$\bar K'$ such that $3(\bar K-\bar K'+2\sqrt{d}\bar\ell)=\eta m
n(k+1)^3\ell_k$
and $\bar\varepsilon$ such that
$(1-\bar\varepsilon)(1-\varepsilon_{k+1})= (1-\frac{\varepsilon
_{k+1}+\varepsilon_k}{2} )$, which gives that
$\bar\varepsilon\geq\frac{\varepsilon_k-\varepsilon
_{k+1}}{2}=\frac
{\varepsilon
}{2(k+1)^2}$. Now, using these values and the fact that $m$ is large
enough, we have that
\[
\bar K' = \ell_k + mk^3\ell_{k-1}
\bigl(\tfrac{5}{3}\eta m n (k+1)^3 + 2\sqrt{d} \bigr)\geq
\ell_k + 2 \eta m nk^3\ell_{k-1}, %
\]
which\vspace*{1pt} is the side length of~$S^\mathrm{ext}_k(\veci)$.
Note also that we have $\bar\Delta\geq\frac{c_1\bar\ell^2}{\bar
\varepsilon^2}$ since $C_\mathrm{mix}\geq4c_1$ in the definition
of~$\beta_{k+1}$.
It remains to check whether $\bar K-\bar K'\geq c_2\sqrt{\bar\Delta
\log(16d\bar\varepsilon^{-1})}$, which is satisfied if the following
is true:
%
%
\begin{equation}
\eta mn(k+1)^3\ell_k \geq4c_2
\sqrt{2C_\mathrm{mix}}\frac{\ell_k(k+1)^2}{\varepsilon} \sqrt {\log \biggl(
\frac{32d(k+1)^2}{\varepsilon} \biggr)}. \label{eqcondm}
\end{equation}
Using the definition of~$\beta_k$ from~(\ref{eqdefbeta}) for $k=1$, we
can write
$C_\mathrm{mix}=\frac{\varepsilon^2m^2\beta}{\ell^2}$, which allows
us to write
the right-hand side of~(\ref{eqcondm}) as
\[
4\sqrt{2}c_2 m \ell_k (k+1)^2 \sqrt{
\frac{\beta}{\ell^2}\log \biggl(\frac
{32d(k+1)^2}{\varepsilon} \biggr)}. %
\]
Note that, in the left-hand side of~(\ref{eqcondm}),
$\eta n$ increases with $m$. So, since $m$ is sufficiently large with
respect to $d$, $\beta/\ell^2$ and $\varepsilon$,
we obtain that $n$ is also sufficiently large
and~(\ref{eqcondm}) is satisfied for all $k$.

Then, we obtain a coupling between the nodes in $S^\mathrm
{base}_k(\veci)$ and
an independent Poisson point process $\Xi$
with intensity
$(1-\bar\varepsilon)(1-\varepsilon_{k+1})\lambda\geq(1-\frac
{\varepsilon
_k}{2}-\frac{\varepsilon_{k+1}}{2} )\lambda$
that succeeds with probability at least
\begin{eqnarray*}
&& 1- \frac{\bar K^d}{\bar\ell^d}\exp \bigl(-c_3\bar\varepsilon^2
\bar\beta\bar\ell^d \bigr)
\\
&&\qquad\geq1- \bigl(1+2\eta mn (k+1)^3 \bigr)^d\exp
\biggl(-c_3\frac{\varepsilon
^2}{4(k+1)^4}(1- \varepsilon_{k+1})\lambda
\ell_k^d \biggr),
\end{eqnarray*}
where $c_3$ is a constant depending on $d$ only. Note that, up to this
moment, we never used the fact that $k\geq2$ and the argument above
holds also
for $k=1$.

Now, for the case $k\geq2$, we define a Poisson point process $\Xi'$
consisting of the nodes of~$\Xi$ whose displacement throughout $[\tau
\beta_k,(\tau+2)\beta_k]$
is in $Q_{\eta m n k^3\ell_{k-1}}$. For each node of~$\Xi$, this
condition is satisfied with probability
$\mathbf{P} (F_{2\beta_k}(\eta m n k^3\ell
_{k-1}) )$,
independently over the nodes of~$\Xi$. Using Lemma~\ref{lemeventF}
and the thinning property of Poisson processes, we have that $\Xi'$ is
a Poisson point
process with intensity
%
%
\begin{eqnarray}
\label{eqpastintensity1} && (1-\bar\varepsilon) (1-\varepsilon _{k+1} )
\mathbf{P} \bigl(F_{2\beta_k} \bigl(\eta m nk^3\ell
_{k-1} \bigr) \bigr)\lambda
\nonumber
\\
&&\qquad\geq \biggl(1-\frac{\varepsilon_k}{2}-\frac{\varepsilon
_{k+1}}{2} \biggr) \biggl(1-d
\exp \biggl(-\frac{(\eta m nk^3\ell_{k-1})^2}{36\beta
_k} \biggr) \biggr)
\nonumber
\\[-8pt]
\\[-8pt]
\nonumber
&&\qquad\geq \biggl(1-\frac{\varepsilon_k}{2}-\frac{\varepsilon
_{k+1}}{2} \biggr)
\biggl(1-d\exp \biggl(-\frac{(\eta m n k\varepsilon)^2}{36C_\mathrm
{mix}} \biggr) \biggr)
\\
&&\qquad\geq \biggl(1-\frac{\varepsilon_k}{2}-\frac{\varepsilon
_{k+1}}{2} \biggr) \biggl(1-d
\exp \biggl(-\frac{(\eta n k)^2}{36(\beta/\ell^2)} \biggr) \biggr),
\nonumber
\end{eqnarray}
where the second inequality follows by the definition of~$\beta_k$
in~(\ref{eqdefbeta}) and the third inequality follows by the
definition of~$C_\mathrm{mix}$.
Then, setting
$m$ sufficiently large with
respect to $d$, $\eta$, $\varepsilon$ and $\beta/\ell^2$, which makes
$\eta n$ sufficiently large, we obtain that
%
%
\begin{equation}
\mbox{intensity of }\Xi' \geq \biggl(1-\frac{3\varepsilon_k}{4}-
\frac
{\varepsilon_{k+1}}{4} \biggr). \label{eqpastintensity}
\end{equation}
Once the coupling is established,
the probability that all
$ (\frac{\ell_k+2\eta m n k^3\ell_{k-1}}{\ell_{k-1}}
)^d=(mk^3+2\eta mn k^3)^d$
subcubes of scale $k-1$ in $S^\mathrm{ext}_k(\veci)$ have at least
$(1-\varepsilon
_{k})\lambda\ell_{k-1}^d$ nodes of~$\Xi'$ is
%
%
\begin{eqnarray}
\label{eqpastchernoff} \qquad && \mathbf{P} \bigl(D^\mathrm{ext}_k(\veci ,
\tau)=1 \mid F \cap \bigl\{D^\mathrm{base}_k(\veci,\tau)=1
\bigr\} \bigr)
\nonumber
\\
&&\qquad\geq1- \bigl(mk^3+2\eta mnk^3
\bigr)^d
\nonumber
\\[-8pt]
\\[-8pt]
&&\quad\qquad{}\times \exp\biggl(-\frac{1}{2} \biggl(\frac
{\varepsilon_{k}-\varepsilon_{k+1}}{4}
\biggr)^2 \biggl(1-\frac{3\varepsilon
_k}{4}-\frac{\varepsilon_{k+1}}{4} \biggr)
\lambda\ell_{k-1}^d \biggr)
\nonumber
\\
&&\qquad\geq1- \bigl(mk^3+2\eta mnk^3
\bigr)^d \exp \biggl(-\frac{1}{2} \biggl(\frac
{\varepsilon^2}{16(k+1)^4}
\biggr) \biggl(1- \frac{15\varepsilon}{16} \biggr) \lambda\ell_{k-1}^d
\biggr),
\nonumber
\end{eqnarray}
where we used the fact that $\varepsilon_k$ is decreasing in $k$ to
infer that
$1-\frac{3\varepsilon_k}{4}-\frac{\varepsilon_{k+1}}{4} \geq
1-\frac
{3\varepsilon}{4}-\frac{\varepsilon_{2}}{4}=1-\frac{15\varepsilon}{16}$.
Now, for large $k$, the result follows since
$\ell_{k-1}=m^{k-2}((k-1)!)^3\ell$ and, for small $k\geq2$, the
result follows since $\varepsilon^2\lambda\ell^d \geq\alpha$ is
large enough.

For part~(i), a similar argument works. Since in this case we want to
bound the unconditioned probability, at time $\tau\beta_k$ the nodes in
$S^\mathrm{ext}_k(\veci)$ consists of a Poisson point process with intensity
$\lambda$. So, using the derivation in~(\ref{eqpastintensity1})
and~(\ref{eqpastintensity}) with \mbox{$\varepsilon_{k+1}=0$},
the nodes of this Poisson point process for which the displacement throughout
$[\tau\beta_k,(\tau+2)\beta_k]$ is in $Q_{\eta m nk^3\ell_{k-1}}$ is
also a Poisson point process with intensity
at least $1-\frac{\varepsilon_k}{2}$. Then, using a derivation similar
to~(\ref{eqpastchernoff}) we have
\begin{eqnarray*}
\mathbf{P} \bigl(D^\mathrm{ext}_k(\veci,\tau)=1 \bigr) &
\geq& 1- \bigl(m k^3+2\eta mn k^3 \bigr)^d
\exp \biggl(-\frac{1}{2}\frac
{\varepsilon
_k^2}{4}(1- \varepsilon_k/2)
\lambda\ell_{k-1}^d \biggr)
\\
&\geq& 1- \bigl(m k^3+2\eta mn k^3 \bigr)^d
\exp \biggl(-\frac{1}{128}\varepsilon^2(1-\varepsilon/2)\lambda
\ell_{k-1}^d \biggr),
\end{eqnarray*}
where we used the fact that $\varepsilon_k = \varepsilon-\sum_{i=2}^k
\frac
{\varepsilon}{i^2}\geq\varepsilon- \frac{\varepsilon
}{4}-\varepsilon\int_2^\infty
\frac{1}{x^2} \,dx=\frac{\varepsilon}{4}$.

Now, for part~(ii) with $k=1$, we again use the Poisson point process
$\Xi$ of intensity at least
\[
1-\frac{\varepsilon_k}{2}-\frac{\varepsilon_{k+1}}{2} = 1- \frac
{7\varepsilon}{8} %
\]
over $S^\mathrm{ext}_1(\veci)$ defined above. We also use the fact
that $E(\veci,\tau)$ is an event restricted to the super cell $\veci
$ (see the definition
of super cells in the beginning of Section~\ref{secfp}) and
$S^\mathrm{ext}
_1(\veci)$ contains the super cell $\veci$
(cf. Remark~\ref{remextcontainssuper}).
Recall that, for the event $E(\veci,\tau)$, we only consider the
nodes of~$\Xi$ whose displacement from time $\tau\beta$ to $(\tau
+\eta
)\beta$ is inside $Q_{w\ell}$.
Let the event that this happens for a given node of~$\Xi$ be denoted
by $F_{\eta\beta}(w\ell)$. Then we apply
Lemma~\ref
{lemeventF} with $\Delta=\eta\beta$ and $z=w\ell$ to obtain
\[
\mathbf{P} \bigl(F_{\eta\beta}(w\ell) \bigr) \geq1-d\exp \biggl(-
\frac{(w\ell)^2}{18\eta\beta} \biggr). %
\]
Using the fact that $w^2\ell\geq18\eta\beta\log(8d\varepsilon
^{-1})$, we have
$\mathbf{P} (F_{\eta\beta}(w\ell) ) \geq
1-\frac{\varepsilon}{8}$.
Therefore, by thinning, we have that the nodes of~$\Xi$ for which
$F_{\eta\beta}(w\ell)$ hold consist of a
Poisson point process with intensity at least $ (1-\frac
{7\varepsilon
}{8} ) (1-\frac{\varepsilon}{8} )\geq1-\varepsilon$.
Since $E(\veci,\tau)$ is increasing, we have that
\[
\mathbf{P} \bigl(E(\veci,\tau)=0 \mid F\cap \bigl\{D^\mathrm
{base}_k( \veci,\tau)=1 \bigr\} \bigr) \leq1-\nu_E
\bigl((1-\varepsilon) \lambda,Q_{w\ell} \bigr) \leq e^{-\alpha},
\]
which establishes the lemma for $k=1$.
\end{pf}

The lemma below gives an upper bound to the probability that a support
connected path is bad.
Here, we use that all cells in a supported connected path are mutually
well separated so that we can apply Lemma~\ref{lempast}.
Henceforth, we consider paths in $\Omega^\mathrm{sup}_{\kappa-1}$
only since
the cells of scale $\kappa$ will be handled in a different way latter;
these cells are
just so large that a much simpler bound can be applied.
%

\begin{lemma}\label{lemscpath}
Assume the conditions in Lemma~\ref{lempast} are satisfied and
let $P\in\Omega^\mathrm{sup}_{\kappa-1}$ be the path $((k_1,\veci
_1,\tau
_1),(k_2,\veci_2,\tau_2),\ldots,(k_z,\veci_z,\tau_z))$.
Then, with $\psi$ defined as in~(\ref{eqpsi}), we have
\[
\mathbf{P} \Biggl(\bigcap_{j=1}^z \bigl
\{A_{k_j}(\veci_j,\tau_j)=0 \bigr\} \Biggr)
\leq\exp \Biggl(-c_2\sum_{j=1}^z
\psi_{k_j} \Biggr). %
\]
\end{lemma}

\begin{pf}
We derive the probability that all cells of~$P$ are bad.
Consider the following order for the cells of~$P$. First, take an
arbitrary order of~$\mathbb{Z}^d$.
Then we say
that $(k_j,\veci_j,\tau_j)$ precedes $(k_{j'},\veci_{j'},\tau_{j'})$
in the order if
$\tau_j \beta_{k_j} < \tau_{j'} \beta_{k_{j'}}$ or if both $\tau_j
\beta_{k_j}=\tau_{j'}\beta_{k_{j'}}$
and $\veci_j$ precedes $\veci_{j'}$ in the order of~$\mathbb{Z}^d$.
Then, for any $j$, we let $J_j$ be a subset of~$\{1,2,\ldots,z\}$
containing all $j'$ for which $(k_{j'},\veci_{j'},\tau_{j'})$ precedes
$(k_j,\veci_j,\tau_j)$ in the order. Using this order, we write
\begin{eqnarray*}
&& \mathbf{P} \Biggl(\bigcap_{j=1}^z \bigl
\{A_{k_j}(\veci_j,\tau_j)=0 \bigr\} \Biggr)
\\
&&\qquad =
\prod_{j=1}^z \mathbf{P}
\biggl(A_{k_j}( \veci_j,\tau_j)=0 \Big| \bigcap
_{j'\in J_j} \bigl\{A_{k_{j'}}(\veci_{j'},
\tau_{j'})=0 \bigr\} \biggr). %
\end{eqnarray*}
Note that, for each $j'\in J_j$, we have that $(k_{j},\veci_j,\tau_j)$
and $(k_{j'},\veci_{j'},\tau_{j'})$ are well separated.
Using the definition of well separated cells~(\ref{eqwellseparated}),
we have that
$R^\mathrm{inf}_{k_{j'}}(\veci_{j'},\tau_{j'}) \nsubseteq R^{\mathrm
{sup}}_{k_j}(\veci
_j,\tau_j)$ and
$R^\mathrm{inf}_{k_{j}}(\veci_{j},\tau_{j}) \nsubseteq R^{\mathrm
{sup}}_{k_{j'}}(\veci
_{j'},\tau_{j'})$.
Hence, by Lemma~\ref{lemsupport}, we obtain $R^\mathrm
{inf}_{k_{j'}}(\veci
_{j'},\tau_{j'}) \cap R^\mathrm{inf}_{k_{j}}(\veci_{j},\tau
_{j})=\varnothing$.
By the ordering of the cells described above, we also have $\tau
_j\beta
_{k_j} \geq\tau_{j'}\beta_{k_{j'}}$, which gives that\vspace*{1pt}
the event $\bigcap_{j'\in J_j} \{A_{k_{j'}}(\veci_{j'},\tau_{j'})=0\}$
is measurable with respect to
$\mathcal{F}_{k_j}(\veci_j,\tau_j)$.
Then we apply Lemma~\ref{lempast} to obtain a positive constant $c_2$
such that
\[
\mathbf{P} \Biggl(\bigcap_{j=1}^z \bigl
\{A_{k_j}(\veci_j,\tau_j)=0 \bigr\} \Biggr)
\leq\exp \Biggl(-c_2\sum_{j=1}^z
\psi_{k_j} \Biggr). %
\]\upqed
\end{pf}

At the end, we will take the union bound over all support connected cells.
For this, we need to obtain an upper bound for the number of support
connected path, which is given in the following lemma.
%

\begin{lemma}\label{lemsccount}
Let $z$ be a positive integer and $k_1,k_2,\ldots,k_z\geq1$ be fixed.
Then, if $\alpha$ is sufficiently large, the total number of support
connected paths containing $z$ cells whose scales are $k_1,k_2,\ldots
,k_z$ is
at most $\exp(\frac{c_2}{2}\sum_{j=1}^{z}\psi_{k_j} )$,
where $c_2$ is the same constant of Lemma~\ref{lemscpath} and $\psi$
is defined in~(\ref{eqpsi}).
\end{lemma}

\begin{pf}
For any $j,j'\geq1$, define
\begin{eqnarray*}
\phi_{j,j'} &=& \max_{(\veci_1,\tau_1)\in\mathcal{R}_j} \bigl\vert \bigl\{ (
\veci_2,\tau_2) \dvtx(j,\veci_1,
\tau_1)\mbox{ is support adjacent to and well}
\\
&&\hspace*{151pt}
\mbox{separated from }\bigl(j',\veci_2,\tau_2 \bigr) \bigr\}
\bigr\vert. %
\end{eqnarray*}
Hence, given a cell at scale $j$, $\phi_{j,j'}$ is an upper bound for
the number of cells of scale $j'$
that are support adjacent to and well separated from
the cell of scale $j$.
Also, let $\chi_j$ be the number of cells of scale $j$ whose extended
support contains $R_1(\origin,0)$.
With this notation, we obtain that
%
%
\begin{equation}\label{eqdecomposition}
\begin{tabular}{p{280pt}}
the number of support connected paths with $z$ cells of scales
$k_1,k_2,\ldots,k_z \leq
\chi_{k_1}\prod_{j=2}^z
\phi_{k_{j-1},k_j}$.
\end{tabular}
\end{equation}
First, we derive a bound for $\chi_j$. Note that, at scale $j$, the
number of cells that have the same extended support is
$ (\frac{\ell_{j+1}}{\ell_j} )^d\frac{\beta
_{j+1}}{\beta
_j}=m^{d+2}j^2(j+1)^{3d+4}$.
Furthermore, the extended support of a cell of scale $j$ contains
exactly $27\cdot(3(2m+1))^d$ different cells of scale $j+1$. Thus,
the number of different extended supports for a cell of scale $j$ that
contains $R_1(\origin,0)$ can be upper bounded by
\[
\chi_j \leq27 \cdot3^d m^{d+2}(2m+1)^dj^2(j+1)^{3d+4}.
\]
In order to derive a bound for $\phi_{j,j'}$, fix a cell $(j,\veci
_1,\tau_1)$ of scale $j$.
Now, a cell of scale $j'$ can only be support adjacent to $(j,\veci
_1,\tau_1)$ if it is inside the region
\[
\bigcup_{x\in R^{\mathrm{2sup}}_{j}(\veci_1,\tau_1)} \bigl(x+ \bigl[-(3m+2)\ell
_{j'+1},(3m+2)\ell_{j'+1} \bigr]^d\times[-14
\beta_{j'+1},14\beta_{j'+1}] \bigr). %
\]
Therefore, if $j\geq j'$, we have
\begin{eqnarray*}
\phi_{j,j'} &\leq& \biggl( \frac{3(2m+1)\ell_{j+1} +2(3m+2)\ell
_{j'+1}}{\ell
_{j'}} \biggr)^d
\biggl(\frac{27 \beta_{j+1}+28\beta_{j'+1}}{\beta_{j'}} \biggr)
\\
&\leq& \Biggl(6 \biggl(m+\frac{1}{2} \biggr)m^{j-j'+1}\sum
_{i=j'+1}^{j+1}i^3+6 \biggl(m+
\frac{2}{3} \biggr)m \bigl(j'+1 \bigr)^3
\Biggr)^d
\\
&&{} \times \Biggl(27m^{2(j-j'+1)}\sum_{i=j'}^ji^2(i+1)^4
+ 28 m^2 {j'}^2 \bigl(j'+1
\bigr)^4 \Biggr).
\end{eqnarray*}
Using that for any $a\in(0,1)$ and any $x\geq1$ it holds that
$x+a\leq
2x$, we have
\begin{eqnarray*}
\phi_{j,j'} &\leq& \Biggl(12m^{j-j'+2}\sum
_{i=j'+1}^{j+1}i^3+12m^2
\bigl(j'+1 \bigr)^3 \Biggr)^d
\\
&&{}\times
\Biggl(27m^{2(j-j'+1)}\sum_{i=j'}^ji^62^4
+ 28 m^2 {j'}^62^4 \Biggr)
\\
&\leq& c_3 m^{(j-j'+2)d}j^{4d} m^{2(j-j'+1)}j^7
\leq c_3 m^{(d+2)(j-j'+2)}j^{4d+7},
\end{eqnarray*}
where in the derivation above we used the definition of~$\ell_j$ and
$\beta_j$ in~(\ref{eqell}) and~(\ref{eqdefbeta}), respectively,
and $c_3$ is a universal positive constant. We then set $c_4$ in such
a way that, for the case $j<j'$, we obtain
\begin{eqnarray*}
\phi_{j,j'} &\leq& \bigl(1+2(3m+2)m(j+1)^3
\bigr)^d \bigl(1+28 m^2 j^2(j+1)^4
\bigr)
\\
&\leq& \bigl(1+12m^2j^32^3
\bigr)^d \bigl(1+28 m^2 j^62^4
\bigr)
\\
&\leq& c_4 m^{2d+2}j^{3d+6}.
\end{eqnarray*}
Now, for any $j,j',j''$ such that $j\geq j'$ and $j \geq j''$, we have that
\[
\phi_{j,j'}\phi_{j'',j} \leq c_3c_4
m^{(d+2)(j-j'+2)+2d+2}j^{4d+7+3d+6} \leq\exp \biggl(\frac{c_2}{2}
\psi_j \biggr), %
\]
where the last inequality holds for all $j\geq1$ since $\alpha$ is
sufficiently large.
Similarly, by having $\alpha$ sufficiently large, we can guarantee
that, for any $j,j'$ such that $j\geq j'$, we have
\[
\chi_{j}\phi_{j',j} \leq c_5 m^{4d+4}
j^{6d+12} \leq\exp \biggl(\frac{c_2}{2} \psi_j \biggr)
\]
for some constant $c_5$.
Then, if we
consider each term $\phi_{k_{j-1},k_j}$ of~(\ref{eqdecomposition})
and use the bounds above for $\phi$ and $\chi$, we establish the lemma.
\end{pf}

For any support connected path $P=((k_1,\veci_1,\tau_1),(k_2,\veci
_2,\tau_2),\ldots,(k_z,\veci_z,\tau_z))$ in $\Omega^\mathrm
{sup}_{\kappa-1}$,
we define the \emph{weight} of~$P$ as
$\sum_{j=1}^z\psi_{k_j}$. When we take the union bound over all support
connected paths later, we will group the paths by their weight.
The lemma below shows that the paths in $\Omega^\mathrm{sup}_{\kappa-1}$,
which are the paths we need to consider, have
a large enough weight.
%

\begin{lemma}\label{lemscweight}
Let $P=((k_1,\veci_1,\tau_1),(k_2,\veci_2,\tau_2),\ldots
,(k_z,\veci
_z,\tau_z))$ be a path in $\Omega^\mathrm{sup}_{\kappa-1}$.
If $\alpha$ is sufficiently large and $\kappa=O (\log t )$,
then there exist a positive constant $c=c(d)$ and a value $C$
independent of~$t$
such that
\[
\sum_{j=1}^z \psi_{k_j} \geq
\cases{ \displaystyle C \frac{\sqrt{t}}{(\log t)^{c}}, &\quad for $d=1$,
\cr
\noalign{
\vspace*{3pt}} \displaystyle C \frac{t}{(\log t)^{c}}, &\quad for $d=2$,
\cr
\noalign{
\vspace*{3pt}} \displaystyle C t, &\quad for $d\geq3$.} %
\]
\end{lemma}

\begin{pf}
Let $\Delta^\mathrm{2sup}_k$ denote the diameter of the extended
support of
a cell of scale $k$. Then
we have
\begin{eqnarray*}
\Delta^\mathrm{2sup}_k &\leq& 3(2m+1)\ell_{k+1}
\sqrt{d}+27\beta_{k+1}
\\
&=& 3(2m+1)m(k+1)^3\ell_{k}
\sqrt{d}+27C_\mathrm{mix}\frac{\ell
_{k}^2(k+1)^4}{\varepsilon^2}. %
\end{eqnarray*}
Using the definition of~$C_\mathrm{mix}$ from~(\ref{eqdefcmix}), we
obtain a
constant $c_3$ that may depend on the ratio $\beta/\ell^2$ such that
\begin{eqnarray*}
\Delta^\mathrm{2sup}_k &\leq& 3(2m+1)m(k+1)^3
\ell_k\sqrt{d} + 27 \frac{\beta}{\ell^2} m^2\ell
_k^2 (k+1)^4
\\
&\leq& c_3
m^2 (k+1)^4 \ell_k^2. %
\end{eqnarray*}
Then, for $k\geq2$, we have
%
%
\begin{equation}
\qquad\psi_k = \cases{ \displaystyle\frac{\varepsilon^2\lambda\ell
_{k-1}}{(k+1)^4} =
\frac{\varepsilon^2\lambda}{(k+1)^4} \biggl(\frac{\ell
_k}{mk^3} \biggr)
 \geq\frac{\varepsilon^2\lambda}{m(k+1)^7}
\biggl(\frac{\sqrt
{c_3m^2(k+1)^4\ell_k^2}}{\sqrt{c_3m^2(k+1)^4}} \biggr)
\cr
\noalign{\vspace*{3pt}}
\displaystyle\hspace*{41pt}\geq\frac
{\varepsilon^2\lambda}{\sqrt{c_3}m^2(k+1)^9} \sqrt{
\Delta^\mathrm{2sup}_k},\hspace*{49pt}\mbox{for $d=1$,}
\cr\noalign{\vspace*{3pt}} \displaystyle\frac{\varepsilon^2\lambda\ell
_{k-1}^{d-2}}{(k+1)^4} \biggl(\frac{\ell_k}{mk^3}
\biggr)^2
\geq\frac{\varepsilon^2\lambda\ell
_{k-1}^{d-2}}{m^2(k+1)^{10}} \biggl(\frac{c_3m^2(k+1)^4\ell
_k^2}{c_3m^2(k+1)^4} \biggr)
\cr
\noalign{\vspace*{3pt}}
\displaystyle\hspace*{81pt}
\geq\frac{\varepsilon^2 \lambda\ell_{k-1}^{d-2}}{c_3m^4(k+1)^{14}} \Delta^\mathrm{2sup}_k,\qquad\mbox{for $d\geq2$.}} \label{eqdiampsi}
\end{equation}
Now since $\kappa=O(\log t)$, there exists a constant $c_4$ for
which $(k+1)^{a} \leq c_4 (\log t)^a$ for all $k\leq\kappa$ and any
$a\geq1$.
We use this fact for dimensions one and two. For three and higher
dimensions, we simply use the fact that $c_4$ can be set large enough
in order to satisfy also
$\frac{\ell_{k-1}^{d-2}}{(k+1)^{14}}\geq\frac{\ell^{d-2}}{mc_4}$ for
all $k\geq1$.
Plugging this into~(\ref{eqdiampsi}), we obtain
%
%
\begin{equation}
\psi_k \geq\cases{ \displaystyle\frac{\varepsilon^2\lambda
}{\sqrt{c_3}c_4m^2}
\frac
{\sqrt
{\Delta^\mathrm{2sup}_k}}{ (\log t )^{9}}, &\quad for $d=1$,
\cr
\noalign{\vspace*{3pt}} \displaystyle
\frac{\varepsilon^2\lambda}{c_3c_4m^4} \frac{\Delta
^\mathrm{2sup}
_k}{ (\log t )^{14}}, &\quad for $d=2$,
\cr
\noalign{
\vspace*{3pt}} \displaystyle\frac{\varepsilon^2\lambda\ell
^{d-2}}{c_3c_4 m^5} \Delta^\mathrm{2sup}_k,
&\quad for $d\geq3$.} \label{eqdiampsi2}
\end{equation}
For $k=1$, we write $\psi_1 \geq c \sqrt{\Delta^\mathrm{2sup}_1}$
for $d=1$
and $\psi_1 \geq c \Delta^\mathrm{2sup}_1$ for $d\geq2$,
where $c$ is some positive value that may\vspace*{1pt} depend on
$\varepsilon,m,\lambda,\ell$ and $\nu_E$.
Then, if a support connected path is such that
$\sum_{j=1}^z \Delta^\mathrm{2sup}_{k_j} < t$, we have that all the
cells of
the path are contained in $\mathcal{R}^t_1$.
Therefore, for $P\in\Omega^\mathrm{sup}_{\kappa-1}$ we have $\sum_{j=1}^z
\Delta^\mathrm{2sup}_{k_j} \geq t$.
With~(\ref{eqdiampsi2}), this implies that there exists a positive
$C$ independent of~$t$ but depending on everything else such that
\[
\sum_{j=1}^z \psi_{k_j} \geq
\cases{ \displaystyle C \frac{\sqrt{\sum_{j=1}^z \Delta^\mathrm
{2sup}_{k_j}}}{
(\log t )^{9}} \geq C \frac{\sqrt{t}}{ (\log t )^{9}}, &\quad for
$d=1$,
\cr
\noalign{\vspace*{3pt}} \displaystyle C \frac{\sum_{j=1}^z \Delta
^\mathrm
{2sup}_{k_j}}{ (\log
t )^{14}} \geq C
\frac{t}{ (\log t )^{14}}, &\quad for $d=2$,
\cr
\noalign{\vspace*{3pt}} \displaystyle C
\sum_{j=1}^z \Delta^\mathrm{2sup}_{k_j}
\geq Ct, &\quad for $d\geq3$.} %
\]\upqed
\end{pf}

For $j\geq2$, we will work with variables $\tilde\psi_j$ for which
$\tilde\psi_j$ can be written as $b_j \tilde\psi_2$ for some positive
\emph{integer} $b_j$.
Set $\tilde\psi_2 = \psi_2=3^{-4}\varepsilon^2\lambda\ell^d$, and for
$j\geq3$,
define
\[
\tilde\psi_j = 2 \tilde\psi_2 m^{(j-2)d}
\bigl((j-1)! \bigr)^{3d-3} \bigl((j-2)! \bigr)^2(j-3)!.
\]
We only introduce the $\tilde\psi$ in order to simplify a combinatorial
argument later for the counting of support connected paths (we will
need to do this in
order to extend Lemma~\ref{lemsccount} to the case where the scales
are not fixed).
The following lemma establishes that $\psi_j$ and $\tilde\psi_j$
differ only by constant factors.

%
\begin{lemma}\label{lemtildepsi}
For all $j\geq2$, it holds that $\tilde\psi_j \leq\psi_j\leq41
\tilde\psi_j$.
\end{lemma}

\begin{pf}
For $j\geq3$, we write
\begin{eqnarray*}
\psi_{j} &=& \frac{\varepsilon^2\lambda\ell^d
m^{(j-2)d}((j-1)!)^{3d}}{(j+1)^4} = 3^4 \tilde
\psi_2 \frac{m^{(j-2)d}((j-1)!)^{3d}}{(j+1)^4}
\\
&=& 3^4 \tilde\psi_2 m^{(j-2)d} \bigl((j-1)!
\bigr)^{3d-3} \bigl((j-2)! \bigr)^2(j-3)! \biggl(
\frac{(j-1)^3(j-2)}{(j+1)^4} \biggr).
\end{eqnarray*}
This implies that $\psi_j \leq\frac{3^4}{2} \tilde\psi_j\leq41
\tilde\psi_j$.
The other direction follows from the fact that $\frac
{(j-1)^3(j-2)}{(j+1)^4}\geq1/32$ for all $j\geq3$.
\end{pf}

\subsection{Completing the proof of Theorem~\texorpdfstring{\protect\ref
{thmfp}}{3.1}}\label{secproof}
We will need the following technical lemma.

%
\begin{lemma}\label{lemtech}
Let $x,y\in\mathbb{Z}_+$. Then, for any $c_1,c_2>1$, we have
\[
\pmatrix{x+y\cr x}e^{-(c_1x+c_2y)} \leq e^{-(c_1-1)x+(c_2-1)y}. %
\]
\end{lemma}

\begin{pf}
Since ${x+y\choose x}={x+y\choose y}$, we can assume that $x\geq y$.
Then we use the inequality ${x+y\choose x} \leq(\frac
{(x+y)e}{x} )^x$ to obtain
\[
\pmatrix{x+y\cr x}e^{-(c_1x+c_2y)} \leq \biggl(1+\frac{y}{x}
\biggr)^xe^{-(c_1-1)x-c_2y} \leq e^{-(c_1-1)x-(c_2-1)y}. %
\]\upqed
\end{pf}

\begin{pf*}{Proof of Theorem~\ref{thmfp}}
First, for any $k$, note that the number of cells in $\mathcal{R}_k$ satisfies
%
%
\begin{equation}
\vert\mathcal{R}_k\vert\leq \biggl(2 \biggl\lceil
\frac{t}{\ell_k} \biggr\rceil \biggr)^d \biggl\lceil1+
\frac{t}{\beta_k} \biggr\rceil. \label{eqsizer}
\end{equation}
Also, using Lemma~\ref{lemkkp}, we have that
\begin{eqnarray*}
&& \mathbf{P} \bigl(K(\veci,\tau) \nsubseteq\mathcal{R}^t_1
\bigr) \leq\mathbf{P} \bigl(K'(\veci,\tau) \nsubseteq
\mathcal{R}^t_1 \bigr)
\\
&&\qquad = \mathbf{P} (\exists P\in\Omega
\mbox{ s.t. all cells of~$P$ have bad ancestry} ). %
\end{eqnarray*}
Then, with this and Lemma~\ref{lemsupportpath}, we obtain
\[
\mathbf{P} \bigl(K(\veci,\tau) \nsubseteq\mathcal{R}^t_1
\bigr) \leq\mathbf{P} \bigl( \exists P\in\Omega^\mathrm
{sup}_\kappa\mbox{ s.t. all cells of~$P$ are bad} \bigr).
\]
Note that the random variable $A_\kappa$ is defined differently for
scale $\kappa$.
We handle the cells of scale $\kappa$ by showing that none of these
cells are bad with high probability.
It follows by Lemma~\ref{lempast}(i),~(\ref{eqsizer}) and the union
bound over all cells in $\mathcal{R}_\kappa$ that
\[
\mathbf{P} \bigl(A_\kappa \bigl(\veci',
\tau' \bigr)=1\mbox{ for all } (\veci,\tau) \in
\mathcal{R}_\kappa \bigr) \geq1 - \vert\mathcal{R}_\kappa \vert
\exp(-c\psi_\kappa) \geq1 - \exp(-c_1t ), %
\]
for some positive constant $c_1$, where the last step follows by
setting $\kappa$ to be the smallest integer such that $\psi_\kappa
\geq t$,
which gives that $\kappa= \Theta(\frac{\log t}{\log\log
t} )$.
Let $H$ be the event that $A_\kappa(\veci,\tau)=1$ for all $(\veci
,\tau
)\in\mathcal{R}_\kappa$.
Then we have that
\begin{eqnarray*}
&&\mathbf{P} \bigl(\exists P\in\Omega^\mathrm{sup}_\kappa
\mbox{ s.t. all cells of~$P$ are bad} \bigr)
\\
&&\qquad\leq\mathbf{P} \bigl(H \cap \bigl\{\exists P\in\Omega
^\mathrm{sup}_\kappa\mbox{ s.t. all cells of~$P$ are bad}
\bigr\} \bigr) + \mathbf{P} \bigl(H^\compl \bigr)
\\
&&\qquad\leq\mathbf{P} \bigl(\exists P\in\Omega^\mathrm
{sup}_{\kappa-1}\mbox{ s.t. all cells of~$P$ are bad} \bigr) +
e^{-c_1t}.
\end{eqnarray*}

In order to get a bound for the term above,
we first fix a support connected path
%
%
\begin{equation}
P= \bigl((k_1,\veci_1,\tau_1),(k_2,
\veci_2,\tau_2),\ldots,(k_z,\veci
_z,\tau_z) \bigr), \label{eqpathp}
\end{equation}
and use Lemma~\ref{lemscpath} (with the fact that the cells in $P$ are
mutually well separated) to get
\[
\mathbf{P} \Biggl(\bigcap_{j=1}^z \bigl
\{A_{k_j}(\veci_j,\tau_j)=0 \bigr\} \Biggr)
\leq\exp \Biggl(-c_2\sum_{j=1}^z
\psi_{k_j} \Biggr). %
\]
Now, taking the union bound over all support connected paths with $z$
cells of scale $k_1,k_2,\ldots,k_z$, and using Lemma~\ref{lemsccount},
we obtain that
\begin{eqnarray*}
&& \mathbf{P} \bigl(\exists P\in\Omega^\mathrm{sup}_{\kappa-1}
\mbox{ s.t. $P$ has $z$ bad cells of scales $k_1,k_2,
\ldots,k_z$} \bigr)
\\
&&\qquad \leq\exp \Biggl(-\frac{c_2}{2}\sum
_{j=1}^z \psi_{k_j} \Biggr). %
\end{eqnarray*}
Note that the upper bound above depends on $z$ and $k_1,k_2,\ldots,k_z$
only through $\sum_{j=1}^z \psi_{k_j}$, which we call the weight of
the path.
We will group the paths by their weight.
Let $W$ be the set of weights for which there exists at least one path
in $\Omega^\mathrm{sup}_{\kappa-1}$ with that weight. Then
%
%
\begin{equation}
\mathbf{P} \bigl(\exists P\in\Omega^\mathrm{sup}_{\kappa-1}
\mbox{ s.t. all cells of~$P$ are bad} \bigr) \leq\sum_{w \in W}
\exp \biggl(-\frac{c_2}{2}w \biggr) M(w), \label{equseless}
\end{equation}
where $M(w)$ is the number of possible ways to choose $z$ and
$k_1,k_2,\ldots,k_z$ such that $\sum_{j=1}^z \psi_{k_j}=w$.

In order to get an upper bound for $M(w)$, we will use the $\tilde\psi$.
Consider the path $P$ in~(\ref{eqpathp}), let $w=\sum_{j=1}^z \psi
_{k_j}$ and take $w_1 = \psi_1 \vert \{j \dvtx k_j=1\}\vert $.
Let $w_2 = w-w_1$, so $w_1$ is the weight given by cells of scale $1$
and $w_2$ is
the weight given by the other cells of the path.
Note that, by Lemma~\ref{lemtildepsi},
$w_2 = \sum_{j \dvtx k_j \geq2}\psi_{k_j} \geq\sum_{j \dvtx k_j
\geq2}\tilde\psi_{k_j}=h_2 \psi_2$ for some nonnegative integer
$h_2$. Similarly,
we have
$w_2 \leq41 h_2\psi_2$ and $w_1 = h_1 \psi_1$ for some nonnegative
integer $h_1$. Then we have
%
%
\begin{equation}
h_1\psi_1 + h_2\psi_2 \geq
\frac{w}{41}. \label{eqwlb}
\end{equation}

Let $w_0$ be the lower bound on the weight of the path given by
Lemma~\ref{lemscweight}, so for all $w\in W$ we have $w\geq w_0$.
Since either $w_1$ or $w_2$ must be larger than $w_0/2$,
we have that either $h_1 \geq\lceil\frac{w_0}{2\psi_1}\rceil$ or $h_2
\geq\lceil\frac{w_0}{2\cdot41\psi_2}\rceil$.
Let $M(h_1,h_2)$ be the number of ways to choose $z$ and
$k_1,k_2,\ldots,k_z$ such that there are
$h_1$ values $j$ with $k_j=1$ and $\sum_{j\dvtx k_j \geq2}\tilde
\psi
_{k_j}=h_2 \psi_2$.
Note that, for any such choice, we have $w = \sum_{j=1}^z \psi_{k_j}
\geq h_1\psi_1+h_2\psi_2$.
Thus, the sum in the right-hand side of~(\ref{equseless}) can be
bounded above by
\begin{eqnarray*}
&& \sum_{h_1=\lceil {w_0}/{(2\psi_1)}\rceil}^\infty\sum
_{h_2=0}^\infty\exp \biggl(-\frac{c_2}{2}
(h_1\psi_1+h_2\psi_2) \biggr)
M(h_1,h_2)
\\
&&\qquad{}+ \sum_{h_1=0}^\infty\sum
_{h_2=\lceil {w_0}/{(82\psi_2)}\rceil
}^\infty\exp \biggl(-\frac{c_2}{2}
(h_1\psi_1+h_2\psi_2) \biggr)
M(h_1,h_2).
\end{eqnarray*}

In order to bound $M(h_1,h_2)$, we will consider the following
pictorial way to define the values of~$z$ and $k_1,k_2,\ldots,k_z$.
Suppose we have $h_1$ blocks of size $\psi_1$ and $h_2$ blocks of size
$\psi_2$.
Consider an ordering of the blocks, but such that permuting blocks of
the same size does not change the order.
Then, for each block
of size~$\psi_2$, we color it either black or white, while
blocks of size $\psi_1$ are not colored.
Now, for each choice of~$z$ and $k_1,k_2,\ldots,k_z$ we associate an
order and coloring of the
blocks as follows.
If $k_1=1$, then the first block is of size $\psi_1$. Otherwise, the
first $\tilde\psi_{k_1}/\psi_2$ blocks are of size $\psi_2$ and have
black color.
Then, if $k_2=1$, the next block is of size $\psi_1$, otherwise the next
$\tilde\psi_{k_2}/\psi_2$ blocks are of size $\psi_2$ and have white
color. We proceed in this way until $k_z$, where whenever $k_i\neq1$
we use the color black
if $i$ is odd and the color white if $i$ is even. Though
there are orders and colorings that are not associated to any choice
of~$z$ and $k_1,k_2,\ldots,k_z$,
each such choice of~$z$ and $k_1,k_2,\ldots,k_z$
corresponds to a unique order and coloring of the blocks.
Therefore, the number of ways to order and color the blocks gives an
upper bound for $M(h_1,h_2)$.
Note that there are ${h_1+h_2\choose h_1}$ ways to order the blocks
and $2^{h_2}$ ways to color the size-$\psi_2$ blocks.
Therefore,
\begin{eqnarray*}
&& \mathbf{P} \bigl(\exists P\in\Omega^\mathrm{sup}_{\kappa-1}
\mbox{ s.t. all cells of~$P$ are bad} \bigr)
\\
&&\qquad\leq\sum_{h_1=\lceil{w_0}/{(2\psi_1)}\rceil}^\infty \sum
_{h_2=0}^\infty\exp \biggl(-\frac{c_2}{2}
(h_1\psi_1+h_2\psi_2) \biggr)
\pmatrix{h_1+h_2\cr h_1}2^{h_2}
\\
&&\quad\qquad{}+ \sum_{h_1=0}^\infty\sum
_{h_2=\lceil{w_0}/{(82\psi_2)}\rceil
}^\infty\exp \biggl(-\frac{c_2}{2}
(h_1\psi_1+h_2\psi_2) \biggr)
\pmatrix{h_1+h_2\cr h_1}2^{h_2}
\\
&&\qquad\leq\sum_{h_1=\lceil{w_0}/{(2\psi_1)}\rceil}^\infty \sum
_{h_2=0}^\infty\exp \biggl(-\frac{c_2(h_1\psi_1+h_2\psi
_2)}{3} \biggr)
\\
&&\qquad\quad{}+
\sum_{h_1=0}^\infty\sum
_{h_2=\lceil{w_0}/{(82\psi_2)}\rceil
}^\infty\exp \biggl(-\frac{c_2(h_1\psi_1+h_2\psi_2)}{3} \biggr)
\\
&&\qquad\leq\sum_{w\in W} \exp(-c_4 w )
\leq\exp(-c w_0 ),
\end{eqnarray*}
where in the second inequality we use Lemma~\ref{lemtech} and the fact
that $\alpha$ is sufficiently large to write $\frac{c_2\psi
_1}{2}-1\geq
\frac{c_2\psi_1}{3}$, and
similarly for $\psi_2$. In the third inequality, we used~(\ref
{eqwlb}). Since $w_0$ is the lower bound in Lemma~\ref{lemscweight},
the proof of Theorem~\ref{thmfp} is completed.
\end{pf*}

\section{Detection}\label{secdetection}
In this section, we use Theorem~\ref{thmfp} to prove Theorem~\ref
{thmdetection}.

\begin{pf*}{Proof of Theorem~\ref{thmdetection}}
Recall that we say that the displacement of a node throughout
$[t_0,t_1]$ is in $Q_z$ if the node never leaves
$x+Q_z$ during the whole of~$[t_0,t_1]$, where $x$ is the position of
the node at time $t_0$.
Fix $\eta=1$ and $\varepsilon=1/2$.
Now, we fix the ratio $\beta/\ell^2$ small enough so that the lower
bound for $w$ in Theorem~\ref{thmfp} is at most $1$ and
we can set $w=1$.
We tessellate $\mathbb{R}^d$ into cells of side length $\ell=\frac
{r}{2\sqrt{d}}$.
Let $S$ be a cell of the tessellation and let $v$ be a node of~$\Pi_0$
that is inside $S$ at some time $s$.
Then, if the displacement of~$v$ during the interval $[s,s+\beta]$ is
in $Q_{w\ell}$, we have that
the distance between $v$ and any point of~$S$ at any time in
$[s,s+\beta]$ is at most $\frac{\ell\sqrt{d}}{2}+\ell\sqrt
{d}=\frac
{3\sqrt{d}\ell}{2}\leq r$.
Therefore, for such a node $v$,
the ball of radius $r$ centered at $v$ covers the whole of~$S$
during the entire duration of the interval $[s,s+\beta]$. Hence, if
$S$ contains at least one such node at time $s$ and
the target enters $S$ during $[s,s+\beta]$, then the target is detected.

Now, we apply Theorem~\ref{thmfp}.
For each $(\veci,\tau)\in\mathcal{R}_1$,
define $E(\veci,\tau)$ to be the event that there is at least one node
in the cube
$S_1(\veci)$ at time $\tau\beta$ for which its displacement from time
$\tau\beta$ to $(\tau+1)\beta$ is inside $Q_\ell$.
This event is clearly increasing.
Let $N$ be a Poisson random variable of mean
$\frac{\lambda\ell^d}{2}$.
Then, using the fact that $\varepsilon=1/2$,
we have that
\[
\nu_E(\lambda/2,Q_\ell) \geq\mathbf{P} (N \geq1 ) = 1-
\exp \biggl(-\frac{\lambda
\ell^d}{2} \biggr). %
\]
Clearly, $\log(\frac{1}{1-\nu_E(\lambda/2,Q_\ell)}
)\geq
\frac{\lambda\ell^d}{2}$, which increases with $\lambda$.
Therefore, we can set $\lambda$ large enough so that $\log
(\frac
{1}{1-\nu_E(\lambda/2,Q_\ell)} )$ and
$\varepsilon^2\lambda\ell^d= \lambda\ell^d/4$ are larger than
$\alpha_0$.
With this, we apply Theorem~\ref{thmfp} to obtain that
%
%
\begin{equation}
\mathbf{P} \bigl(K(\origin,0)\subseteq\mathcal{R}^t_1
\bigr) \leq\cases{ \displaystyle\exp \biggl(-C \frac{\sqrt
{t}}{(\log t)^{c}} \biggr), &\quad
for $d=1$,
\cr
\noalign{\vspace*{3pt}} \displaystyle\exp \biggl(-C
\frac{t}{(\log t)^{c}} \biggr), &\quad for $d=2$,
\cr
\noalign{\vspace*{3pt}}
\displaystyle\exp(-Ct ), &\quad for $d\geq3$.} \label{eqfinalstep}
\end{equation}
Note that the target is not detected at time $0$ only if $E(\origin
,0)=0$. Then, in this case, since $K(\origin,0)$ is contained in
$\mathcal{R}^t_1$ and
$\mathcal{R}^t_1$ contains all the cells that are contained in the
space--time region $(-t,t)\times[0,t)$ we have that the target must be
at some time in the
interval $[0,t]$ inside a cell $(\veci,\tau)$ of scale 1 for which
$E(\veci,\tau)=1$ and, therefore, the target is detected.

This shows that the probability that the target is able to evade
detection up to time $t$ is given by~(\ref{eqfinalstep}), which
completes the proof of
Theorem~\ref{thmdetection}.
\end{pf*}

The same proof as above can be used to establish that, for any $\lambda
>0$, there exists a value $\bar r=\bar r(\lambda)$ so that with high
probability the target will
eventually get within distance $\bar r$ from at least one node. We
state this slight generalization below. In Theorem~\ref{thmdetection},
we require $\lambda$ to
be large enough so that $\bar r \leq r$.

%
%
\begin{theorem}
In dimensions $d\geq2$, there exist an explicit constant $c=c(d)$ and
a positive $C$ independent of~$t$ so that the following holds for all
large enough~$t$.
For any $\lambda>0$, there exists $\bar r=\bar r(\lambda)>0$ so that
the probability that there exists a trajectory $g$ for the target so that
for all $s\in[0,t]$ the ball $B(g(s),\bar r)$ contains no node of~$\Pi
_s$ is at most $\exp(-C\frac{t}{(\log t)^c} )$ in $d=2$ and
at most $\exp(-Ct )$ in $d\geq3$.
\end{theorem}

\begin{appendix}\label{app}
\section*{Appendix: Standard large deviation results}
\setcounter{theorem}{0}

We use the following standard Chernoff bounds and large deviation results.

%

\begin{lemma}[(Chernoff bound for Poisson)]\label{lemcbpoisson}
Let $P$ be a Poisson random variable with mean $\lambda$. Then, for any
$0<\varepsilon<1$,
\[
\mathbf{P} \bigl(P \geq(1+\varepsilon) \lambda \bigr) \leq\exp \biggl(-
\frac{\lambda\varepsilon
^2}{2}(1-\varepsilon/3) \biggr) %
\]
and
\[
\mathbf{P} \bigl(P \leq(1-\varepsilon) \lambda \bigr) \leq\exp \biggl(-
\frac{\lambda\varepsilon
^2}{2} \biggr). %
\]
\end{lemma}




%
%
\begin{lemma}[(Gaussian tail bound~\cite{BM}, Theorem~12.9)]\label{lemgaussiantail}
Let $X$ be a normal random variable with mean $0$ and variance $\sigma
^2$. Then, for any $R \geq\sigma$ we have that
$\mathbf{P} (X \geq R ) \leq\frac{\sigma
}{\sqrt{2\pi}R}\exp(-\frac
{R^2}{2\sigma^2} )$.
\end{lemma}
\end{appendix}

\section*{Acknowledgements}
We are grateful to Yuval Peres, Alistair Sinclair, Perla Sousi and Omer
Tamuz for useful discussions.
We are also thankful to Ron Peled for valuable suggestions that
improved the
bounds in~(\ref{eqtaildetection}), and to
Alistair Sinclair for helpful comments on a previous version of the paper.


%

\printaddresses

\begin{thebibliography}{32}

\bibitem{Bal}
%
\begin{bincollection}[auto:STB|2014/08/04|07:23:14]
\bauthor{\bsnm{Balister},~\bfnm{P.}\binits{P.}},
\bauthor{\bsnm{Zheng},~\bfnm{Z.}\binits{Z.}},
\bauthor{\bsnm{Kumar},~\bfnm{S.}\binits{S.}} \AND
\bauthor{\bsnm{Sinha},~\bfnm{P.}\binits{P.}}
(\byear{2009}).
\btitle{Trap coverage: Allowing coverage holes of bounded diameter in
wireless sensor networks}.
In \bbooktitle{Proceedings of the 28th IEEE Conference on Computer
Communications}
\bpages{19--25}.
\bpublisher{IEEE}, \blocation{Piscataway, NJ}.
\end{bincollection}
%
\bptok{imsref}%
\endbibitem

\bibitem{BenjaminiStauffer11}
%
\begin{barticle}[mr]
\bauthor{\bsnm{Benjamini},~\bfnm{Itai}\binits{I.}} \AND
\bauthor{\bsnm{Stauffer},~\bfnm{Alexandre}\binits{A.}}
(\byear{2013}).
\btitle{Perturbing the hexagonal circle packing: A percolation perspective}.
\bjournal{Ann. Inst. Henri Poincar\'e Probab. Stat.}
\bvolume{49}
\bpages{1141--1157}.
\bid{doi={10.1214/12-AIHP524}, issn={0246-0203}, mr={3127917}}
\end{barticle}
%
\bptok{imsref}%
\endbibitem

\bibitem{CPPR10}
%
\begin{barticle}[mr]
\bauthor{\bsnm{Chatterjee},~\bfnm{Sourav}\binits{S.}},
\bauthor{\bsnm{Peled},~\bfnm{Ron}\binits{R.}},
\bauthor{\bsnm{Peres},~\bfnm{Yuval}\binits{Y.}} \AND
\bauthor{\bsnm{Romik},~\bfnm{Dan}\binits{D.}}
(\byear{2010}).
\btitle{Phase transitions in gravitational allocation}.
\bjournal{Geom. Funct. Anal.}
\bvolume{20}
\bpages{870--917}.
\bid{doi={10.1007/s00039-010-0090-7}, issn={1016-443X}, mr={2729280}}
\end{barticle}
%
\bptok{imsref}%
\endbibitem

\bibitem{CC89}
%
\begin{barticle}[mr]
\bauthor{\bsnm{Chayes},~\bfnm{J.~T.}\binits{J.~T.}} \AND
\bauthor{\bsnm{Chayes},~\bfnm{L.}\binits{L.}}
(\byear{1989}).
\btitle{The large-{$N$} limit of the threshold values in
{M}andelbrot's fractal percolation process}.
\bjournal{J. Phys. A}
\bvolume{22}
\bpages{L501--L506}.
\bid{issn={0305-4470}, mr={1003259}}
\end{barticle}
%
\bptok{imsref}%
\endbibitem

\bibitem{CCD88}
%
\begin{barticle}[mr]
\bauthor{\bsnm{Chayes},~\bfnm{J.~T.}\binits{J.~T.}},
\bauthor{\bsnm{Chayes},~\bfnm{L.}\binits{L.}} \AND
\bauthor{\bsnm{Durrett},~\bfnm{R.}\binits{R.}}
(\byear{1988}).
\btitle{Connectivity properties of {M}andelbrot's percolation process}.
\bjournal{Probab. Theory Related Fields}
\bvolume{77}
\bpages{307--324}.
\bid{doi={10.1007/BF00319291}, issn={0178-8051}, mr={0931500}}
\end{barticle}
%
\bptok{imsref}%
\endbibitem

\bibitem{CCGS91}
%
\begin{barticle}[mr]
\bauthor{\bsnm{Chayes},~\bfnm{J.~T.}\binits{J.~T.}},
\bauthor{\bsnm{Chayes},~\bfnm{L.}\binits{L.}},
\bauthor{\bsnm{Grannan},~\bfnm{E.}\binits{E.}} \AND
\bauthor{\bsnm{Swindle},~\bfnm{G.}\binits{G.}}
(\byear{1991}).
\btitle{Phase transitions in {M}andelbrot's percolation process in
three dimensions}.
\bjournal{Probab. Theory Related Fields}
\bvolume{90}
\bpages{291--300}.
\bid{doi={10.1007/BF01193747}, issn={0178-8051}, mr={1133368}}
\end{barticle}
%
\bptok{imsref}%
\endbibitem

\bibitem{clementi}
%
\begin{bincollection}[mr]
\bauthor{\bsnm{Clementi},~\bfnm{Andrea~E.~F.}\binits{A.~E.~F.}},
\bauthor{\bsnm{Pasquale},~\bfnm{Francesco}\binits{F.}} \AND
\bauthor{\bsnm{Silvestri},~\bfnm{Riccardo}\binits{R.}}
(\byear{2009}).
\btitle{M{ANETS}: High mobility can make up for low transmission power}.
In \bbooktitle{Automata, Languages and Programming. {P}art {II}}.
\bpages{387--398}.
\bpublisher{Springer},
\blocation{Berlin}.
\bid{doi={10.1007/978-3-642-02930-1_32}, mr={2544811}}
\end{bincollection}
%
\bptok{imsref}%
\endbibitem

\bibitem{DiazFinal}
%
\begin{barticle}[auto]
\bauthor{\bsnm{D{\'{\i}}az},~\bfnm{Josep}\binits{J.}},
\bauthor{\bsnm{Mitsche},~\bfnm{Dieter}\binits{D.}} \AND
\bauthor{\bsnm{P{\'e}rez-Gim{\'e}nez},~\bfnm{Xavier}\binits{X.}}
(\byear{2009}).
\btitle{Large connectivity for dynamic random geometric graphs}.
\bjournal{IEEE Trans. Mob. Comput.}
\bvolume{8}
\bpages{821--835}.
\end{barticle}
%
\bptok{imsref}%
\endbibitem

\bibitem{Dousse}
%
\begin{bincollection}[auto:STB|2014/08/04|07:23:14]
\bauthor{\bsnm{Dousse},~\bfnm{O.}\binits{O.}},
\bauthor{\bsnm{Tavoularis},~\bfnm{C.}\binits{C.}} \AND
\bauthor{\bsnm{Thiran},~\bfnm{P.}\binits{P.}}
(\byear{2006}).
\btitle{Delay of intrusion detection in wireless sensor networks}.
In \bbooktitle{Proceedings of the 7th ACM International Conference on
Mobile Computing and Networking (MobiCom)}
\bpages{155--165}.
\bpublisher{ACM}, \blocation{New York, NY}.
\end{bincollection}
%
\bptok{imsref}%
\endbibitem

\bibitem{DGRS}
%
\begin{bincollection}[auto:STB|2014/08/04|07:23:14]
\bauthor{\bsnm{Drewitz},~\bfnm{A.}\binits{A.}},
\bauthor{\bsnm{G{\"a}rtner},~\bfnm{J.}\binits{J.}},
\bauthor{\bsnm{Ram{\'i}rez},~\bfnm{A.~F.}\binits{A.~F.}} \AND
\bauthor{\bsnm{Sun},~\bfnm{R.}\binits{R.}}
(\byear{2012}).
\btitle{Survival probability of a random walk among a Poisson system
of moving traps}.
In \bbooktitle{Probability in Complex Physical Systems, Springer
Proceedings in Mathematics}
\bpages{119--158}.
\bpublisher{Springer},
\blocation{Berlin}.
\end{bincollection}
%
\bptok{imsref}%
\endbibitem

\bibitem{elgamal}
%
\begin{bincollection}[auto:STB|2014/08/04|07:23:14]
\bauthor{\bsnm{El Gamal},~\bfnm{A.}\binits{A.}},
\bauthor{\bsnm{Mammen},~\bfnm{J.}\binits{J.}},
\bauthor{\bsnm{Prabhakar},~\bfnm{B.}\binits{B.}} \AND
\bauthor{\bsnm{Shah},~\bfnm{D.}\binits{D.}}
(\byear{2004}).
\btitle{Throughput-delay trade-off in wireless networks}.
In \bbooktitle{Proceedings of the 23rd IEEE Conference on Computer
Communications}
\bpages{464--475}.
\bpublisher{IEEE}, \blocation{Piscataway, NJ}.
\end{bincollection}
%
\bptok{imsref}%
\endbibitem

\bibitem{FG92}
%
\begin{barticle}[mr]
\bauthor{\bsnm{Falconer},~\bfnm{K.~J.}\binits{K.~J.}} \AND
\bauthor{\bsnm{Grimmett},~\bfnm{G.~R.}\binits{G.~R.}}
(\byear{1992}).
\btitle{On the geometry of random {C}antor sets and fractal percolation}.
\bjournal{J. Theoret. Probab.}
\bvolume{5}
\bpages{465--485}.
\bid{doi={10.1007/BF01060430}, issn={0894-9840}, mr={1176432}}
\end{barticle}
%
\bptok{imsref}%
\endbibitem

\bibitem{FSS13}
%
\begin{barticle}[mr]
\bauthor{\bsnm{Friedrich},~\bfnm{Tobias}\binits{T.}},
\bauthor{\bsnm{Sauerwald},~\bfnm{Thomas}\binits{T.}} \AND
\bauthor{\bsnm{Stauffer},~\bfnm{Alexandre}\binits{A.}}
(\byear{2013}).
\btitle{Diameter and broadcast time of random geometric graphs in
arbitrary dimensions}.
\bjournal{Algorithmica}
\bvolume{67}
\bpages{65--88}.
\bid{doi={10.1007/s00453-012-9710-y}, issn={0178-4617}, mr={3072817}}
\end{barticle}
%
\bptok{imsref}%
\endbibitem

\bibitem{Kesidis}
%
\begin{bincollection}[auto:STB|2014/08/04|07:23:14]
\bauthor{\bsnm{Kesidis},~\bfnm{G.}\binits{G.}},
\bauthor{\bsnm{Konstantopoulos},~\bfnm{T.}\binits{T.}} \AND
\bauthor{\bsnm{Phoha},~\bfnm{S.}\binits{S.}}
(\byear{2003}).
\btitle{Surveillance coverage of sensor networks under a random
mobility strategy}.
In \bbooktitle{Proceedings of the 2nd IEEE International Conference on
Sensors}.
\bpublisher{IEEE}, \blocation{Piscataway, NJ}.
\end{bincollection}
%
\bptok{imsref}%
\endbibitem

\bibitem{KS05}
%
\begin{barticle}[mr]
\bauthor{\bsnm{Kesten},~\bfnm{Harry}\binits{H.}} \AND
\bauthor{\bsnm{Sidoravicius},~\bfnm{Vladas}\binits{V.}}
(\byear{2005}).
\btitle{The spread of a rumor or infection in a moving population}.
\bjournal{Ann. Probab.}
\bvolume{33}
\bpages{2402--2462}.
\bid{doi={10.1214/009117905000000413}, issn={0091-1798}, mr={2184100}}
\end{barticle}
%
\bptok{imsref}%
\endbibitem

\bibitem{KS04}
%
\begin{barticle}[mr]
\bauthor{\bsnm{Kesten},~\bfnm{Harry}\binits{H.}} \AND
\bauthor{\bsnm{Sidoravicius},~\bfnm{Vladas}\binits{V.}}
(\byear{2006}).
\btitle{A phase transition in a model for the spread of an infection}.
\bjournal{Illinois J. Math.}
\bvolume{50}
\bpages{547--634}.
\bid{issn={0019-2082}, mr={2247840}}
\end{barticle}
%
\bptok{imsref}%
\endbibitem

\bibitem{Konst}
\begin{barticle}[auto:STB|2014/08/04|07:23:14]
\bauthor{\bsnm{Konstantopoulos},~\bfnm{T.}\binits{T.}}
(\byear{2010}).
\btitle{Response to Prof. Baccelli's lecture on modelling of wireless communication networks by stochastic geometry}.
\bjournal{The Computer Journal}
\bvolume{53}
\bpages{612--614}.
\end{barticle}
\bptok{imsref}%
\endbibitem

\bibitem{LLMSW12}
%
\begin{bincollection}[auto:STB|2014/08/04|07:23:14]
\bauthor{\bsnm{Lam},~\bfnm{H.}\binits{H.}},
\bauthor{\bsnm{Liu},~\bfnm{Z.}\binits{Z.}},
\bauthor{\bsnm{Mitzenmacher},~\bfnm{M.}\binits{M.}},
\bauthor{\bsnm{Sun},~\bfnm{X.}\binits{X.}} \AND
\bauthor{\bsnm{Wang},~\bfnm{Y.}\binits{Y.}}
(\byear{2012}).
\btitle{Information dissemination via random walks in $d$-dimensional space}.
In \bbooktitle{Proceedings of the 23th {ACM-SIAM} Symposium on
Discrete Algorithms (SODA)}.
\bpublisher{SIAM}, \blocation{Philadelphia, PA}.
\end{bincollection}
%
\bptok{imsref}%
\endbibitem

\bibitem{Liu}
%
\begin{bincollection}[auto:STB|2014/08/04|07:23:14]
\bauthor{\bsnm{Liu},~\bfnm{B.}\binits{B.}},
\bauthor{\bsnm{Brass},~\bfnm{P.}\binits{P.}},
\bauthor{\bsnm{Dousse},~\bfnm{O.}\binits{O.}},
\bauthor{\bsnm{Nain},~\bfnm{P.}\binits{P.}} \AND
\bauthor{\bsnm{Towsley},~\bfnm{D.}\binits{D.}}
(\byear{2005}).
\btitle{Mobility improves coverage of sensor networks}.
In \bbooktitle{Proceedings of the 6th ACM International Conference on
Mobile Computing and Networking (MobiCom)}.
\bpublisher{ACM}, \blocation{New York, NY}.
\end{bincollection}
%
\bptok{imsref}%
\endbibitem

\bibitem{MR}
%
\begin{bbook}[mr]
\bauthor{\bsnm{Meester},~\bfnm{Ronald}\binits{R.}} \AND
\bauthor{\bsnm{Roy},~\bfnm{Rahul}\binits{R.}}
(\byear{1996}).
\btitle{Continuum Percolation}.
\bpublisher{Cambridge Univ. Press},
\blocation{Cambridge}.
\bid{doi={10.1017/CBO9780511895357}, mr={1409145}}
\end{bbook}
%
\bptok{imsref}%
\endbibitem

\bibitem{MPV01}
%
\begin{barticle}[mr]
\bauthor{\bsnm{Menshikov},~\bfnm{M.~V.}\binits{M.~V.}},
\bauthor{\bsnm{Popov},~\bfnm{S.~Yu.}\binits{S.~Yu.}} \AND
\bauthor{\bsnm{Vachkovskaia},~\bfnm{M.}\binits{M.}}
(\byear{2001}).
\btitle{On the connectivity properties of the complementary set in
fractal percolation models}.
\bjournal{Probab. Theory Related Fields}
\bvolume{119}
\bpages{176--186}.
\bid{doi={10.1007/PL00008757}, issn={0178-8051}, mr={1818245}}
\end{barticle}
%
\bptok{imsref}%
\endbibitem

\bibitem{MOBC}
%
\begin{barticle}[auto:STB|2014/08/04|07:23:14]
\bauthor{\bsnm{Moreau},~\bfnm{M.}\binits{M.}},
\bauthor{\bsnm{Oshanin},~\bfnm{G.}\binits{G.}},
\bauthor{\bsnm{B{\'e}nichou},~\bfnm{O.}\binits{O.}} \AND
\bauthor{\bsnm{Coppey},~\bfnm{M.}\binits{M.}}
(\byear{2004}).
\btitle{Lattice theory of trapping reactions with mobile species}.
\bjournal{Phys. Rev. E}
\bvolume{69}
\bpages{046101}.
\end{barticle}
%
\bptok{imsref}%
\endbibitem

\bibitem{BM}
%
\begin{bbook}[mr]
\bauthor{\bsnm{M{\"o}rters},~\bfnm{Peter}\binits{P.}} \AND
\bauthor{\bsnm{Peres},~\bfnm{Yuval}\binits{Y.}}
(\byear{2010}).
\btitle{Brownian Motion}.
\bpublisher{Cambridge Univ. Press},
\blocation{Cambridge}.
\bid{doi={10.1017/CBO9780511750489}, mr={2604525}}
\end{bbook}
%
\bptok{imsref}%
\endbibitem

\bibitem{O96}
%
\begin{barticle}[mr]
\bauthor{\bsnm{Orzechowski},~\bfnm{M.~E.}\binits{M.~E.}}
(\byear{1996}).
\btitle{On the phase transition to sheet percolation in random
{C}antor sets}.
\bjournal{J. Stat. Phys.}
\bvolume{82}
\bpages{1081--1098}.
\bid{doi={10.1007/BF02179803}, issn={0022-4715}, mr={1372435}}
\end{barticle}
%
\bptok{imsref}%
\endbibitem

\bibitem{Penrose}
%
\begin{bbook}[mr]
\bauthor{\bsnm{Penrose},~\bfnm{Mathew}\binits{M.}}
(\byear{2003}).
\btitle{Random Geometric Graphs}.
\bpublisher{Oxford Univ. Press},
\blocation{Oxford}.
\bid{doi={10.1093/acprof:oso/9780198506263.001.0001}, mr={1986198}}
\end{bbook}
%
\bptok{imsref}%
\endbibitem

\bibitem{PSSS11}
%
\begin{barticle}[mr]
\bauthor{\bsnm{Peres},~\bfnm{Yuval}\binits{Y.}},
\bauthor{\bsnm{Sinclair},~\bfnm{Alistair}\binits{A.}},
\bauthor{\bsnm{Sousi},~\bfnm{Perla}\binits{P.}} \AND
\bauthor{\bsnm{Stauffer},~\bfnm{Alexandre}\binits{A.}}
(\byear{2013}).
\btitle{Mobile geometric graphs: Detection, coverage and percolation}.
\bjournal{Probab. Theory Related Fields}
\bvolume{156}
\bpages{273--305}.
\bid{doi={10.1007/s00440-012-0428-1}, issn={0178-8051}, mr={3055260}}
\end{barticle}
%
\bptok{imsref}%
\endbibitem

\bibitem{PeresSousi11}
%
\begin{barticle}[mr]
\bauthor{\bsnm{Peres},~\bfnm{Yuval}\binits{Y.}} \AND
\bauthor{\bsnm{Sousi},~\bfnm{Perla}\binits{P.}}
(\byear{2012}).
\btitle{An isoperimetric inequality for the {W}iener sausage}.
\bjournal{Geom. Funct. Anal.}
\bvolume{22}
\bpages{1000--1014}.
\bid{doi={10.1007/s00039-012-0184-5}, issn={1016-443X}, mr={2984124}}
\end{barticle}
%
\bptok{imsref}%
\endbibitem

\bibitem{PPPU11}
%
\begin{bincollection}[auto:STB|2014/08/04|07:23:14]
\bauthor{\bsnm{Pettarin},~\bfnm{A.}\binits{A.}},
\bauthor{\bsnm{Pietracaprina},~\bfnm{A.}\binits{A.}},
\bauthor{\bsnm{Pucci},~\bfnm{G.}\binits{G.}} \AND
\bauthor{\bsnm{Upfal},~\bfnm{E.}\binits{E.}}
(\byear{2011}).
\btitle{Tight bounds on information dissemination in sparse mobile networks}.
In \bbooktitle{30th Annual {ACM SIGACT-SIGOPS} Symposium on Principles
of Distributed Computing}.
\bpublisher{ACM}, \blocation{New York, NY}.
\end{bincollection}
%
\bptok{imsref}%
\endbibitem

\bibitem{SinclairStauffer10}
%
\begin{bmisc}[auto:STB|2014/08/04|07:23:14]
\bauthor{\bsnm{Sinclair},~\bfnm{A.}\binits{A.}} \AND
\bauthor{\bsnm{Stauffer},~\bfnm{A.}\binits{A.}}
(\byear{2010}).
\bhowpublished{Mobile geometric graphs, and detection and
communication problems in mobile wireless networks.
Preprint. Available at \arxivurl{arXiv:1005.1117v1}.}
\end{bmisc}
%
\bptok{imsref}%
\endbibitem

\bibitem{SKM95}
%
\begin{bbook}[auto]
\bauthor{\bsnm{Stoyan},~\bfnm{D.}\binits{D.}},
\bauthor{\bsnm{Kendall},~\bfnm{W.~S.}\binits{W.~S.}} \AND
\bauthor{\bsnm{Mecke},~\bfnm{J.}\binits{J.}}
(\byear{1995}).
\btitle{Stochastic Geometry and Its Applications},
\bedition{2nd} ed.
\bpublisher{Wiley},
\blocation{Chichester}.
\end{bbook}
%
\bptok{imsref}%
\endbibitem

\bibitem{vandenberg}
%
\begin{barticle}[mr]
\bauthor{\bsnm{van~den Berg},~\bfnm{J.}\binits{J.}},
\bauthor{\bsnm{Meester},~\bfnm{Ronald}\binits{R.}} \AND
\bauthor{\bsnm{White},~\bfnm{Damien~G.}\binits{D.~G.}}
(\byear{1997}).
\btitle{Dynamic {B}oolean models}.
\bjournal{Stochastic Process. Appl.}
\bvolume{69}
\bpages{247--257}.
\bid{doi={10.1016/S0304-4149(97)00044-6}, issn={0304-4149}, mr={1472953}}
\end{barticle}
%
\bptok{imsref}%
\endbibitem

\bibitem{YG10}
%
\begin{bincollection}[auto:STB|2014/08/04|07:23:14]
\bauthor{\bsnm{Yanmaz},~\bfnm{E.}\binits{E.}} \AND
\bauthor{\bsnm{Guclu},~\bfnm{H.}\binits{H.}}
(\byear{2010}).
\btitle{Stationary and mobile target detection using mobile wireless
sensor networks}.
In \bbooktitle{Proceedings of the 29th Conference on Computer
Communications (INFOCOM)}.
\bpublisher{IEEE}, \blocation{Piscataway, NJ}.
\end{bincollection}
%
\bptok{imsref}%
\endbibitem

\end{thebibliography}
\end{document}